# On Topics Related to Tensor Products with Abelian von Neumann Algebras

David Gao

**Abstract:** In this note, we developed several results concerning abelian von Neumann algebras, their spectrums, and their tensor products with other von Neumann algebras. In particular, we developed a theory connecting elements of the spectrum of $L^\infty(\Omega, \mu)$ with ultrafilters on a Boolean algebra associated with $(\Omega, \mu)$, generalizing the case where $\Omega$ is a discrete space and its spectrum is identified with the space of ultrafilters on $\Omega$. These ultrafilters were also used to define a new kind of tracial ultrapowers. We also developed a theory representing elements of $L^\infty(\Omega, \mu) \overline{\otimes} M$ as functions that are continuous from the spectrum to $M$ equipped with the weak* topology, analogous to the classical direct integral theory, but which can apply to non-separable von Neumann algebras. We then applied this theory to prove several structural results of isomorphisms $L^\infty(\Omega, \mu) \overline{\otimes} M \to L^\infty(\Omega, \mu) \overline{\otimes} N$ where $M$ and $N$ are $II_1$ factors and placed several constraints on $M$ and $N$. The ultrafilter representation of the spectrum was then used to perform certain calculations to provide specific examples regarding the spectrum representations of tensor products as well as the structures of isomorphisms $L^\infty(\Omega, \mu) \overline{\otimes} M \to L^\infty(\Omega, \mu) \overline{\otimes} N$.

## Section I: Tensor Products, Separability, and Countable Decomposability

**Definition:** Let $M$, $N$ be von Neumann algebras, $\varphi \in M_*$, $a \in M \overline{\otimes} N$. Then $\varphi \otimes Id(a) \in N = (N_*)^*$ is given by $(\varphi \otimes Id(a))(\psi) = \varphi \otimes \psi(a)$ for all $\psi \in N_*$. We also similarly define $Id \otimes \psi(a)$ for $\psi \in N_*$.

**Proposition 1:** Let $M$, $N$ be von Neumann algebras, $a \in M \overline{\otimes} N$, $Q$ be a weak*-closed subspace of $N$. If $a \in M \overline{\otimes} Q$, then for all $\varphi \in M_*$, $\varphi \otimes Id(a) \in Q$.

**Proof:** Assume to the contrary that there exists $\varphi \in M_*$ with $\varphi \otimes Id(a) \notin Q$. By Hahn-Banach, there exists $\psi \in N_*$ s.t. $\psi(Q) = 0$ but $1 = \psi(\varphi \otimes Id(a)) = \varphi \otimes \psi(a)$. However, $a \in M \overline{\otimes} Q$, so $\varphi \otimes \psi(a) = 0$, a contradiction! ∎

**Remark:** For the following theorem, we need to first observe a few technical points regarding abelian von Neumann algebra and Radon measures. In general, we may define Radon measures to be either outer regular on all Borel sets and inner regular only on open sets, or inner regular on all Borel sets. The two definitions are not equivalent, but Radon measures defined one way are in one-to-one correspondence with Radon measures defined the other way (see [Sch73, §I.2] and [Fol99, §7.2, Exercise 14]). Now, in [Tak79, Theorem III.1.18], it is shown that all abelian von Neumann algebras are of the form $L^\infty(\Omega, \mu)$ for some locally compact Hausdorff space $\Omega$ and some Radon measure $\mu$ on $\Omega$. Going through the proof, one sees that the resulting Radon measure has the property of being inner regular on all Borel sets, but not necessarily being outer regular on all Borel sets. Indeed, the proof proceeds by assembling $\Omega$ as a disjoint union of clopen subsets $\Omega_i$ of the spectrum, each one of which is the support of a normal finite measure $\mu_i$. Then the proof proceeds to identify all functions $L^\infty(\Omega, \mu)$ with a continuous function on the spectrum (i.e., an element of the algebra we started with) by performing such identification on each of the $\Omega_i$. For this to work, then, we must have whenever $B_i \subseteq \Omega_i$ are $\mu_i$-null, then $\mu(\bigcup_i B_i) = 0$. If there are only countably many $\Omega_i$, this would not be an issue. However, if there are uncountably many $\Omega_i$, and for each $\Omega_i$ we have a non-empty Borel $\mu_i$-null $B_i \subseteq \Omega_i$, then it is easy to see that requiring $\mu$ to be outer regular on all Borel sets would mean $\mu(\bigcup_i B_i) = \infty$. This could happen, for example, when the algebra we started with is $l^\infty(\kappa) \overline{\otimes} L^\infty([0,1], \lambda)$, where $\kappa$ is an uncountable cardinal and $\lambda$ is the Lebesgue measure on $[0,1]$. Instead, for the proof to work, we must require $\mu(B) = \sum_i \mu_i(B \cap \Omega_i)$ for all Borel sets $B$, at the expense of no longer being outer regular. However, this definition does mean that $\mu$ is inner regular on all Borel sets, which follows from the fact that a finite Radon measure, defined either way, is always regular (see [Sch73, §I.2] or [Fol99, Proposition 7.5]), and that any compact set can only intersect finitely many $\Omega_i$. For this reason, when we mention Radon measures below, we shall always mean a Borel measure on a locally compact set such that it is locally finite and inner regular on all Borel sets.

We observe that this definition has the added benefit that $\mu$ is semifinite, i.e., for all measurable set $E$ with $\mu(E) = \infty$, there exists a measurable $F \subseteq E$ with $0 < \mu(F) < \infty$ [Fol99, §1.3], which follows from local finiteness and inner regularity. We further observe that $\mu$ is decomposable (for definition, see [Fol99, §3.2, Exercise 15]), which follows from the fact that $\Omega$ is the disjoint union of $\Omega_i$, on which $\mu$ restricts to the finite measure $\mu_i$, and that $\mu(B) = \sum_i \mu_i(B \cap \Omega_i)$. These conditions are important. Indeed, we see from [Fol99, §6.2, Exercises 23, 24, & 25] that semifiniteness and decomposability together ensure that $L^1(\Omega, \mu)^* = L^\infty(\Omega, \mu)$ in the usual way, thereby ensuring $L^\infty(\Omega, \mu)$ is actually a von Neumann algebra and its usual representation on $L^2(\Omega, \mu)$ is a normal faithful representation.

Both semifiniteness and decomposability are necessary for this. Indeed, let $I, J$ be uncountable sets and let $\Omega = I \times J$. Let all subsets of $\Omega$ be measurable. Let $\mu(X) = 0$ if $X$ is countable and $\mu(X) = \infty$ otherwise. This measure space is decomposable, as the decomposition into singletons meets all the required conditions. It is clearly not semifinite. We claim that $L^\infty(\Omega, \mu)$ is not a von Neumann algebra. Indeed, for any finite subset $F \subseteq J$, let $P_F$ be the projection onto $I \times F$. Then $P_F$ forms an increasing net of projections. For any projection $P$ with its corresponding measurable set $X$, if $P$ majorizes all $P_F$, then for all $j \in J$, $I \times \{j\} \subseteq X$ up to some countable set, i.e., $I_j \times \{j\} \subseteq X$ for some cocountable $I_j \subseteq I$. It is easy to see that this condition is also sufficient for $P$ to majorize $P_F$. Now, if $L^\infty(\Omega, \mu)$ is a von Neumann algebra, then $P_F$ has a supremum which is a projection. Let it be denoted by $P$ and its corresponding measurable set be $X$. Then as we have seen, for each $j \in J$ there is a cocountable $I_j \subseteq I$ s.t. $X = \bigcup_{j \in J} I_j \times \{j\}$. Since $I_j$ is cocountable, we may pick $i_j \in I_j$ for each $j$. Let $I'_j = I_j \setminus \{i_j\}$, which is still cocountable. Let $X' = \bigcup_{j \in J} I'_j \times \{j\}$. As we have seen, $P_{X'}$ still majorizes all $P_F$, so it must majorizes $P$ as well. As $X' \subseteq X$, we also have $P_{X'} \leq P$, so $P_{X'} = P$. But $X \setminus X' = \bigcup_{j \in J} \{i_j\} \times \{j\}$, which is uncountable, so $\mu(X \setminus X') = \infty$, a contradiction.

On the other hand, let $\Omega$ be an uncountable set. Let the $\sigma$-algebra $\Sigma$ of measurable sets consists of countable and cocountable sets. Let $\mu$ be the restriction of the counting measure. This is clearly semifinite. It is not decomposable, since any decomposition $\mathcal{F} \subseteq \Sigma$ must consist solely of finite sets. By the definition of decomposability, if $X \subseteq \Omega$ is such that $X \cap F$ is measurable for all $F \in \mathcal{F}$, then $X$ is measurable. But as all $F \in \mathcal{F}$ are finite, $X \cap F$ is always measurable, so this would imply all $X \subseteq \Omega$ are measurable, which is not the case. Now, we claim $L^\infty(\Omega, \mu)$ is not a von Neumann algebra. Indeed, fix $X \subseteq \Omega$ which is neither countable nor cocountable. For any finite subset $F \subseteq X$, let $P_F$ be the projection onto $F$. Then $P_F$ forms an increasing net of projections. For any projection $P$ with its corresponding measurable set $Y$, $P$ majorizes all $P_F$ iff $F \subseteq Y$ for all such $F$ iff $X \subseteq Y$. Now, if $L^\infty(\Omega, \mu)$ is a von Neumann algebra, then $P_F$ has a supremum which is a projection. Let it be denoted by $P$ and its corresponding measurable set be $Y$. Then as we have seen, $Y$ is a cocountable set containing $X$. As $X$ is not cocountable, $Y \setminus X$ is nonempty and we may choose $y \in Y \setminus X$. Then $Y' = Y \setminus \{y\}$ is cocountable and $P_{Y'} \leq P$. But as $P_{Y'}$ majorizes all $P_F$, we have $P_{Y'} \geq P$, so $P_{Y'} = P$. But $\mu(Y \setminus Y') = \mu(\{y\}) = 1$, a contradiction.

We also recall here that simple functions are dense in the $L^1$ space for any measure space [Fol99, Proposition 6.7]. ∎

**Theorem 2:** Assume $M$ is either abelian or $M = \mathbb{B}(H)$, then the converse to the above proposition holds, i.e., if for all $\varphi \in M_*$, $\varphi \otimes Id(a) \in Q$, then $a \in M \overline{\otimes} Q$.

**Proof:** Assume to the contrary that $a \notin M \overline{\otimes} Q$. WLOG assume $\|a\| = 1$. By Hahn-Banach, there exists $\Phi \in \left(M \overline{\otimes} N\right)_*$ s.t. $\Phi(M \overline{\otimes} Q) = 0$ but $\Phi(a) = 1$.

We first consider the case where $M$ is abelian. We now represent $M = L^\infty(\Omega, \mu)$ for some locally compact Hausdorff space $\Omega$ and some decomposable Radon measure $\mu$ on $\Omega$. Then $M_* = L^1(\Omega, \mu)$, in which simple functions are dense. As such, we may choose $\Psi \in \left(M \overline{\otimes} N\right)_*$ with $\|\Phi - \Psi\| < \frac{1}{2}$ s.t. $\Psi$ is of the form,

$$\Psi = \sum_{i=1}^{n} E_{X_i} \otimes \psi_i$$

Where $X_i \subseteq \Omega$ is measurable with $0 < \mu(X_i) < \infty$, $E_{X_i} \in M_*$ is given by $E_{X_i}(f) = \int f 1_{X_i} d\mu$, and $\psi_i \in N_*$. We may further assume $X_i$'s are pairwise disjoint. By assumption, $q_i = E_{X_i} \otimes Id(a) \in Q$. Since $\|E_{X_i}\| = \mu(X_i)$, we have $\|q_i\| \leq \mu(X_i)$. Let

$$a' = \sum_{i=1}^n 1_{X_i} \otimes \frac{q_i}{\mu(X_i)} \in M \overline{\otimes} Q$$

Then $\|a'\| \leq 1$. We have,

$$\Psi(a') = \sum_{i=1}^n \sum_{j=1}^n E_{X_i}(1_{X_j}) \psi_i\left(\frac{q_i}{\mu(X_i)}\right)$$

$$= \sum_{i=1}^n \psi_i(q_i)$$

$$= \sum_{i=1}^n \psi_i\left(E_{X_i} \otimes Id(a)\right)$$

$$= \sum_{i=1}^n E_{X_i} \otimes \psi_i(a)$$

$$= \Psi(a)$$

Now, $\Phi(a) = 1$, $\|a\| = 1$, and $\|\Phi - \Psi\| < \frac{1}{2}$, so $|\Psi(a')| = |\Psi(a)| > 1 - \frac{1}{2} = \frac{1}{2}$. But $a' \in M \overline{\otimes} Q$, $\|a'\| \leq 1$, and $\Phi(M \overline{\otimes} Q) = 0$, so $|\Psi(a')| < \frac{1}{2}$, a contradiction!

We now consider the case where $M = \mathbb{B}(H)$. Then $M_* = S^1(H)$ in which finite rank operators are dense. As such, we may choose $\Psi \in (M \overline{\otimes} N)_*$ with $\|\Phi - \Psi\| < \frac{1}{2}$ s.t. $\Psi$ is of the form,

$$\Psi = \sum_{i=1}^n \varphi_i \otimes \psi_i$$

Where $\varphi_i \in M_*$ is represented by a finite rank operator $X_i$ and $\psi_i \in N_*$. Let $H_0 \subseteq H$ be the span of $s_l(X_i)$ and $s_r(X_i)$, $1 \leq i \leq n$. Then $H_0$ is finite-dimensional. Let $\{e_k\}_{k=1}^m$ be an orthonormal basis of $H_0$, $e_{kl}$ be the rank-one operator $|e_k\rangle\langle e_l|$, $E_{kl} \in M_*$ be represented by $e_{lk}$. Since all $X_i$'s are linear combinations of $e_{kl}$'s, we may rewrite $\Psi$ in the form:

$$\Psi = \sum_{1 \leq k,l \leq m} E_{kl} \otimes \psi_{kl}$$

Where $\psi_{kl} \in N_*$. Let $P \in \mathbb{B}(H)$ be the projection onto $H_0$, $a' = PaP \in \mathbb{B}(H_0) \otimes N \subseteq M \overline{\otimes} N$. We observe that $\|a'\| \leq 1$. It is easy to see that in matrix form, the $k, l$-entry of $a'$ is $E_{kl} \otimes Id(a)$, which by assumption is in $Q$. Thus, $a' \in \mathbb{B}(H_0) \overline{\otimes} Q$, whence we may write,

$$a' = \sum_{1 \leq k,l \leq m} e_{kl} \otimes q_{kl}$$

Where $q_{kl} = E_{kl} \otimes Id(a) \in Q$. Now,

$$\Psi(a') = \sum_{1 \leq k,l \leq m} \sum_{1 \leq k',l' \leq m} E_{kl}(e_{k'l'}) \psi_{kl}(q_{k'l'})$$

$$= \sum_{1 \leq k,l \leq m} \psi_{kl}(q_{kl})$$

$$= \sum_{1 \leq k,l \leq m} \psi_{kl}(E_{kl} \otimes Id(a))$$

$$= \sum_{1 \leq k,l \leq m} E_{kl} \otimes \psi_{kl}(a)$$

$$= \Psi(a)$$

The remainder of the argument is the same as in the abelian case. ∎

**Corollary 3:** Let $M, N$ be von Neumann algebras, $a \in M \overline{\otimes} N$, $Q$ be a weak*-closed subalgebra of $N$. If for all $\varphi \in M_*$, $\varphi \otimes Id(a) \in Q$, then $a \in M \overline{\otimes} Q$.

**Proof:** Let $M$ acts on a Hilbert space $H$. Then by [Tak79, Corollary IV.5.10], we have $M \overline{\otimes} Q = (\mathbb{B}(H) \overline{\otimes} Q) \cap (M \overline{\otimes} N)$. Thus, assume to the contrary that $a \notin M \overline{\otimes} Q$, then we must have $a \notin \mathbb{B}(H) \overline{\otimes} Q$. By Theorem 2, there exists $\varphi \in \mathbb{B}(H)_*$ s.t. $\varphi \otimes Id(a) \notin Q$. But $\varphi$ restricts to a normal linear functional on $M$, so $\varphi \otimes Id(a) = \varphi|_M \otimes Id(a) \in Q$, a contradiction! ∎

**Definition:** For $a \in M \overline{\otimes} N$, we define $E_N^{alg}(a)$ to be the smallest weak*-closed subalgebra $Q$ of $N$ s.t. $a \in M \overline{\otimes} Q$. Then the results above show that $E_N^{alg}(a) = W^*(\{\varphi \otimes Id(a) | \varphi \in M_*\})$. If $M$ is abelian or $\mathbb{B}(H)$, we define $E_N^{sp}(a)$ to be the smallest weak*-closed subspace $Q$ of $N$ s.t. $a \in M \overline{\otimes} Q$. The results above show that $E_N^{sp}(a) = \overline{\{\varphi \otimes Id(a) | \varphi \in M_*\}}^{weak^*}$.

**Corollary 4:** If $M$ is separable, then for all $a \in M \overline{\otimes} N$, $E_N^{alg}(a)$ is countably generated as a von Neumann algebra.

**Proof:** Observe that $\|\varphi \otimes Id(a)\| \leq \|\varphi\| \|a\|$. Thus, by choosing a countable dense subset of $M_*$, which is possible since $M$ is separable, we see that $\{\varphi \otimes Id(a) | \varphi \in M_*\}$ is separable. Hence, $E_N^{alg}(a) = W^*(\{\varphi \otimes Id(a) | \varphi \in M_*\})$ is countably generated. ∎

While $E_N^{alg}(a)$ is always countably generated in this case, it is not necessary that $E_N^{alg}(a)$ is separable, as we shall see from the following example.

**Definition ([HVVW16, Definition 1.1.27]):** Let $M$ be a von Neumann algebra acting on $H$, $(\Omega, \Sigma, \mu)$ be a probability space. We say a function $f: \Omega \to M$ is *strongly measurable* if,

1. For all $h, k \in H$, $\Omega \ni \omega \mapsto \langle h, f(\omega)k \rangle$ is measurable;
2. For all $h \in H$, there exists a separable subspace $H_0 \subseteq H$ s.t. $f(\omega)h \in H_0$ $\mu$-a.e.

**Definition:** Let $H$ be a Hilbert space, $(\Omega, \Sigma, \mu)$ be a probability space. Let $L^2(\Omega, \mu; H)$ be the Hilbert space defined as follows: it consists of equivalence classes of all functions $f: \Omega \to H$ s.t.

1. The range is contained in a separable subspace of $H$ $\mu$-a.e.;
2. For any vector $h \in H$, $\Omega \ni \omega \mapsto \langle h, f(\omega) \rangle$ is measurable;
3. $\int_\Omega \|f(\omega)\|^2 \, d\mu < \infty$.

Where we observe that the third condition makes sense as $\|f(\omega)\|^2$ is always measurable assuming $f$ satisfies the first two conditions. Indeed, suppose the range of $f$ lies in $H_0 \subseteq H$ $\mu$-a.e. and $H_0$ is separable. Let $\{e_i\}_{i \in \mathbb{N}}$ be an orthonormal basis of $H_0$. Then $\|f(\omega)\|^2 = \sum_{i=1}^\infty |\langle e_i, f(\omega) \rangle|^2$ $\mu$-a.e. Since $\langle e_i, f(\omega) \rangle$ is measurable, $\|f(\omega)\|^2$ is as well.

We identify two such functions $f, g: \Omega \to H$ if they coincide $\mu$-a.e. An inner product can be defined on $L^2(\Omega, \mu; H)$ by $\langle f, g \rangle = \int_\Omega \langle f(\omega), g(\omega) \rangle \, d\mu$. This is well-defined since, following the same method we used to show $\|f(\omega)\|^2$ is measurable, it is easy to see $\langle f(\omega), g(\omega) \rangle$ is measurable. Furthermore, by Cauchy-Schwarz, it is integrable.

**Remark:** We observe that,

$$L^2(\Omega, \mu) \otimes H = \bigcup_{\substack{H_0 \subseteq H \\ H_0 \text{ separable}}} L^2(\Omega, \mu) \otimes H_0$$

Recall that elements of $L^2(\Omega, \mu) \otimes H_0$ can be canonically represented by measurable functions from $\Omega$ to $H_0$ (see [Con99, Definition 52.3] and [Con99, §52, Exercise 1]). Thus, $L^2(\Omega, \mu; H)$ is canonically isomorphic to $L^2(\Omega, \mu) \otimes H$. As such, it is easy to see that a bounded strongly measurable function $f: \Omega \to M$ acts on $L^2(\Omega, \mu) \otimes H \cong L^2(\Omega, \mu; H)$ by pointwise multiplication. It is easy to see that all operators arising this way commute with both $L^\infty(\Omega, \mu) \otimes 1$ and $1 \otimes M'$. As such, all such operators lie within $\left(L^\infty(\Omega, \mu) \overline{\otimes} M'\right)' = L^\infty(\Omega, \mu) \overline{\otimes} M$, where the equality follows from [Tak79, Theorem IV.5.9]. ∎

Consider $M = l^\infty\big((0, 1]\big)$ where $(0, 1]$ is understood to be a discrete set. It naturally acts on $H = l^2\big((0, 1]\big)$. For $t \in [0, 1]$, let $E_t \in M = l^\infty\big((0, 1]\big)$ be defined by $E_t(s) = 1$ whenever $s > t$ and $E_t(s) = 0$ otherwise. Let $\lambda$ be the Lebesgue measure on $[0, 1]$, $f: ([0, 1], \lambda) \to M$ be defined by $f(t) = E_t$.

**Lemma 5:** $f$ is strongly measurable.

**Proof:** We first show that for all $h, k \in H$, $[0, 1] \ni t \mapsto \langle h, f(t)k \rangle$ is measurable. By taking linear combinations and limits of sequences, it suffices to consider the case where $h = e_u$, $k = e_v$, $u, v \in (0, 1]$, where $e_u$ is the standard basis element of $l^2((0,1])$ supported on $\{u\}$ and similarly for $e_v$. Then we see that $\langle h, f(t)k\rangle = \delta_{u=v} 1_{t<v}$, which is clearly measurable.

Now, fix $h \in H = l^2((0,1])$. It must be supported on a countable subset $S \subseteq (0,1]$. It is then easy to see that $f(t)h \subseteq l^2(S)$, which is separable. ∎

As $f$ is clearly bounded, we may therefore regard $f$ as an element of $L^\infty([0,1], \lambda) \overline{\otimes} M$.

**Proposition 6:** $E_M^{sp}(f) = M$.

**Proof:** Fix a continuous function $g$ on $[0,1]$, which we shall regard as an element of $L^1([0,1], \lambda) = L^\infty([0,1], \lambda)_*$. We now calculate $g \otimes Id(f)$. Observe that,

$$\langle e_u, g \otimes Id(f) e_v\rangle = (g \otimes |e_v\rangle\langle e_u|)(f)$$

$$= (|g\rangle\langle 1| \otimes |e_v\rangle\langle e_u|)(f)$$

$$= \langle 1 \otimes e_u, f(g \otimes e_v)\rangle$$

As a measurable function, $(f(g \otimes e_v))(t) = g(t) E_t e_v = 1_{t<v} g(t) e_v$. As such,

$$\langle e_u, g \otimes Id(f) e_v\rangle = \langle 1 \otimes e_u, f(g \otimes e_v)\rangle = \delta_{u=v} \int_0^v g(t)\, dt$$

Thus, $g \otimes Id(f)$ is a diagonal operator, sending $e_v$ to $\left(\int_0^v g(t)\, dt\right) e_v$. As a function of $v$, we see that $\int_0^v g(t)\, dt$ is a $C^1$-function. In fact, if we let $g$ range over all continuous function on $[0,1]$, then $\int_0^v g(t)\, dt$ ranges over all $C^1$-functions on $[0,1]$ with value $0$ at $v=0$. Such functions form a $*$-algebra, so the weak*-closure $Q$ of the space of all such functions in $M = l^\infty((0,1])$ is a von Neumann subalgebra of $M$. We also have $Q \subseteq E_M^{sp}(f)$.

Now, if we choose $g$ to be the constant function 1, then $g \otimes Id(f)$ is a diagonal operator, sending $e_v$ to $v e_v$. But then for each $v \in (0,1]$, the spectral projection of $g \otimes Id(f)$ associated to $\{v\}$ is the projection onto $e_v$. Such projections generate $M$. Since $Q$ contains $W^*(g \otimes Id(f))$, we conclude that $Q = E_M^{sp}(f) = M$. ∎

**Corollary 7:** Let $M$ be a von Neumann algebra which is not countably decomposable. Then there exists an element $a \in L^\infty([0,1], \lambda) \overline{\otimes} M$ s.t. $E_M^{sp}(a)$ is a non-separable von Neumann subalgebra of $M$, in particular, the predual of $E_M^{sp}(a)$ is unique and not separable.

**Proof:** Since $M$ is not countably decomposable, it admits a collection $\{p_i\}_{i \in I}$ of non-trivial mutually orthogonal projections where $I$ is uncountable. By adding $1 - \sum_{i \in I} p_i$ to one of the projections, we may assume $\sum_{i \in I} p_i = 1$. Furthermore, if the cardinality of $I$ is strictly larger than the continuum, we may choose $J \subseteq I$ whose cardinality is the continuum. Then by adding $\sum_{i \in I \setminus J} p_i$ to $p_j$ for some $j \in J$, we may assume the cardinality of our collection is smaller than or equal to the continuum.

As such, $l^\infty(I) \subseteq M$ and by dualizing an injection $I \hookrightarrow (0, 1]$, we may define a normal *-homomorphism $\pi: l^\infty((0, 1]) \to M$ whose range is $l^\infty(I)$. Let $f \in L^\infty([0, 1], \lambda) \overline{\otimes} l^\infty((0, 1])$ be as defined before. Let $a = Id \otimes \pi(f) \in L^\infty([0, 1], \lambda) \overline{\otimes} M$. We observe that, for $\varphi \in L^\infty([0, 1], \lambda)_*$, $\psi \in M_*$,

$$\psi\big(\varphi \otimes Id(Id \otimes \pi(f))\big) = \varphi \otimes \psi(Id \otimes \pi(f))$$

$$= \varphi \otimes \pi^*\psi(f)$$

$$= \pi^*\psi(\varphi \otimes Id(f))$$

$$= \psi\big(\pi(\varphi \otimes Id(f))\big)$$

That is, $\varphi \otimes Id(a) = \varphi \otimes Id(Id \otimes \pi(f)) = \pi(\varphi \otimes Id(f))$. Since $\overline{\{\varphi \otimes Id(f)\}}^{weak^*} = l^\infty((0, 1])$, we see that $E_M^{sp}(a) = \pi\big(l^\infty((0, 1])\big) = l^\infty(I)$, which is non-separable. ∎

**Remark:** If we instead want to ask about $E_M^{alg}(a)$, the proof can be a lot easier. Indeed, as we have seen, $l^\infty((0, 1])$ is generated as a von Neumann algebra by a single operator $t$, defined by $te_u = ue_u$. Given any normal state $\varphi$ on $L^\infty([0, 1], \lambda)$, we have $\varphi \otimes Id(1 \otimes t) = t$. So, $E_M^{alg}(1 \otimes \pi(t)) = \pi\big(l^\infty((0, 1])\big) = l^\infty(I)$. ∎

We see here a strange phenomenon that $l^\infty((0, 1])$ is countably generated, in fact singly generated, so that it is separable under the weak* topology, yet nevertheless it is not separable in the sense that its predual is not norm-separable. In fact, we have the following result:

**Theorem 8:** Let $M$ be a von Neumann algebra, then the following are equivalent,

1. $M$ admits a faithful representation on a separable Hilbert space;
2. $M$ is separable, i.e., $M_*$ is norm-separable;
3. $M$ is countably decomposable and $M$ is separable under the weak* topology;
4. $M$ is countably decomposable and $M$ admits a norm-separable $C^*$-subalgebra which is weak*-dense in $M$;
5. $M$ is countably decomposable and $M$ is countably generated as a von Neumann algebra.

**Proof:** (3 ⇔ 4 ⇔ 5) Obvious.

(1 ⇒ 2) Let $M$ acts on a separable Hilbert space $H$. Then $M_*$ is a quotient of $\mathbb{B}(H)_* = S^1(H)$. The latter space is norm-separable as $H$ is separable, so $M_*$ is norm-separable.

(2 ⇒ 5) We first show that $M$ is countably decomposable. Indeed, since $M_*$ is norm-separable, the space of normal states on $M$ is norm-separable as well. Thus, we may choose a countable collection $\{\varphi_i\}_{i=1}^\infty$ of normal states which are dense in all normal states. Then $\varphi = \sum_{i=1}^\infty 2^{-i}\varphi_i$ is a normal state and clearly faithful.

Now, to show that $M$ is separable under the weak* topology, let $S(M_*) = \{\varphi \in M_* | \|\varphi\| = 1\}$ be the unit sphere of $M_*$. Since $M_*$ is norm-separable, we may pick a countable dense subset $T \subseteq S(M_*)$. By Hahn-

Banach, for each $t \in T$, there exists an $x_t \in M$ with $\|x_t\| = 1$ and $t(x_t) = 1$. We claim that all such $x_t$ generate $M$. In fact, the weak*-closure of the span of all such $x_t$ is $M$. Assume otherwise. Then by Hahn-Banach there exists $\varphi \in M_*$ with $\varphi(x_t) = 0$ for all $t \in T$ but $\|\varphi\| = 1$. Then there exists a sequence $t_i \to \varphi$. In particular, there exists a $t_i$ s.t. $\|t_i - \varphi\| < 1$. Thus, $|t_i(x_{t_i})| < 1$. But by definition $t_i(x_{t_i}) = 1$, a contradiction!

($4 \Rightarrow 1$) As $M$ is countably decomposable, there exists a normal faithful state $\varphi$. Then the GNS representation $\pi_\varphi$ is normal and faithful. Thus, it suffices to prove $H_\varphi$ is separable. Let $A \subseteq M$ be a norm-separable $C^*$-subalgebra which is weak*-dense in $M$, then $\pi_\varphi(A)$ is SOT-dense in $\pi_\varphi(M)$. Since $A$ is norm-separable, there is a countable subset $S \subseteq A$ s.t. $\pi_\varphi(S)$ is SOT-dense in $\pi_\varphi(M)$. But then $\pi_\varphi(S)\hat{1}$ is dense in $\overline{\pi_\varphi(M)\hat{1}} = H_\varphi$, whence $H_\varphi$ is separable. ∎

**Remark:** As we have seen, the condition that $M$ is countably decomposable is necessary in 3 to 5, since $l^\infty((0,1])$ is singly generated, but, failing to be countably decomposable, its predual is not norm-separable. ∎

Combining Corollary 4, Corollary 7, the remark after Corollary 7, and Theorem 8, we see that,

**Theorem 9:** Let $N$ be a von Neumann algebra. The following are equivalent:

1. $N$ is countably decomposable;
2. For any separable von Neumann algebra $M$, any $a \in M \overline{\otimes} N$, there exists a separable subalgebra $Q \subseteq N$ s.t. $a \in M \overline{\otimes} Q$;
3. There exists a separable von Neumann algebra $M$, s.t. for any $a \in M \overline{\otimes} N$, there exists a separable subalgebra $Q \subseteq N$ s.t. $a \in M \overline{\otimes} Q$;
4. For any $a \in L^\infty([0,1], \lambda) \overline{\otimes} N$, there exists a weak*-closed subspace $Q \subseteq M$ with all its preduals being separable, s.t. $a \in L^\infty([0,1], \lambda) \overline{\otimes} Q$.

In case all algebras involved are tracial, the $1 \Rightarrow 2$ direction can be improved so that the assumption that $M$ is separable can be removed:

**Theorem 10:** Let $M$ and $N$ be tracial von Neumann algebras. Then for any $a \in M \overline{\otimes} N$, there exists separable $M_0 \subseteq M$, $N_0 \subseteq N$ s.t. $a \in M_0 \overline{\otimes} N_0$.

**Proof:** Let $\{e_i\}_{i \in I}$ be a basis of $L^2(M)$ and $\{f_j\}_{j \in J}$ be a basis of $L^2(N)$. Then $\{e_i \otimes f_j\}_{i \in I, j \in J}$ is a basis of $L^2(M \overline{\otimes} N) = L^2(M) \otimes L^2(N)$. Thus, there exists countable $I_0 \subseteq I, J_0 \subseteq J$ s.t. $a = \sum_{i \in I_0, j \in J_0} c_{ij} e_i \otimes f_j$ where the convergence is in $L^2(M)$. For each $i \in I$, $e_i$ is the $L^2$-limit of a sequence of elements in $M$, so it is contained in the $L^2$ space of a countably generated subalgebra of $M$. As $I_0$ is countable, there exists countably generated $M_0 \subseteq M$ s.t. $e_i \in L^2(M_0)$ for all $i \in I_0$. By Theorem 8, $M_0$ is separable. Similarly, there exists separable $N_0 \subseteq N$ s.t. $f_j \in L^2(N_0)$ for all $j \in J_0$. Hence, $a \in L^2(M_0) \otimes L^2(N_0) = L^2(M_0 \overline{\otimes} N_0)$. Finally, we observe that $(M \overline{\otimes} N) \cap L^2(M_0 \overline{\otimes} N_0) = M_0 \overline{\otimes} N_0$, whence the result follows. (It is a general fact that for a tracial algebra $P$ and a subalgebra $Q \subseteq P$, $L^2(Q) \cap P = Q$. Indeed, let $x \in L^2(Q) \cap P$. Let $E_Q$ be the conditional expectation onto $Q$. Then $E_Q$ is $L^2$-norm-decreasing and therefore extends to the orthogonal projection from $L^2(P)$ onto $L^2(Q)$. Hence, as $x \in L^2(Q)$, $x = E_Q(x) \in Q$.) ∎

# Section II: Ultrafilters, Ultrapowers, and the Spectrum of an Abelian von Neumann Algebra

**Definition:** Let $(\Omega, \Sigma, \mu)$ be a measure space. Define $\mathfrak{B}((\Omega, \mu))$ to be the Boolean algebra whose elements are equivalence classes of measurable sets, where two measurable sets are identified if their symmetric difference is $\mu$-null. Meet, join, and complement are given by intersections, unions, and set complements (of representatives of the equivalence classes), respectively. 0 and 1 are given by the equivalence classes of $\emptyset$ and $\Omega$, respectively.

**Remark:** The Boolean prime ideal theorem implies that any Boolean algebra admits an ultrafilter. ∎

**Notation:** While elements of $\mathfrak{B}((\Omega, \mu))$ are technically equivalence classes, we shall, by a slight abuse of notation, mean $F$ is a representative of an equivalence class in $\mathfrak{B}((\Omega, \mu))$ when we write $F \in \mathfrak{B}((\Omega, \mu))$. The same shall hold for ultrafilters on $\mathfrak{B}((\Omega, \mu))$ or any other subset of $\mathfrak{B}((\Omega, \mu))$, unless otherwise indicated.

**Definition:** Let $(\Omega, \Sigma, \mu)$ be a semifinite, decomposable measure space, $\mathcal{U}$ be an ultrafilter on $\mathfrak{B}((\Omega, \mu))$. Let $f \in L^\infty(\Omega, \mu)$, $a \in \mathbb{C}$. We say $\lim_{\mathcal{U}} f = a$ if, for any $\varepsilon > 0$, there exists $F \in \mathcal{U}$ s.t. for all measurable $X \subseteq F$ with $0 < \mu(X) < \infty$, we have,

$$\left|\left(\frac{1}{\mu(X)}\int_X f \, d\mu\right) - a\right| < \varepsilon$$

**Remark:** Observe that when $\Omega$ is a discrete space and $\mu$ is the counting measure, this is the same as the usual definition of ultrafilter limit. We also note that any $\sigma$-finite measure space is automatically semifinite and decomposable. ∎

We shall assume, in the rest of this section, that any measure space mentioned is semifinite and decomposable.

**Lemma 1:** $\lim_{\mathcal{U}} f$, if it exists, is unique.

**Proof:** Suppose $\lim_{\mathcal{U}} f = a_1$ and $\lim_{\mathcal{U}} f = a_2$ but $a_1 \neq a_2$. Then we may pick $\varepsilon > 0$ s.t. $2\varepsilon < |a_1 - a_2|$. By definition, then, there exists $F_i \in \mathcal{U}$ satisfying the condition in the definition of ultrafilter limit for $a_i$, for each $i = 1, 2$ respectively. But then $F_1 \cap F_2 \in \mathcal{U}$, so in particular $F_1 \cap F_2$ is not null. By semifiniteness, we may pick a measurable $X \subseteq F_1 \cap F_2$ with $0 < \mu(X) < \infty$. But then $\left|\left(\frac{1}{\mu(X)}\int_X f \, d\mu\right) - a_1\right| < \varepsilon$ and $\left|\left(\frac{1}{\mu(X)}\int_X f \, d\mu\right) - a_2\right| < \varepsilon$ simultaneously. Since $2\varepsilon < |a_1 - a_2|$, this yields a contradiction. ∎

**Lemma 2:** For any $\varepsilon > 0$ and any measurable $X$ with $\mu(X) > 0$, there exists $S \subseteq X$ with $\mu(S) > 0$ and $|f(s) - f(t)| < \varepsilon$ whenever $s, t \in S$.

**Proof:** We observe that the complex plane can be covered by countably many open balls of radius $\frac{\varepsilon}{2}$. Thus, the inverse image under $f|_X$ of at least one of these open balls must have positive measure. It is easy to see that such an inverse image satisfies the requirement of the lemma. ∎

**Remark:** Note that $\mathfrak{B}((\Omega, \mu))$ is naturally isomorphic to the Boolean algebra of projections in $L^\infty(\Omega, \mu)$. Since we are assuming semifiniteness and decomposability, $L^\infty(\Omega, \mu)$ is a von Neumann algebra, so $\mathfrak{B}((\Omega, \mu))$ is a complete lattice. This fact is important in proving the following theorem. ∎

**Theorem 3:** $\lim_{\mathcal{U}} f$ exists for all $f \in L^\infty(\Omega, \mu)$.

**Proof:** WLOG we may assume $\|f\| \leq 1$. Let $D$ be the unit disc of the complex plane. By making alterations on a null set if necessary, we may assume $f$ always takes value in $D$. We claim that there exists an $a \in D$ s.t. $\lim_{\mathcal{U}} f = a$. Assume to the contrary, then by definition, for each $a \in D$, there exists an $\varepsilon_a > 0$ s.t., for all $F \in \mathcal{U}$, there exists a measurable $X \subseteq F$ with $0 < \mu(X) < \infty$ and,

$$\left| \left( \frac{1}{\mu(X)} \int_X f \, d\mu \right) - a \right| > \varepsilon_a$$

Since $B\left(a, \frac{\varepsilon_a}{2}\right)$ forms an open cover of $D$, there exists $\{a_i\}_{i=1}^n \subseteq D$ s.t. $B\left(a_i, \frac{\varepsilon_{a_i}}{2}\right)$'s cover $D$. Now, for each $i = 1, \cdots, n$, let $F_i \in \mathfrak{B}((\Omega, \mu))$ be defined by,

$$F_i = \bigvee \left\{ X \in \mathfrak{B}((\Omega, \mu)) : \forall \text{measurable } Y \subseteq X \text{ with } 0 < \mu(Y) < \infty, \left| \left( \frac{1}{\mu(Y)} \int_Y f \, d\mu \right) - a_i \right| > \varepsilon_{a_i} \right\}$$

Which exists by the fact that $\mathfrak{B}((\Omega, \mu))$ is a complete lattice. We now show that $\neg F_i \notin \mathcal{U}$. Assume to the contrary that $\neg F_i \in \mathcal{U}$. Then there exists a measurable $X \subseteq \neg F_i$ with $0 < \mu(X) < \infty$ and $\left| \left( \frac{1}{\mu(X)} \int_X f \, d\mu \right) - a_i \right| > \varepsilon_{a_i}$. Let $\eta = \frac{1}{2} \left( \left| \left( \frac{1}{\mu(X)} \int_X f \, d\mu \right) - a_i \right| - \varepsilon_{a_i} \right) > 0$. Let,

$$A = \{ t \in X : |f(t) - a_i| > \varepsilon_{a_i} + \eta \}$$

$A$ is measurable. Furthermore, it must be of positive measure, as otherwise $\left| \left( \frac{1}{\mu(X)} \int_X f \, d\mu \right) - a_i \right| \leq \varepsilon_{a_i} + \eta < \left| \left( \frac{1}{\mu(X)} \int_X f \, d\mu \right) - a_i \right|$, a contradiction. Thus, by Lemma 2, there exists $Y \subseteq A$ with $\mu(Y) > 0$ and $|f(s) - f(t)| < \eta$ whenever $s, t \in Y$. Fix $t \in Y$. Then for all $S \subseteq Y$ with $0 < \mu(S) < \infty$,

$$\left| \left( \frac{1}{\mu(S)} \int_S f \, d\mu \right) - a_i \right| \geq |f(t) - a_i| - \left| \left( \frac{1}{\mu(S)} \int_S f \, d\mu \right) - f(t) \right| > \varepsilon_{a_i} + \eta - \eta = \varepsilon_{a_i}$$

By definition, $Y \subseteq F_i$. But we also have $Y \subseteq A \subseteq X \subseteq \neg F_i$. As $\mu(Y) > 0$, this is a contradiction. Hence, $\neg F_i \notin \mathcal{U}$, so as $\mathcal{U}$ is an ultrafilter, $F_i \in \mathcal{U}$. Thus, $\bigwedge_{i=1}^n F_i \in \mathcal{U}$.

Now, fix an $\varepsilon < \frac{1}{2} \min\{\varepsilon_{a_1}, \cdots, \varepsilon_{a_n}\}$. Since $\bigwedge_{i=1}^n F_i \in \mathcal{U}$, $\mu(\bigwedge_{i=1}^n F_i) > 0$. We may apply Lemma 2 to obtain $S \subseteq \bigwedge_{i=1}^n F_i$ with $\mu(S) > 0$ and $|f(s) - f(t)| < \varepsilon$ whenever $s, t \in S$. Fix $t \in S$. $f(t) \in D$, so there exists $i \in \{1, \cdots, n\}$ s.t. $f(t) \in B\left(a_i, \frac{\varepsilon_{a_i}}{2}\right)$. Thus, all $s \in S$ satisfies $|f(s) - a_i| < \varepsilon + \frac{\varepsilon_{a_i}}{2}$. Therefore, for any measurable $Y \subseteq S$ with $0 < \mu(Y) < \infty$, we must have,

$$\left| \left( \frac{1}{\mu(Y)} \int_Y f \, d\mu \right) - a_i \right| \leq \varepsilon + \frac{\varepsilon_{a_i}}{2} < \varepsilon_{a_i}$$

Now, fix any $X \in \mathcal{B}((\Omega, \mu))$ satisfying the condition in the definition of $F_i$, i.e., for all measurable $Y \subseteq X$ with $0 < \mu(Y) < \infty$, we have $\left|\left(\frac{1}{\mu(Y)} \int_Y f \, d\mu\right) - a_i\right| > \varepsilon_{a_i}$. But then clearly $\mu(X \cap S) = 0$, i.e., $S \subseteq \neg X$. But this means $S$ is contained in the infimum of all such $\neg X$, which is the complement of $F_i$, i.e., $S \subseteq \neg F_i$, which contradicts $S \subseteq \bigwedge_{i=1}^n F_i$. ∎

By Lemma 1 and Theorem 3 then, $\lim_{\mathcal{U}} f$ is always a well-defined number. We record a few easy observations about it:

**Proposition 4:**

1. $\lim_{\mathcal{U}} f$ is linear in $f$, i.e., $\lim_{\mathcal{U}}(cf + g) = c \lim_{\mathcal{U}} f + \lim_{\mathcal{U}} g$ for all $f, g \in L^\infty(\Omega, \mu)$, $c \in \mathbb{C}$;
2. $\lim_{\mathcal{U}} \bar{f} = \overline{\lim_{\mathcal{U}} f}$;
3. $\lim_{\mathcal{U}} 1 = 1$;
4. $\left|\lim_{\mathcal{U}} f\right| \leq \|f\|$.

Hence, $\lim_{\mathcal{U}} f$ defines a state on $L^\infty(\Omega, \mu)$, which we shall denote as $\omega_{\mathcal{U}}$.

**Theorem 5:** $\lim_{\mathcal{U}} fg = \left(\lim_{\mathcal{U}} f\right)\left(\lim_{\mathcal{U}} g\right)$ for all $f, g \in L^\infty(\Omega, \mu)$, i.e., $\omega_{\mathcal{U}}$ is a character.

**Proof:** Let $a = \lim_{\mathcal{U}} f$, $b = \lim_{\mathcal{U}} g$. Let $\varepsilon > 0$. We fix $\varepsilon_0 > 0$ s.t. $(|a| + |b| + \varepsilon_0)\varepsilon_0 < \varepsilon$. Let $F = f^{-1}(B(a, \varepsilon_0))$ and $G = g^{-1}(B(b, \varepsilon_0))$. We claim that $F \in \mathcal{U}$. Assume otherwise. Then $\neg F = f^{-1}(\mathbb{C} \setminus \overline{B(a, \varepsilon_0)}) \in \mathcal{U}$. Since $a = \lim_{\mathcal{U}} f$, there exists $E \in \mathcal{U}$ s.t. for all measurable $Y \subseteq E$ with $0 < \mu(Y) < \infty$, we have $\left|\left(\frac{1}{\mu(Y)} \int_Y f \, d\mu\right) - a\right| < \varepsilon_0$. Since $E \cap \neg F \in \mathcal{U}$, it has positive measure.

Now, $\mathbb{C} \setminus \overline{B(a, \varepsilon_0)}$ is an open subset of $\mathbb{C}$, so it is a union of countably many open balls $\{B_i\}_{i=1}^\infty$. Thus, $E \cap \neg F = E \cap f^{-1}(\bigcup_{i=1}^\infty B_i) = \bigcup_{i=1}^\infty (E \cap f^{-1}(B_i))$. As such, for some $i$, $E \cap f^{-1}(B_i)$ must have positive measure. Fix such an $i$ and by semifiniteness choose $S \subseteq E \cap f^{-1}(B_i)$ s.t. $0 < \mu(S) < \infty$. Then we see that $\frac{1}{\mu(S)} \int_S f \, d\mu \in \overline{B_i} \subseteq \mathbb{C} \setminus B(a, \varepsilon_0)$, so $\left|\left(\frac{1}{\mu(S)} \int_S f \, d\mu\right) - a\right| \geq \varepsilon_0$, contradicting the condition on $E$.

Hence, $F \in \mathcal{U}$. Similarly, $G \in \mathcal{U}$. Consider $F \cap G \in \mathcal{U}$. For all measurable $Y \subseteq F \cap G$ with $0 < \mu(Y) < \infty$, we have,

$$\left|\left(\frac{1}{\mu(Y)} \int_Y fg \, d\mu\right) - ab\right| \leq \frac{1}{\mu(Y)} \int_Y |f(t) - a||g(t)| \, d\mu(t) + \frac{|a|}{\mu(Y)} \int_Y |g(t) - b| \, d\mu(t)$$

$$\leq \varepsilon_0(|b| + \varepsilon_0) + |a|\varepsilon_0$$

$$< \varepsilon$$

Since $\varepsilon > 0$ is arbitrary, this concludes the proof. ∎

Let the spectrum of $L^\infty(\Omega, \mu)$ be $\Lambda$. Then Theorem 5 shows that $\omega_{\mathcal{U}} \in \Lambda$. Now, we observe that clopen subsets of $\Lambda$ correspond to projections in $C(\Lambda) = L^\infty(\Omega, \mu)$, which then correspond to elements of $\mathcal{B}((\Omega, \mu))$. Let this one-to-one correspondence be denoted by $\mathfrak{b}$: {clopen subsets of $\Lambda$} $\to \mathcal{B}((\Omega, \mu))$.

**Proposition 6:** Let $\omega \in \Lambda$. Then the following defines an ultrafilter on $\mathfrak{B}((\Omega, \mu))$:

$$\mathcal{U}_\omega = \mathfrak{b}(\{\text{clopen neighborhoods of } \omega\}) \subseteq \mathfrak{B}((\Omega, \mu))$$

**Proof:** This follows from observing that $\mathfrak{b}$ is an isomorphism of Boolean algebras and that $\{\text{clopen neighborhoods of } \omega\}$ is easily seen to be an ultrafilter in the Boolean algebras of all clopen subsets of $\Lambda$. ∎

We shall now show that ultrafilters on $\mathfrak{B}((\Omega, \mu))$ and elements of the spectrum are in one-to-one correspondence with each other, via the maps we just defined: $\mathcal{U} \mapsto \omega_\mathcal{U}$ and $\omega \mapsto \mathcal{U}_\omega$. Before we prove it, we need to recall certain facts about the spectrum of an abelian von Neumann algebras, namely that it is a hyperstonean space [Tak79, Theorem III.1.18]. In particular, being a stonean space, it is a compact Hausdorff space in which the closure of any open set is open. It is easy to see that in such a space, clopen sets form a basis for the topology.

**Theorem 7:** Ultrafilters on $\mathfrak{B}((\Omega, \mu))$ and elements of $\Lambda$ are in one-to-one correspondence with each other via the following two maps, which are inverses of each other:

$$\{\text{Ultrafilters on } \mathfrak{B}((\Omega, \mu))\} \ni \mathcal{U} \mapsto \omega_\mathcal{U} \in \Lambda$$

$$\Lambda \ni \omega \mapsto \mathcal{U}_\omega \in \{\text{Ultrafilters on } \mathfrak{B}((\Omega, \mu))\}$$

**Proof:** We need to show $\omega_{\mathcal{U}_\omega} = \omega$ and $\mathcal{U}_{\omega_\mathcal{U}} = \mathcal{U}$ for all $\omega \in \Lambda$ and ultrafilter $\mathcal{U}$ on $\mathfrak{B}((\Omega, \mu))$. For the first part, let $\omega \in \Lambda$. We need to show for all $f \in L^\infty(\Omega, \mu)$, $\lim_{\mathcal{U}_\omega} f = \omega(f)$. Fix $\varepsilon > 0$. Consider $f$ as a continuous function on $\Lambda$. Noting that clopen sets form a basis for the topology on $\Lambda$, we may choose a clopen neighborhood $F$ of $\omega$ s.t. $|f(t) - \omega(f)| < \varepsilon$ for all $t \in F$. Since $F$ is clopen, $1_F$ is a continuous function, so we have $\|f 1_F - \omega(f) 1_F\| < \varepsilon$. Now, $\mathfrak{b}(F) \in \mathcal{U}_\omega$. We observe that, by definition, $1_F$ and $1_{\mathfrak{b}(F)}$ are the same element. Therefore, for all measurable $X \subseteq \mathfrak{b}(F)$ with $0 < \mu(X) < \infty$, we have,

$$\left|\left(\frac{1}{\mu(X)} \int_X f \, d\mu\right) - \omega(f)\right| \leq \frac{1}{\mu(X)} \int_X |f - \omega(f)| \, d\mu$$

$$= \frac{1}{\mu(X)} \int_X |f 1_F - \omega(f) 1_F| \, d\mu$$

$$\leq \|f 1_F - \omega(f) 1_F\|$$

$$< \varepsilon$$

Since $\varepsilon > 0$ is arbitrary, $\lim_{\mathcal{U}_\omega} f = \omega(f)$.

For the second part of the theorem, let $\mathcal{U}$ be an ultrafilter on $\mathfrak{B}((\Omega, \mu))$. Assume to the contrary that $\mathcal{U}_{\omega_\mathcal{U}} \neq \mathcal{U}$. Then there exists $F \in \mathfrak{B}((\Omega, \mu))$ s.t. $F \in \mathcal{U}_{\omega_\mathcal{U}}$ but $\neg F \in \mathcal{U}$. Consider $1_F \in L^\infty(\Omega, \mu)$. As a continuous function on $\Lambda$, it is $1_{\mathfrak{b}^{-1}(F)}$, since $F \in \mathcal{U}_{\omega_\mathcal{U}}$, by definition $\mathfrak{b}^{-1}(F)$ is a clopen neighborhood of $\omega_\mathcal{U}$, so $\omega_\mathcal{U}(1_F) = 1$. On the hand, for all measurable $X \subseteq \neg F$ with $0 < \mu(X) < \infty$, we have $\frac{1}{\mu(X)} \int_X 1_F \, d\mu = 0$. As $\neg F \in \mathcal{U}$, this implies $\omega_\mathcal{U}(1_F) = \lim_\mathcal{U} f = 0$, a contradiction! ∎

**Remark:** In general, there are two kinds of ultrafilters on $\mathfrak{B}\big((\Omega,\mu)\big)$. Assume $(\Omega,\mu)$ admits atoms, i.e., a measurable set $A \subseteq \Omega$ such that $\mu(A) > 0$ and for all measurable $B \subseteq A$ with $\mu(B) < \mu(A)$, we must have $\mu(B) = 0$. Note that by assumption of semifiniteness, this implies $\mu(A) < \infty$ and that for all measurable $B \subseteq A$, either $\mu(B) = 0$ or $\mu(A \setminus B) = 0$, so the atoms correspond to minimal projections in $L^\infty(\Omega,\mu)$. In such cases we can have, for each atom $A$, the *principal ultrafilter at $A$*, which consists of all measurable sets majorizing $A$. The ultrafilter limit then is simply integration over $A$, divided by $\mu(A)$. (Alternatively, it is easy to see any measurable function must be constant on $A$, $\mu$-a.e., and that constant is the ultrafilter limit.) These are the normal characters on $L^\infty(\Omega,\mu)$. Indeed, characters are extreme points of all states (see [Mur90, Theorem 5.1.6(2)] and [Mur90, Theorem 5.1.8]), so a normal character would have to be an extreme point of the unit ball of $L^1(\Omega,\mu)$. It is then easy to see that such an extreme point has to be supported on an atom. It is also easy to see that these ultrafilter limits would correspond to isolated points in the spectrum.

We shall call all other ultrafilters *free*. Note that decomposability guarantees that the union of any collection of pairwise essentially disjoint atoms is measurable, which follows from the fact that a set of finite measure can only contain countably many atoms. As such, we may consider $\Omega$ minus all the atoms. If it has positive measure, then we say $(\Omega,\mu)$ has an *atomless part*, which corresponds to a diffuse summand of $L^\infty(\Omega,\mu)$. In this case a free ultrafilter exists by applying Zorn's lemma. The *atomic* part, i.e., the union of all atoms, corresponds to a summand of $L^\infty(\Omega,\mu)$ isomorphic to $l^\infty(\kappa)$, where $\kappa$ is the cardinality of the set of equivalence classes of all atoms. As long as $\kappa$ is not finite, a free ultrafilter also exists and can be chosen supported solely in the atomic part, i.e., the corresponding character factors through $l^\infty(\kappa)$. In this case they simply correspond to free ultrafilters on the discrete set $\kappa$. Again, this depends on Zorn's lemma. It is therefore all but impossible to tangibly exhibit a free ultrafilter. For this decomposition into the atomless part and the atomic part, see also [Tak79, Proposition III.1.19]. ∎

Since $L^\infty(\Omega,\mu) \cong L^\infty(\Omega,\mu')$ naturally when $\mu$ and $\mu'$ are equivalent and the map $\Lambda \ni \omega \mapsto \mathcal{U}_\omega$ in Proposition 6 and Theorem 7 also has no dependence on the measure chosen, we have,

**Corollary 8:** Ultrafilter limits do not change if $\mu$ is replaced by an equivalent semifinite decomposable measure.

In light of [Tak79, Theorem III.1.18] and the remark preceding Theorem I.2, we also have the following corollary:

**Corollary 9:** Given any hyperstonean space $\Lambda$, for the Boolean algebra of clopen subsets of $\Lambda$, the map sending any given element of $\Lambda$ to the principal ultrafilter at that point is a bijection, i.e., given any ultrafilter, there must be a unique $\omega \in \Lambda$ s.t. the ultrafilter consists of all clopen neighborhoods of $\omega$.

**Remark:** We observe that, for an ultrafilter $\mathcal{U}$ on $\mathfrak{B}\big((\Omega,\mu)\big)$ and $F \in \mathfrak{B}\big((\Omega,\mu)\big)$, $\omega_\mathcal{U} \in \mathfrak{b}^{-1}(F)$ iff $F \in \mathcal{U}$. Indeed, this follows from Theorem 7, since elements of $\mathcal{U} = \mathcal{U}_{\omega_\mathcal{U}}$ correspond to clopen neighborhoods of $\omega_\mathcal{U}$. As clopen sets form a topological basis for $\Lambda$, we see that $\Lambda$ is naturally homeomorphic to the Stone space of $\mathfrak{B}\big((\Omega,\mu)\big)$. In that light, Corollary 9 is then simply a special case of the Stone representation theorem (see [BS81, Theorem IV.4.6]). ∎

Similar to how ultrapowers are defined for a tracial algebra and an ultrafilter on a discrete set, we shall do the same in the context of ultrafilters on $\mathfrak{B}\big((\Omega,\mu)\big)$:

**Definition:** Let $M$ be a tracial von Neumann algebra with normal faithful tracial state $\tau$, $\mathcal{U}$ be an ultrafilter on $\mathfrak{B}\big((\Omega,\mu)\big)$, we define the *(tracial) ultrapower of $(M,\tau)$ along $\mathcal{U}$*, denoted by $(M,\tau)^\mathcal{U}$, to be,

$$(M,\tau)^{\mathcal{U}} = \left(L^{\infty}(\Omega,\mu) \overline{\otimes} M\right)/I_{\mathcal{U},\tau}$$

Where the ideal $I_{\mathcal{U},\tau} \subseteq L^{\infty}(\Omega,\mu) \overline{\otimes} M$ is defined by,

$$I_{\mathcal{U},\tau} = \left\{a \in L^{\infty}(\Omega,\mu) \overline{\otimes} M : \lim_{\mathcal{U}}(Id \otimes \tau(a^*a)) = 0\right\}$$

If $\omega$ is a character on $L^{\infty}(\Omega,\mu)$, then we define the (*tracial*) *ultrapower of* $(M,\tau)$ *along* $\omega$, denoted by $(M,\tau)^{\omega}$, to be $(M,\tau)^{\mathcal{U}_\omega}$. When $\tau$ is clear from the context, we shall omit $\tau$ and simply write $M^{\mathcal{U}}$ and $M^{\omega}$ for $(M,\tau)^{\mathcal{U}}$ and $(M,\tau)^{\omega}$, respectively.

This definition is consistent with the usual tracial ultrapowers for ultrafilters on a discrete set. To show this is well-defined, we need to prove $I_{\mathcal{U},\tau}$ is indeed a norm-closed ideal:

**Lemma 10:** $I_{\mathcal{U},\tau}$ is a norm-closed ideal.

**Proof:** Item 4 of Proposition 4 implies $I_{\mathcal{U},\tau}$ is norm-closed. To show it is an ideal, we prove $Id \otimes \tau$ is positive and tracial, i.e., $Id \otimes \tau(a) \geq 0$ if $a \geq 0$ and $Id \otimes \tau(ab) = Id \otimes \tau(ba)$. For the first statement, let $\psi$ be a normal positive linear functional on $L^{\infty}(\Omega,\mu)$, then $\psi \otimes \tau$ is positive, so $\psi(Id \otimes \tau(a)) = \psi \otimes \tau(a) \geq 0$ when $a \geq 0$. As $\psi$ is arbitrary, $Id \otimes \tau(a) \geq 0$.

For the second statement, let $\psi \in L^{\infty}(\Omega,\mu)_*$. Then as $L^{\infty}(\Omega,\mu)$ is abelian and $\tau$ is tracial, we see that $\psi \otimes \tau$ is tracial, so $\psi(Id \otimes \tau(ab)) = \psi \otimes \tau(ab) = \psi \otimes \tau(ba) = \psi(Id \otimes \tau(ba))$. As $\psi$ is arbitrary, $Id \otimes \tau(ab) = Id \otimes \tau(ba)$.

Since $\lim_{\mathcal{U}}$ defines a character, we see that $\lim_{\mathcal{U}}(Id \otimes \tau(a^*b)) = 0$ defines a positive sesquilinear form, so by Cauchy-Schwarz, $I_{\mathcal{U},\tau}$ is a linear space. Positivity then implies $I_{\mathcal{U},\tau}$ is a left ideal. As $\lim_{\mathcal{U}}(Id \otimes \tau(a^*a)) = \lim_{\mathcal{U}}(Id \otimes \tau(aa^*))$, it is also a right ideal. ∎

**Notation:** For an element $x \in (M,\tau)^{\mathcal{U}}$, we write $x = (a)$ for $a \in L^{\infty}(\Omega,\mu) \overline{\otimes} M$ to mean $a$ is a representative of $x \in (M,\tau)^{\mathcal{U}} = \left(L^{\infty}(\Omega,\mu) \overline{\otimes} M\right)/I_{\mathcal{U},\tau}$.

**Definition:** We define a tracial state $\tau_{\mathcal{U}}$ on $(M,\tau)^{\mathcal{U}}$ by $\tau_{\mathcal{U}}((a)) = \lim_{\mathcal{U}}(Id \otimes \tau(a))$. It is easy to see that this is a well-defined faithful tracial state on $(M,\tau)^{\mathcal{U}}$.

We shall prove in Section IV that $(M,\tau)^{\mathcal{U}}$ is a tracial von Neumann algebra and that $\tau_{\mathcal{U}}$ is a normal faithful tracial state on it. Furthermore, we will show that if $M$ is a factor, so is $M^{\mathcal{U}}$.

# Section III: More on the Spectrum of $L^\infty([0, 1], \lambda)$

As mentioned in the remark after Theorem II.7, an abelian von Neumann algebra can be decomposed into the direct sum of its atomic part and its atomless (diffuse) part. If the algebra is assumed to be separable, then the atomic part, if it exists, must be isomorphic to $l^\infty$ on a finite set or $l^\infty(\mathbb{N})$. (As we have observed, the atomic part is $l^\infty(\kappa)$, where $\kappa$ is the cardinality of the set of equivalence classes of all atoms. If $\kappa > |\mathbb{N}|$, then the algebra is no longer separable.) $l^\infty$ on a finite set has no free ultrafilters. Its spectrum is simply the finite set itself. The spectrum of $l^\infty(\mathbb{N})$ is the Stone-Čech compactification of $\mathbb{N}$. In either case, the spectrum is relatively well-understood. The atomless part, if it exists, must be isomorphic to $L^\infty([0, 1], \lambda)$ [AP18, Theorem 3.2.4]. Since the spectrum is the disjoint union of the spectrum of the atomic part and the spectrum of the atomless part, the task of understanding the spectrums of separable abelian von Neumann algebras mostly reduced to the issue of understanding the spectrum of $L^\infty([0, 1], \lambda)$.

Before we begin to undertake this task, we first need a better understanding of filters on $\mathfrak{B}((\Omega, \mu))$. In the following, we shall continue the assumptions in Section II, i.e., we shall assume $(\Omega, \Sigma, \mu)$ is a semifinite, decomposable measure space and $\Lambda$ is the spectrum of $L^\infty(\Omega, \mu)$. Recall the notation from Proposition II.6 that the one-to-one correspondence between clopen subsets of $\Lambda$ and elements of $\mathfrak{B}((\Omega, \mu))$ shall be denoted by $\mathfrak{b}: \{\text{clopen subsets of } \Lambda\} \to \mathfrak{B}((\Omega, \mu))$. A filter shall be called *principal at $A$* for some $A \in \mathfrak{B}((\Omega, \mu))$ if the filter consists of all measurable sets majorizing $A$. It is easy to see that if an ultrafilter is principal at $A$, then $A$ must be an atom, so it is the principal ultrafilter at $A$, in the sense of the definition in the remark after Theorem II.7. Note that principal filters include the nonproper filter $\mathfrak{B}((\Omega, \mu))$, when $A = \emptyset$. We also note that the principal filter at $A \in \mathfrak{B}((\Omega, \mu))$ and the principal filter at $B \in \mathfrak{B}((\Omega, \mu))$ are equal iff $A = B$ up to a null set difference.

**Definition:** For a filter $\mathcal{F} \subseteq \mathfrak{B}((\Omega, \mu))$, we define $\Omega_\mathcal{F} = \bigcap \{\mathfrak{b}^{-1}(F): F \in \mathcal{F}\} \subseteq \Lambda$. For a closed subset $K \subseteq \Lambda$, we define $\mathcal{F}_K = \{\mathfrak{b}(F): F \subseteq \Lambda \text{ is clopen}, K \subseteq F\}$.

**Theorem 1:** Filters on $\mathfrak{B}((\Omega, \mu))$ and closed subsets of $\Lambda$ are in one-to-one correspondence with each other via the following two maps, which are inverses of each other:

$$\{\text{Filters on } \mathfrak{B}((\Omega, \mu))\} \ni \mathcal{F} \mapsto \Omega_\mathcal{F} \in \{\text{Closed subsets of } \Lambda\}$$

$$\{\text{Closed subsets of } \Lambda\} \ni K \mapsto \mathcal{F}_K \in \{\text{Filters on } \mathfrak{B}((\Omega, \mu))\}$$

Furthermore, $\mathcal{F}$ is proper iff $\Omega_\mathcal{F} \neq \emptyset$. $\mathcal{F}$ is an ultrafilter iff $\Omega_\mathcal{F}$ is a singleton, in which case $\Omega_\mathcal{F} = \{\omega_\mathcal{F}\}$. $\mathcal{F}$ is principal iff $\Omega_\mathcal{F}$ is clopen, in which case $\mathcal{F}$ is principal at $\mathfrak{b}(\Omega_\mathcal{F})$.

**Proof:** We first observe that $\Omega_\mathcal{F}$, being an intersection of clopen sets, is indeed closed. That $\mathcal{F}_K$ is a filter is also clear. We now show $\Omega_{\mathcal{F}_K} = K$. This simply follows from the observation that, as clopen sets form a basis for the topology, any closed set is the intersection of all its clopen neighborhoods. We then show $\mathcal{F}_{\Omega_\mathcal{F}} = \mathcal{F}$. It is clear that $\mathcal{F} \subseteq \mathcal{F}_{\Omega_\mathcal{F}}$. Assume to the contrary that there exists $G \in \mathcal{F}_{\Omega_\mathcal{F}} \setminus \mathcal{F}$. Then the following is a filter: $\mathcal{F}(\neg G) = \{X \in \mathfrak{B}((\Omega, \mu)): \exists F \in \mathcal{F} \text{ s.t. } F \cap \neg G \subseteq X\}$. We observe that it must be proper, as otherwise $F \cap \neg G = \emptyset$ for some $F \in \mathcal{F}$, so $F \subseteq G$ and $G \in \mathcal{F}$, a contradiction. Then there must be an ultrafilter $\mathcal{U}$ extending $\mathcal{F}(\neg G)$. Since $\neg G$ is an element of $\mathcal{U} \supseteq \mathcal{F}(\neg G)$, by the remark following Corollary II.9, $\omega_\mathcal{U} \in \mathfrak{b}^{-1}(\neg G) = \Lambda \setminus \mathfrak{b}^{-1}(G)$. But by definition of $\mathcal{F}_{\Omega_\mathcal{F}}$, $\mathfrak{b}^{-1}(G)$ is a clopen neighborhood of $\Omega_{\mathcal{F}}$, so $\omega_{\mathcal{U}} \notin \Omega_{\mathcal{F}}$. On the other hand, as $\mathcal{F} \subseteq \mathcal{F}(\neg G) \subseteq \mathcal{U}$, for any $F \in \mathcal{F}$ we similarly have $\omega_{\mathcal{U}} \in \mathfrak{b}^{-1}(F)$, so $\omega_{\mathcal{U}} \in \Omega_{\mathcal{F}}$, a contradiction! This shows $\mathcal{F}_{\Omega_{\mathcal{F}}} = \mathcal{F}$.

Now, if $\Omega_{\mathcal{F}} = \emptyset$, then $\mathcal{F} = \mathcal{F}_{\Omega_{\mathcal{F}}} \ni \mathfrak{b}(\emptyset) = \emptyset$. Conversely, if $\emptyset \in \mathcal{F}$, then $\Omega_{\mathcal{F}} \subseteq \mathfrak{b}^{-1}(\emptyset) = \emptyset$. If $\mathcal{F}$ is an ultrafilter, then $\omega \in \Omega_{\mathcal{F}}$ iff $\omega \in \mathfrak{b}^{-1}(F)$ for all $F \in \mathcal{F}$, iff $F \in \mathcal{U}_\omega$ for all $F \in \mathcal{F}$, iff $\mathcal{F} \subseteq \mathcal{U}_\omega$. Since both $\mathcal{F}$ and $\mathcal{U}_\omega$ are ultrafilters, this happens iff $\mathcal{F} = \mathcal{U}_\omega$. So $\Omega_{\mathcal{F}} = \{\omega_{\mathcal{F}}\}$. Conversely, if $\Omega_{\mathcal{F}}$ is a singleton, then clopen neighborhoods of $\Omega_{\mathcal{F}}$ is easily seen to be an ultrafilter in the Boolean algebra of clopen subsets of $\Lambda$. As $\mathcal{F} = \mathcal{F}_{\Omega_{\mathcal{F}}}$ is the image of this ultrafilter under $\mathfrak{b}$, a Boolean algebra isomorphism, $\mathcal{F}$ is an ultrafilter. Finally, if $\mathcal{F}$ is principal at $A \in \mathfrak{B}\big((\Omega, \mu)\big)$, then for all $F \in \mathcal{F}$, $\mathfrak{b}^{-1}(A) \subseteq \mathfrak{b}^{-1}(F)$. We also have $A \in \mathcal{F}$, so $\Omega_{\mathcal{F}} = \cap\{\mathfrak{b}^{-1}(F): F \in \mathcal{F}\} = \mathfrak{b}^{-1}(A)$. Conversely, if $\Omega_{\mathcal{F}}$ is clopen, then $\mathcal{F} = \mathcal{F}_{\Omega_{\mathcal{F}}}$ consists of the images under $\mathfrak{b}$ of all clopen sets containing $\Omega_{\mathcal{F}}$, which are exactly all measurable sets majorizing $\mathfrak{b}(\Omega_{\mathcal{F}})$. ∎

**Remark:** The one-to-one correspondence above is order-reversing. In particular, for an ultrafilter $\mathcal{U}$ and a filter $\mathcal{F}$ on $\mathfrak{B}\big((\Omega, \mu)\big)$, $\mathcal{F} \subseteq \mathcal{U}$ iff $\omega_{\mathcal{U}} \in \Omega_{\mathcal{F}}$. Also, as the intersection of two closed sets is closed, correspondingly two filters $\mathcal{F}$ and $\mathcal{G}$ have a join (supremum), $\mathcal{F} \vee \mathcal{G} = \{F \cap G: F \in \mathcal{F}, G \in \mathcal{G}\}$, with $\Omega_{\mathcal{F} \vee \mathcal{G}} = \Omega_{\mathcal{F}} \cap \Omega_{\mathcal{G}}$. In case $\mathcal{G}$ is principal at $A \in \mathfrak{B}\big((\Omega, \mu)\big)$, we let $\mathcal{F}(A) = \mathcal{F} \vee \mathcal{G} = \{X \in \mathfrak{B}\big((\Omega, \mu)\big): \exists F \in \mathcal{F} \text{ s.t. } F \cap A \subseteq X\}$. The union of two closed sets is closed. So the two filters also have a meet (infimum), $\mathcal{F} \wedge \mathcal{G} = \mathcal{F} \cap \mathcal{G}$, with $\Omega_{\mathcal{F} \wedge \mathcal{G}} = \Omega_{\mathcal{F}} \cup \Omega_{\mathcal{G}}$.

We also note that, just like Corollary II.9, if we regard $\Lambda$ as the Stone space of $\mathfrak{B}\big((\Omega, \mu)\big)$, this is simply a special case of the Stone representation theorem (see [BS81, Exercise IV.4.7]). ∎

Another tool that will be useful as well as be of independent interest is the notion of a lifting map of $L^\infty(\Omega, \mu)$, defined as follows,

**Definition:** Let $M^\infty(\Omega, \mu)$ be the $C^*$-algebra of all bounded measurable functions from $(\Omega, \mu)$ to $\mathbb{C}$, equipped with the supremum norm. (Note that as opposed to in $L^\infty(\Omega, \mu)$, we are not identifying functions that coincide a.e., and correspondingly the norm is the supremum norm instead of the essential supremum norm.) Let $\pi: M^\infty(\Omega, \mu) \to L^\infty(\Omega, \mu)$ be the canonical quotient map. A *lifting map* $\rho: L^\infty(\Omega, \mu) \to M^\infty(\Omega, \mu)$ is a *-homomorphism s.t. $\pi \circ \rho = Id_{L^\infty(\Omega, \mu)}$.

**Theorem 2:** There exists a lifting map $\rho: L^\infty(\Omega, \mu) \to M^\infty(\Omega, \mu)$.

**Proof:** This is [II69, Theorem IV.3]. The assumption needed there is that $(\Omega, \mu)$ is strictly localizable, which is defined in [II69, §I.8]. It is not hard to verify using results listed there and in previous sections of [II69] that this condition is equivalent to $(\Omega, \mu)$ being semifinite and decomposable. ∎

**Remark:** We observe that, for any lifting map $\rho: L^\infty(\Omega, \mu) \to M^\infty(\Omega, \mu)$, as it is a *-homomorphism, $\rho$ followed by pointwise evaluation at $\omega \in \Omega$ is a character of $L^\infty(\Omega, \mu)$, which we denote by $\phi_\rho(\omega)$. This defines a map $\phi_\rho: \Omega \to \Lambda$. ∎

We now specialize to the case where $\Omega$ is a compact Hausdorff space and $\mu$ is a faithful Radon probability measure. Then $\mu$ induces a faithful normal Radon probability measure on $\Lambda$, which we shall denote by $\mu_\Lambda$. (These assumptions shall be in effect for the rest of this section.) We note that, as $C(\Omega) \subseteq L^\infty(\Omega, \mu)$, characters on $L^\infty(\Omega, \mu)$ restricts to characters on $C(\Omega)$, which are naturally in bijection with $\Omega$. Hence, $\Lambda$ is partitioned into collections of characters that extend evaluation at a fixed $\omega \in \Omega$. We shall now produce an ultrafilter version of this partition.

**Definition:** For any $\omega \in \Omega$, we define a filter $\mathcal{D}_\omega$ on $\mathfrak{B}\big((\Omega, \mu)\big)$ by,

$$\mathcal{D}_\omega = \{X \in \mathfrak{B}((\Omega, \mu)): X \text{ almost contains an open neighborhood of } \omega\}$$

Since $\mu$ has full support, this is a proper filter. As $\Omega$ is Hausdorff, it is easy to see that $\mathcal{D}_{\omega_1} \vee \mathcal{D}_{\omega_2} = \mathfrak{B}((\Omega, \mu))$ whenever $\omega_1 \neq \omega_2$, so no ultrafilter can contain both $\mathcal{D}_{\omega_1}$ and $\mathcal{D}_{\omega_2}$. Hence $\Omega_{\mathcal{D}_\omega}$'s are disjoint closed subsets of $\Lambda$. We shall now see that they actually form a partition,

**Proposition 3:** Given any ultrafilter $\mathcal{U}$ on $\mathfrak{B}((\Omega, \mu))$, it must contain $\mathcal{D}_\omega$ for some $\omega \in \Omega$.

**Proof:** Let $\mathbb{K}$ be the collection of all compact sets representing any element of $\mathcal{U}$. Since $\Omega \in \mathcal{U}$, it is nonempty. Any finite intersection of elements of $\mathbb{K}$ represents an element of $\mathcal{U}$ and therefore has positive measure and is nonempty, so compactness implies $\cap \mathbb{K}$ is nonempty. Let $\omega \in \cap \mathbb{K}$. We claim that for all open neighborhood $U$ of $\omega$, $U \in \mathcal{U}$, i.e., $\mathcal{D}_\omega \subseteq \mathcal{U}$. Indeed, assume otherwise, then $\Omega \setminus U \in \mathcal{U}$ for some $U$, but $\Omega \setminus U$ is compact, so it is in $\mathbb{K}$. Therefore $\omega \notin \cap \mathbb{K}$, a contradiction! ∎

We now demonstrate that ultrafilters in $\mathcal{D}_\omega$ do extend limits approaching $\omega$. To prove a stronger version of this that will be useful in later sections, we shall introduce the following concept,

**Definition:** Let $M$ be a von Neumann algebra, $a \in L^\infty(\Omega, \mu) \overline{\otimes} M$, $\varphi \in L^\infty(\Omega, \mu)^*$. Then $\varphi \otimes Id(a) \in M = (M_*)^*$ is given by $(\varphi \otimes Id(a))(\psi) = \varphi(Id \otimes \psi(a))$ for all $\psi \in M_*$. Note that this definition is consistent with the previous definition when $\varphi \in L^\infty(\Omega, \mu)_*$. If $\mathcal{U}$ is an ultrafilter on $\mathfrak{B}((\Omega, \mu))$, then we let $\lim_\mathcal{U} a = \omega_\mathcal{U} \otimes Id(a)$.

**Remark:** This definition applies, in general, when $(\Omega, \Sigma, \mu)$ is a semifinite, decomposable measure space, which will be the context in which this definition is applied in later sections. For this section, we shall continue to assume $\Omega$ is a compact Hausdorff space and $\mu$ is a faithful Radon probability measure. ∎

**Theorem 4:** Let $a \in L^\infty(\Omega, \mu) \overline{\otimes} M$ be given by a bounded strongly measurable function $f: (\Omega, \mu) \to M$. If for some $\omega_0 \in \Omega$, $\lim_{\omega \to \omega_0} f(\omega)$ exists where the limit is in the weak* topology, then for any ultrafilter $\mathcal{U}$ containing $\mathcal{D}_{\omega_0}$, $\lim_\mathcal{U} a = \lim_{\omega \to \omega_0} f(\omega)$.

**Proof:** Fix any $\varphi \in M_*$, we first show that $Id \otimes \varphi(a) = \varphi \circ f$. We only need to consider the case where $\varphi = |k\rangle\langle h|$ for some $h, k \in H$, since taking linear combinations and limits then yields the general result. Then, for any measurable $F \subseteq \Omega$ with positive measure,

$$1_F(Id \otimes \varphi(a)) = \langle 1 \otimes h, a(1_F \otimes k)\rangle = \int_F \langle h, f(\omega)k\rangle \, d\mu = \int_F \varphi \circ f(\omega) \, d\mu = 1_F(\varphi \circ f)$$

This proves the claim. We observe that $\varphi\left(\lim_{\omega \to \omega_0} f(\omega)\right) = \lim_{\omega \to \omega_0} \varphi \circ f(\omega)$. Thus, for any $\varepsilon > 0$, there exists an open neighborhood $U$ of $\omega_0$ s.t. for all $\omega \in U$, $\left|\varphi \circ f(\omega) - \varphi\left(\lim_{\omega \to \omega_0} f(\omega)\right)\right| < \varepsilon$. It is then easy to see that $\varphi\left(\lim_\mathcal{U} a\right) = \lim_\mathcal{U}(Id \otimes \varphi(a)) = \lim_\mathcal{U}(\varphi \circ f) = \varphi\left(\lim_{\omega \to \omega_0} f(\omega)\right)$. Since $\varphi \in M_*$ is arbitrary, $\lim_\mathcal{U} a = \lim_{\omega \to \omega_0} f(\omega)$. ∎

**Remark:** When $M = \mathbb{C}$ and $a \in C(\Omega)$, the above demonstrates that for any ultrafilter $\mathcal{U}$ containing $\mathcal{D}_\omega$, $\lim_\mathcal{U} a = a(\omega)$. Hence, a character is in $\Omega_{\mathcal{D}_\omega}$ iff it extends evaluation at $\omega$ on $C(\Omega)$. ∎

We shall now use this result to obtain a characterization of $\phi_\rho: \Omega \to \Lambda$ when $\rho: L^\infty(\Omega, \mu) \to M^\infty(\Omega, \mu)$ is a lifting map. To do so, we need one extra assumption, namely $C(\Omega)$ is norm-separable. (This assumption shall be in effect for the rest of this section.) This includes the most interesting case $(\Omega, \mu) = ([0, 1], \lambda)$.

We first observe that $M^\infty(\Omega, \mu) \subseteq l^\infty(\Omega)$ (where $l^\infty(\Omega)$ is the algebra of all bounded functions from $\Omega$ to $\mathbb{C}$). It is easy to see that the definition of $\phi_\rho: \Omega \to \Lambda$ apply to the more general case where $\rho$ is a *-homomorphism $L^\infty(\Omega, \mu) \to l^\infty(\Omega)$. In fact, we see that such *-homomorphisms are in bijection with maps $\phi: \Omega \to \Lambda$.

**Remark:** Before we prove this characterization result, we need to recall certain facts regarding Baire sets. Recall that for a locally compact Hausdorff space $X$, the $\sigma$-algebra of Baire sets is the smallest $\sigma$-algebra that makes all compactly supported continuous functions on $X$ measurable. Equivalently, it is the $\sigma$-algebra generated by all compact $G_\delta$ sets. (See [Fol99, Exercise 7.1.4].) When we say a function $\phi: (\Omega, \mu) \to X$ is *Baire-measurable*, we shall mean it is a measurable map when the former is equipped with the $\sigma$-algebra of all $\mu$-measurable map and the latter is equipped with the $\sigma$-algebra of Baire sets. We shall further observe that, on a hyperstonean space such as $\Lambda$, as the algebra of continuous functions is a von Neumann algebra, all continuous functions are norm-limits of linear combinations of projections. Since projections in $C(\Lambda)$ are indicator functions of clopen sets, we see that the Baire $\sigma$-algebra is generated by clopen subsets of $\Lambda$. ∎

**Theorem 5:** Let $\rho: L^\infty(\Omega, \mu) \to l^\infty(\Omega)$ be a *-homomorphism. Then the following are equivalent,

1. $\rho$ has range in $M^\infty(\Omega, \mu)$ and is a lifting map $L^\infty(\Omega, \mu) \to M^\infty(\Omega, \mu)$;
2. $\phi_\rho: \Omega \to \Lambda$ satisfies all of the three conditions below:
   A. $\phi_\rho(\omega) \in \Omega_{\mathcal{D}_\omega}$ for a.e.-$\omega$;
   B. $\phi_\rho: (\Omega, \mu) \to \Lambda$ is Baire-measurable;
   C. For any compact $G_\delta$ set $K \subseteq \Lambda$ with no interior, $\phi_\rho^{-1}(K)$ is null.

**Proof:** $(2 \Rightarrow 1)$ We first observe that, regarding elements of $L^\infty(\Omega, \mu)$ as continuous functions on $\Lambda$, then $\rho(f) = f \circ \phi_\rho$. As all continuous functions are measurable under the Baire $\sigma$-algebra, we see that condition B implies $\rho(f)$ is always measurable, i.e., $\rho$ has range in $M^\infty(\Omega, \mu)$. Let $\pi: M^\infty(\Omega, \mu) \to L^\infty(\Omega, \mu)$ be the canonical quotient map. We need to show $\pi \circ \rho = Id_{L^\infty(\Omega, \mu)}$.

Let $E = \{\omega \in \Omega: \phi_\rho(\omega) \in \Omega_{\mathcal{D}_\omega}\}$. Then $E$ is co-null by condition A. For any $f \in C(\Omega) \subseteq L^\infty(\Omega, \mu)$, $\omega \in E$, by Theorem 4, $\rho(f)(\omega) = f(\phi_\rho(\omega)) = f(\omega)$. Hence, for all continuous functions $\rho(f)(\omega) = f(\omega)$ a.e., so $\pi \circ \rho|_{C(\Omega)} = Id_{C(\Omega)}$. We observe that, by Lusin's Theorem [Fol99, Theorem 7.10] and Dominated Convergence Theorem, we can prove $C(\Omega)$ is weak*-dense in $L^\infty(\Omega, \mu)$, so it suffices to show $\pi \circ \rho$ is normal.

To do so, we first observe that, for any closed subset $K$ of $\Lambda$ with no interior, there is a compact $G_\delta$ set $K'$ with no interior that contains $K$. Indeed, by [Tak79, Proposition III.1.11], $\mu_\Lambda(K) = 0$. By outer regularity, there is a sequence of open sets $O_n \supseteq K$ s.t. $\mu_\Lambda(O_n) \to 0$. We observe that $\overline{O_n}$ is clopen and $\mu_\Lambda(\overline{O_n}) = \mu_\Lambda(O_n) \to 0$, so $K'$, the intersection of all $\overline{O_n}$, is a compact $G_\delta$ set containing $K$ that is $\mu_\Lambda$-null. Since $\mu_\Lambda$ is faithful, $K'$ has no interior.

Let $\mu' = (\phi_\rho)_*(\mu)$ be the push-forward of $\mu$ by $\phi_\rho$, which is defined on all Baire sets. Being a finite Baire measure defined on a compact Hausdorff space, it is well-known that it extends uniquely to a Radon measure, which we shall also denote by $\mu'$. (See, for example, [Hal74, Theorem 54.D].) By the result in last paragraph and condition C, we see that all nowhere dense set is $\mu'$-null. Again by [Tak79, Proposition III.1.11], $\mu'$ is normal. Regarding $\mu'$ as a linear functional on $C(\Lambda) = L^\infty(\Omega, \mu)$ and $\mu$ as a linear functional on $L^\infty(\Omega, \mu)$, we see that $\mu' = \mu \circ \pi \circ \rho$. Since $\mu$ is a normal faithful linear functional, this is sufficient in showing $\pi \circ \rho$ is normal.

(1 ⇒ 2) By definition of $\rho$, for any fixed $f \in C(\Omega)$ we have $\phi_\rho(\omega)(f) = f(\omega)$ a.e. Since $C(\Omega)$ is norm-separable, there is a shared co-null set $E \subseteq \Omega$ s.t. $\phi_\rho(\omega)(f) = f(\omega)$ for all $\omega \in E$ and $f \in C(\Omega)$. The remark after Theorem 4 then implies $\phi_\rho(\omega) \in \Omega_{\mathcal{D}_\omega}$ for all $\omega \in E$, whence condition A follows.

To prove condition B, as the Baire $\sigma$-algebra is generated by clopen sets, it suffices to show $\phi_\rho^{-1}(K)$ is measurable for all clopen $K \subseteq \Lambda$. But it is easy to see that $\phi_\rho^{-1}(K)$ is the set on which $\rho(1_K)$ takes value 1. As $\rho(1_K) \in M^\infty(\Omega, \mu)$, the result follows.

We further observe that, as $\rho$ is a lifting map, for a clopen $K \subseteq \Lambda$, $\phi_\rho^{-1}(K) \subseteq \Omega$ is not only measurable but equals the measurable set represented by $K$, up to a null set. Hence, $\mu\left(\phi_\rho^{-1}(K)\right) = \mu_\Lambda(K)$. As has been seen in the proof of the 2 ⇒ 1 direction, any compact $K \subseteq \Lambda$ with no interior is contained in $\cap K_n$ where $K_n$ is a sequence of clopen sets whose measures converge to 0. Thus, $\mu\left(\phi_\rho^{-1}(K)\right) \leq \mu\left(\phi_\rho^{-1}(K_n)\right) = \mu_\Lambda(K_n) \to 0$, so $\phi_\rho^{-1}(K)$ is null, which is condition C. ∎

**Remark:** Since a lifting map $\rho: L^\infty(\Omega, \mu) \to M^\infty(\Omega, \mu)$ exists, by the above theorem we may, if necessary, alter $\phi_\rho$ on a null set so that $\phi_\rho(\omega) \in \Omega_{\mathcal{D}_\omega}$ for all $\omega$. That is, there exists a lifting map $\rho: L^\infty(\Omega, \mu) \to M^\infty(\Omega, \mu)$ s.t. $\phi_\rho(\omega) \in \Omega_{\mathcal{D}_\omega}$ for all $\omega$. We shall call such lifting maps *standard*. ∎

**Corollary 6:** If $\rho: L^\infty(\Omega, \mu) \to M^\infty(\Omega, \mu)$ is a lifting map, then $\phi_\rho: (\Omega, \mu) \to \Lambda$ is measurable with $\Lambda$ equipped with the $\sigma$-algebra of $\mu_\Lambda$-measurable sets. Furthermore, $\mu\left(\phi_\rho^{-1}(E)\right) = \mu_\Lambda(E)$ for all measurable $E \subseteq \Lambda$.

**Proof:** By [Tak79, Corollary III.1.13], every measurable $E \subseteq \Lambda$ differs from a clopen set by a $\mu_\Lambda$-null set. As we have seen in the proof of the 1 ⇒ 2 direction of Theorem 5, $\mu\left(\phi_\rho^{-1}(E)\right) = \mu_\Lambda(E)$ holds for all clopen $E \subseteq \Lambda$. Hence, it suffices to show $\phi_\rho^{-1}(E)$ is null for all $\mu_\Lambda$-null set $E$. Again by [Tak79, Corollary III.1.13], $\bar{E}$ is null and therefore has no interior. We have seen in the proof of the 1 ⇒ 2 direction of Theorem 5 that $\phi_\rho^{-1}(\bar{E})$ is null, whence the result follows. ∎

**Remark:** If $(\Omega, \mu)$ is a finite Radon measure on a locally compact (but not compact) Hausdorff space with full support, then we can simply extend it to the one-point compactification by letting $\mu(\{\infty\}) = 0$. By renormalization, we may assume $\mu$ is a probability measure. This doesn't change $L^\infty(\Omega, \mu)$ in any way and we can then apply all the above results. (Theorem 5 and Corollary 6 require $C(\Omega \cup \{\infty\})$ to be norm-separable. Since $C(\Omega \cup \{\infty\}) = C_0(\Omega) + \mathbb{C}$, this is equivalent to $C_0(\Omega)$ being norm-separable.) For example, if $\Omega = \mathbb{N}$ and $\mu$ is a faithful probability measure on $\mathbb{N}$, then $\mathcal{D}_a$ is the principal ultrafilter at $a$ when $a \in \mathbb{N}$. $\mathcal{D}_\infty$ is the Fréchet filter (i.e., the filter consisting of all cofinite subsets) on $\mathbb{N}$, so $\Omega_{\mathcal{D}_\infty}$ consists of all free ultrafilters on $\mathbb{N}$. ∎

Condition A in Theorem 5 is necessary. Indeed, let $(\Omega, \mu) = ([0, 1], \lambda)$. Fix a standard lifting map $\rho: L^\infty([0, 1], \lambda) \to M^\infty([0, 1], \lambda)$. Then define $\varphi: M^\infty([0, 1], \lambda) \to M^\infty([0, 1], \lambda)$ by $\varphi(f)(t) = f(1 - t)$. We then have $\phi_{\varphi \circ \rho}(t) \in \Omega_{\mathcal{D}_{1-t}}$, so $\phi_{\varphi \circ \rho}(t) \in \Omega_{\mathcal{D}_t}$ only when $t = \frac{1}{2}$ and condition A fails for $\varphi \circ \rho$. But as

$\varphi$ corresponds to a measure-preserving automorphism of $([0,1], \lambda)$, we see that conditions B and C are satisfied by $\varphi \circ \rho$, so condition A is necessary.

Moreover, condition A does not imply either condition B or condition C, as we shall see later. Neither does condition B imply condition C. Indeed, suppose $(\Omega, \mu)$ is atomless, then there are no isolated point in $\Lambda$. Hence, it is easy to verify any constant map $\phi: (\Omega, \mu) \to \Lambda$ satisfies condition B but not C. Whether condition C is necessary in general is unknown, though in some cases it easily follows from condition A. For example, when $\Omega = \mathbb{N}$ and $\mu$ is a faithful probability measure on $\mathbb{N}$. (More precisely, $\Omega$ is the one-point compactification of $\mathbb{N}$ and $\mu$ is extended by letting $\mu(\{\infty\}) = 0$.) Readers may easily verify in that case condition B is automatic and condition A implies condition C.

To further analyze the structure of $\Lambda$, we need a few definitions,

**Definition:** Let $F \subseteq \Omega$ be a measurable subset. We define the *proper interior* of $F$, denoted by $^p\mathbf{int}(F)$ to be $^p\mathbf{int}(F) = \{\omega \in \Omega: \exists U \subseteq \Omega, U \text{ is an open neighborhood of } \omega, U \subseteq F \text{ a.e.}\}$. It is easy to see that $^p\mathbf{int}(F)$ is the union of all open sets $U \subseteq \Omega$ which are almost contained in $F$ and is therefore open. We define the *weak extension* of $F$, denoted by $^w\mathbf{ext}(F)$, to be $^w\mathbf{ext}(F) = \bigcup_{\omega \in F} \Omega_{\mathcal{D}_\omega} \subseteq \Lambda$. We define the *strong extension* of $F$, denoted by $^s\mathbf{ext}(F)$, to be $^s\mathbf{ext}(F) = \bigcup \{\Omega_{\mathcal{D}_\omega}: \Omega_{\mathcal{D}_\omega} \cap \mathfrak{b}^{-1}(F) \neq \emptyset\} \subseteq \Lambda$. We see that $^p\mathbf{int}(F)$ and $^s\mathbf{ext}(F)$ do not change while $^w\mathbf{ext}(F)$ does when $F$ is altered by a null set. For a subset $K \subseteq \Lambda$, we also define its *intension*, denoted by $\mathbf{ints}(K)$, to be $\mathbf{ints}(K) = \{\omega \in \Omega: \Omega_{\mathcal{D}_\omega} \cap K \neq \emptyset\}$.

**Lemma 7:** Let $F \subseteq \Omega$ be measurable, $K \subseteq \Lambda$ be a clopen subset,

1. $^w\mathbf{ext}(^p\mathbf{int}(\Omega \setminus F)) = \Lambda \setminus {}^s\mathbf{ext}(F)$ for all measurable $F \subseteq \Omega$;
2. $\mathbf{ints}(K) = \Omega \setminus {}^p\mathbf{int}(\Omega \setminus \mathfrak{b}^{-1}(K))$, so $\mathbf{ints}(K)$ is closed and we have, if $F = \mathfrak{b}(K)$, then $^s\mathbf{ext}(F) = {}^w\mathbf{ext}(\mathbf{ints}(K))$.

**Proof:**

1. We observe that, since $\Omega_{\mathcal{D}_\omega}$ form a partition of $\Lambda$, $\Lambda \setminus {}^s\mathbf{ext}(F) = \bigcup\{\Omega_{\mathcal{D}_\omega}: \Omega_{\mathcal{D}_\omega} \cap \mathfrak{b}^{-1}(F) = \emptyset\} = \bigcup\{\Omega_{\mathcal{D}_\omega}: \Omega_{\mathcal{D}_\omega} \subseteq \mathfrak{b}^{-1}(\Omega \setminus F)\}$. By Theorem 1 and the remark after it, $\Omega_{\mathcal{D}_\omega} \subseteq \mathfrak{b}^{-1}(\Omega \setminus F)$ is equivalent to the principal filter at $\Omega \setminus F$ being contained in $\mathcal{D}_\omega$, or equivalently $\Omega \setminus F$ almost contains an open neighborhood of $\omega$ by the definition of $\mathcal{D}_\omega$. By the definition of $^p\mathbf{int}(\Omega \setminus F)$ this is further equivalent to $\omega \in {}^p\mathbf{int}(\Omega \setminus F)$. Hence, $\Lambda \setminus {}^s\mathbf{ext}(F) = \bigcup\{\Omega_{\mathcal{D}_\omega}: \omega \in {}^p\mathbf{int}(\Omega \setminus F)\} = {}^w\mathbf{ext}({}^p\mathbf{int}(\Omega \setminus F))$.
2. If $F = \mathfrak{b}(K)$, then $^s\mathbf{ext}(F) = \bigcup_{\omega \in \mathbf{ints}(K)} \Omega_{\mathcal{D}_\omega} = \Lambda \setminus {}^w\mathbf{ext}({}^p\mathbf{int}(\Omega \setminus F)) = \Lambda \setminus \left(\bigcup_{\omega \in {}^p\mathbf{int}(\Omega \setminus F)} \Omega_{\mathcal{D}_\omega}\right) = \bigcup_{\omega \in \Omega \setminus {}^p\mathbf{int}(\Omega \setminus F)} \Omega_{\mathcal{D}_\omega}$. Hence $\mathbf{ints}(K) = \Omega \setminus {}^p\mathbf{int}(\Omega \setminus \mathfrak{b}^{-1}(K))$. That $^s\mathbf{ext}(F) = {}^w\mathbf{ext}(\mathbf{ints}(K))$ is then just a restatement of definitions. ∎

**Corollary 8:** If $K \subseteq \Lambda$ is closed, then $\mathbf{ints}(K)$ is closed.

**Proof:** Let $N$ be the collection of clopen neighborhoods of $K$, then $K = \bigcap N$. We claim that $\mathbf{ints}(K) = \bigcap_{F \in N} \mathbf{ints}(F)$, whence the result follows from the proposition above. That $\mathbf{ints}(K) \subseteq \bigcap_{F \in N} \mathbf{ints}(F)$ is obvious. Conversely, suppose $\omega \in \bigcap_{F \in N} \mathbf{ints}(F)$. Then $\Omega_{\mathcal{D}_\omega}$ intersects with $F$ for all $F \in N$. As $K = \bigcap N$, we have $\Omega_{\mathcal{D}_\omega} \cap K = \bigcap_{F \in N}(\Omega_{\mathcal{D}_\omega} \cap F)$. Since $N$ is closed under taking finite intersections, $\Omega_{\mathcal{D}_\omega} \cap F \neq \emptyset$ for all $F \in N$, and $\Omega_{\mathcal{D}_\omega}$ is closed, we see that compactness implies $\Omega_{\mathcal{D}_\omega} \cap K \neq \emptyset$. Thus, $\omega \in \mathbf{ints}(K)$. ∎

**Proposition 9:** If $K \subseteq \Lambda$ is co-null, then $\mathbf{ints}(K)$ is as well.

**Proof:** Fix a standard lifting $\rho: L^\infty(\Omega, \mu) \to M^\infty(\Omega, \mu)$. Then by Corollary 6, $\phi_\rho^{-1}(K)$ is co-null. By the definition of standard lifting it is easy to see that $\phi_\rho^{-1}(K) \subseteq \mathbf{ints}(K)$, so the result follows. ∎

We now turn to our main task of analyzing the case where $(\Omega, \mu) = ([0,1], \lambda)$. For the rest of this section, we shall use $\Lambda$ to denote specifically the spectrum of $L^\infty([0,1], \lambda)$. To ensure notational consistency, we write the induced measure on $\Lambda$ to be $\lambda_\Lambda$.

We first refine our partition of $\Lambda$ into $\Omega_{\mathcal{D}_t}$ by differentiating between left and right limits, as follows,

**Definition:** For any $t \in (0, 1]$, we define $\mathcal{L}_t = \{X \in \mathcal{B}(([0,1], \lambda)): \exists \varepsilon > 0 \text{ s.t. } (t - \varepsilon, t) \cap [0,1] \text{ is almost contained in } X\}$. For any $t \in [0, 1)$, we define $\mathcal{R}_t = \{X \in \mathcal{B}(([0,1], \lambda)): \exists \varepsilon > 0 \text{ s.t. } (t, t + \varepsilon) \cap [0,1] \text{ is almost contained in } X\}$. They are proper filters extending $\mathcal{D}_t$. We observe that $\mathcal{D}_0 = \mathcal{R}_0$ and $\mathcal{D}_1 = \mathcal{L}_1$.

**Proposition 10:** For any fixed $0 < a < 1$, if $\mathcal{U}$ contains $\mathcal{D}_a$, then it must contain exactly one of $\mathcal{L}_a$ or $\mathcal{R}_a$.

**Proof:** We observe that, since $(a - \varepsilon, a) \cap (a, a + \varepsilon) = \emptyset$ for any $\varepsilon > 0$, $\mathcal{L}_a \vee \mathcal{R}_a = \mathcal{B}(([0,1], \lambda))$, so no ultrafilter can contain both. Now, let $\mathcal{U}$ contains $\mathcal{D}_a$ but not $\mathcal{R}_a$, it suffices to show $\mathcal{L}_a \subseteq \mathcal{U}$. Since $\mathcal{U}$ does not contain $\mathcal{R}_a$, there exists $\varepsilon_0 > 0$ s.t. $(a, a + \varepsilon_0) \cap [0,1] \notin \mathcal{U}$. Then for any $\varepsilon \leq \varepsilon_0$, $(a, a + \varepsilon) \cap [0,1] \notin \mathcal{U}$, so $[0,1] \setminus (a, a + \varepsilon) \in \mathcal{U}$. Since $\mathcal{D}_a \subseteq \mathcal{U}$, $(a - \varepsilon, a + \varepsilon) \cap [0,1] \in \mathcal{U}$, so their intersection, which up to a null set (namely $\{a\}$) is $(a - \varepsilon, a) \cap [0,1]$, belongs to $\mathcal{U}$. As this holds for all $\varepsilon \leq \varepsilon_0$, we have $\mathcal{L}_a \subseteq \mathcal{U}$. ∎

Hence, for any ultrafilter $\mathcal{U}$ on $\mathcal{B}(([0,1], \lambda))$, it contains exactly one of $\mathcal{L}_a$, $0 < a \leq 1$, or $\mathcal{R}_a$, $0 \leq a < 1$. By essentially the same proof as Theorem 4, we may refine it as follows,

**Theorem 11:** Let $a \in L^\infty([0,1], \lambda) \overline{\otimes} M$ be given by a bounded strongly measurable function $f: ([0,1], \lambda) \to M$. If for some $0 < t \leq 1$, $\lim_{s \uparrow t} f(s)$ exists where the limit is in the weak* topology, then for any ultrafilter $\mathcal{U}$ containing $\mathcal{L}_a$, $\lim_\mathcal{U} a = \lim_{s \uparrow t} f(s)$. If for some $0 \leq t < 1$, $\lim_{s \downarrow t} f(s)$ exists where the limit is in the weak* topology, then for any ultrafilter $\mathcal{U}$ containing $\mathcal{R}_a$, $\lim_\mathcal{U} a = \lim_{s \downarrow t} f(s)$.

**Example:** As an example demonstrating the usage of this theorem, we shall consider $f \in L^\infty([0,1], \lambda) \overline{\otimes} l^\infty((0,1])$ defined in Lemma I.5. We calculate the collection of all $\lim_\mathcal{U} f$ for ultrafilters $\mathcal{U}$. For simplicity, we shall denote the projection onto $X \subseteq (0,1]$ in $l^\infty((0,1])$ by $E_X$. Then as a strongly measurable function $([0,1], \lambda) \to l^\infty((0,1])$, $f$ is defined by $f(s) = E_{(s, 1]}$. We then see that, for any $0 < t \leq 1$, $\lim_{s \uparrow t} f(s) = E_{[t, 1]}$ weak*, and for any $0 \leq t < 1$, $\lim_{s \downarrow t} f(s) = E_{(t, 1]}$. Thus, Theorem 11 implies the collection of all $\lim_\mathcal{U} f$ is $\{E_{[t, 1]}: 0 < t \leq 1\} \cup \{E_{(t, 1]}: 0 \leq t < 1\}$.

We note that none of $\mathcal{L}_a$ or $\mathcal{R}_a$ is an ultrafilter itself. To show this, we first note that these filters are all related to each other via Borel automorphisms of $([0,1], \lambda)$. Indeed, the translation automorphisms sending $a \in [0, 1)$ to $a + x \pmod 1 \in [0, 1)$ for $x \in \mathbb{R}$ send $\mathcal{L}_a$ to each other and $\mathcal{R}_a$ to each other. The reflection automorphism sending $a \in [0, 1]$ to $1 - a \in [0, 1]$ sends $\mathcal{L}_a$ to $\mathcal{R}_{1-a}$ and $\mathcal{R}_a$ to $\mathcal{L}_{1-a}$. Hence, we only need to show any one of these is not an ultrafilter. We need the following lemma:

**Lemma 12:** Given any compact $F \subseteq [0, 1]$ with positive Lebesgue measure, there exists $a \in F$ s.t. for all $\varepsilon > 0$, $(a - \varepsilon, a + \varepsilon) \cap F$ has positive measure, i.e., $\mathcal{D}_a(F)$ is proper.

**Proof:** Assume the result does not hold, then for all $a \in F$, there exists $\varepsilon_a > 0$ s.t. $(a - \varepsilon_a, a + \varepsilon_a) \cap F$ is null. $(a - \varepsilon_a, a + \varepsilon_a)$ for all $a \in F$ forms an open cover of $F$, so there exists $\{a_i\}_{i=1}^n \subseteq F$ s.t. $F \subseteq \bigcup_{i=1}^n (a_i - \varepsilon_{a_i}, a_i + \varepsilon_{a_i})$. Thus,

$$\lambda(F) = \lambda\left(F \cap \left(\bigcup_{i=1}^n (a_i - \varepsilon_{a_i}, a_i + \varepsilon_{a_i})\right)\right) \leq \sum_{i=1}^n \lambda\left(F \cap (a_i - \varepsilon_{a_i}, a_i + \varepsilon_{a_i})\right) = 0$$

But $F$ is assumed to have positive measure, a contradiction! ∎

Recall that it is possible to have a compact $F \subseteq [0, 1]$ with no interior but positive measure. Indeed, fix a bijection $k: \mathbb{Q} \cap [0, 1] \to \mathbb{N}$, let,

$$O = \bigcup_{q \in \mathbb{Q} \cap [0,1]} (q - 2^{-k(q)-2}, q + 2^{-k(q)-2}) \cap [0, 1]$$

Then $O$ is an open dense subset of $[0, 1]$ with measure smaller than or equal to $\frac{1}{2}$, so $F = [0, 1] \setminus O$ is compact, nowhere dense, but has measure at least $\frac{1}{2}$. By Lemma 12, there exists $a \in F$ with $\mathcal{D}_a(F)$ being proper, so there is an ultrafilter $\mathcal{U}$ extending it. If either $\mathcal{L}_a$ or $\mathcal{R}_a$ is an ultrafilter (and therefore both of them are), then by Proposition 10, $\mathcal{U}$ is either $\mathcal{L}_a$ or $\mathcal{R}_a$, so $F \in \mathcal{L}_a$ or $F \in \mathcal{R}_a$. This is however impossible. Indeed, if either of these is true, then $F$ contains a nonempty open subset $U \subseteq [0, 1]$ up to a null set (either $(a - \varepsilon, a) \cap [0, 1]$ or $(a, a + \varepsilon) \cap [0, 1]$ for some $\varepsilon > 0$). That is, $\lambda(U \setminus F) = \lambda(U \cap O) = 0$. Since $U \cap O$ is open, this is only possible if $U \cap O = \emptyset$, so $U \subseteq F$, which contradicts the assumption that $F$ has no interior. So we have,

**Theorem 13:** None of $\mathcal{L}_a$, $0 < a \leq 1$, or $\mathcal{R}_a$, $0 \leq a < 1$, is an ultrafilter.

Another way of proving this is to consider the cardinality of $\Lambda$. If any of the $\mathcal{L}_a$ or $\mathcal{R}_a$ is an ultrafilter, they all are, so by Proposition 10, $\Lambda$ has the cardinality of continuum, $2^{\aleph_0}$. However, we have,

**Theorem 14:** $|\Lambda| = 2^{2^{\aleph_0}}$.

**Proof:** We first show that $|\Lambda| \leq 2^{2^{\aleph_0}}$. This follows from the following observation: any Lebesgue measurable set equals a Borel set up to a null set, so any element of $\mathfrak{B}(([0, 1], \lambda))$ can be represented by a Borel set. By [Fol99, remark after Proposition 1.23], as the Borel $\sigma$-algebra is generated by intervals, there are exactly continuum many Borel sets, so $|\mathfrak{B}(([0, 1], \lambda))| = 2^{\aleph_0}$. As ultrafilters are subsets of $\mathfrak{B}(([0, 1], \lambda))$, we see that $|\Lambda| \leq 2^{2^{\aleph_0}}$.

For the other direction, we show that there is an injection $\beta \mathbb{N} \setminus \mathbb{N} \hookrightarrow \Lambda$. Since by [BS74, Theorem 6.1.5], $|\beta \mathbb{N} \setminus \mathbb{N}| = 2^{2^{\aleph_0}}$, this yields the desired result. Now, to construct this injection, we divide $[0, 1]$ into countably infinite many disjoint measurable sets $\{X_i\}_{i=1}^\infty$, each of positive measure. Then $L^\infty([0, 1], \lambda) = \prod_{i=1}^\infty L^\infty(X_i, \lambda)$. For each $i$, fix a character $\omega_i$ on $L^\infty(X_i, \lambda)$. Let $\omega \in \beta \mathbb{N} \setminus \mathbb{N} \subseteq l^\infty(\mathbb{N})^*$. We define,

$$\phi_\omega((a_i)_{i=1}^\infty) = \omega\left((\omega_i(a_i))_{i=1}^\infty\right), (a_i)_{i=1}^\infty \in \prod_{i=1}^\infty L^\infty(X_i, \lambda)$$

It is easy to see that $\phi_\omega$ is a character on $\prod_{i=1}^\infty L^\infty(X_i, \lambda) = L^\infty([0,1], \lambda)$. Furthermore, let $\pi: l^\infty(\mathbb{N}) \to \prod_{i=1}^\infty L^\infty(X_i, \lambda)$ be defined by sending $(c_i)_{i=1}^\infty \in l^\infty(\mathbb{N})$ to $(c_i 1)_{i=1}^\infty \in \prod_{i=1}^\infty L^\infty(X_i, \lambda)$, then $\omega = \phi_\omega \circ \pi$. Thus, the map sending $\omega$ to $\phi_\omega$ is injective, concluding the proof. ∎

Proposition 10 means that $\Lambda$ is divided into continuum many pairwise disjoint closed sets, $\Omega_{\mathcal{L}_a}$, $0 < a \leq 1$, and $\Omega_{\mathcal{R}_a}$, $0 \leq a < 1$. Since automorphisms bring the corresponding filters to one another, these closed sets are homeomorphic to each other. In particular, in light of Theorem 14, all of them have cardinality $2^{2^{\aleph_0}}$. Another fact is that they have no interiors. One way to see this is to note that since the automorphisms used to bring the filters to each other are measure-preserving, all these closed sets must have the same measure as well. Being measurable subsets of the probability space $(\Lambda, \lambda_\Lambda)$ then implies they all have measure zero, so they could have no interior. A perhaps more illuminating way to prove this is to note that, were, say, $\Omega_{\mathcal{L}_a}$ for some $0 < a \leq 1$ to have interior, then the interior is clopen, so it corresponds to a measurable subset $A \subseteq [0,1]$ with positive measure. By the remark following Theorem 1, $\mathcal{L}_a$ is contained in the principal filter at $A$. Then $A$ is almost contained in $(a - \varepsilon, a) \cap [0,1]$ for all $\varepsilon > 0$. But as $\varepsilon \to 0$, $\lambda\big((a - \varepsilon, a) \cap [0,1]\big) \to 0$, so $\lambda(A) = 0$, a contradiction. Similarly, $\Omega_{\mathcal{R}_a}$, $0 \leq a < 1$ has no interior either. Hence,

**Theorem 15:** None of $\Omega_{\mathcal{L}_a}$, $0 < a \leq 1$, or $\Omega_{\mathcal{R}_a}$, $0 \leq a < 1$ has interior.

One more fact regarding this partition is that, while the division of ultrafilters into left and right limits is quite natural, this division is, perhaps surprisingly, not a measurable one, i.e., the union of all $\Omega_{\mathcal{L}_a}$ is not $\lambda_\Lambda$-measurable, and the same goes for the union of all $\Omega_{\mathcal{R}_a}$, as we shall now see.

**Theorem 16:** Neither $\bigcup_{0 < a \leq 1} \Omega_{\mathcal{L}_a}$ nor $\bigcup_{0 \leq a < 1} \Omega_{\mathcal{R}_a}$ is $\lambda_\Lambda$-measurable.

**Proof:** Assume either is measurable, then since they are complements of each other, both are measurable. As the reflection automorphism sending $a \in [0,1]$ to $1 - a \in [0,1]$ sends $\mathcal{L} = \bigcup_{0 < a \leq 1} \Omega_{\mathcal{L}_a}$ to $\mathcal{R} = \bigcup_{0 \leq a < 1} \Omega_{\mathcal{R}_a}$ and vice versa, and as the reflection automorphism is measure-preserving, we necessarily have both $\mathcal{L}$ and $\mathcal{R}$ are of measure $\frac{1}{2}$.

Fix a standard lifting map $\rho: L^\infty(\Omega, \mu) \to M^\infty(\Omega, \mu)$. By Corollary 6, $\phi_\rho^{-1}(\mathcal{L})$ is measurable with measure $\frac{1}{2}$. We consider the measure of $\phi_\rho^{-1}(\mathcal{L}) \cap (a, b)$ where $0 \leq a < b \leq 1$. $(a, b)$ equals $\phi_\rho^{-1}\big(\flat^{-1}((a,b))\big)$ a.e., so we may replace $(a, b)$ with $\phi_\rho^{-1}\big(\flat^{-1}((a,b))\big)$. Then $\phi_\rho^{-1}(\mathcal{L}) \cap \phi_\rho^{-1}\big(\flat^{-1}((a,b))\big) = \phi_\rho^{-1}\big(\mathcal{L} \cap \flat^{-1}((a,b))\big)$, so $\lambda\big(\phi_\rho^{-1}(\mathcal{L}) \cap (a,b)\big) = \lambda_\Lambda\big(\mathcal{L} \cap \flat^{-1}((a,b))\big)$. By Theorem 11, we may verify that,

$$\flat^{-1}((a,b)) = \left(\bigcup_{t \in (a,b)} \Omega_{\mathcal{D}_t}\right) \cup \Omega_{\mathcal{R}_a} \cup \Omega_{\mathcal{L}_b}$$

So,

$$\mathcal{L} \cap \flat^{-1}((a,b)) = \bigcup_{t \in (a,b]} \Omega_{\mathcal{L}_t}$$

Similarly, $\lambda\left(\phi_\rho^{-1}(\mathcal{R}) \cap (a,b)\right) = \lambda_\Lambda\left(\mathcal{R} \cap \mathfrak{b}^{-1}((a,b))\right)$ and $\mathcal{R} \cap \mathfrak{b}^{-1}((a,b)) = \bigcup_{t \in [a,b)} \Omega_{\mathcal{R}_t}$. We observe that the reflection automorphism brings $\bigcup_{t \in (a,b]} \Omega_{\mathcal{L}_t}$ to $\bigcup_{t \in [1-b, 1-a)} \Omega_{\mathcal{R}_t}$. The translation automorphism sending $t \in [0,1)$ to $t + x \pmod{1} \in [0,1)$ for $x = a+b-1$ then sends $\bigcup_{t \in [1-b, 1-a)} \Omega_{\mathcal{R}_t}$ to $\bigcup_{t \in [a,b)} \Omega_{\mathcal{R}_t}$. As both automorphisms are measure-preserving, we see that $\mathcal{L} \cap \mathfrak{b}^{-1}((a,b))$ and $\mathcal{R} \cap \mathfrak{b}^{-1}((a,b))$ has the same measure. Hence, both have measure $\frac{b-a}{2}$. Whence $\lambda\left(\phi_\rho^{-1}(\mathcal{L}) \cap (a,b)\right) = \frac{1}{2}\lambda((a,b))$.

Let $L$ be the collection of all measurable $F \subseteq [0,1]$ for which $\lambda(\phi_\rho^{-1}(\mathcal{L}) \cap F) = \frac{1}{2}\lambda(F)$. As $\phi_\rho^{-1}(\mathcal{L})$ has measure $\frac{1}{2}$, it is easy to see that $L$ is a $\lambda$-system, i.e., it contains the empty set, is closed under taking complements, and is closed under taking countable pairwise disjoint unions. Let $P$ contain all $(a,b)$ with $0 \leq a < b \leq 1$ as well as all null sets. We see that $P$ generates the $\sigma$-algebra of all measurable sets. We also note that $P$ is a $\pi$-system, i.e., it is closed under taking finite intersections. We claim that $P \subseteq L$. Indeed, we have seen that all $(a,b)$ are in $L$. That null sets are in $L$ are obvious. Hence, by the $\pi$-$\lambda$ theorem, all measurable sets are in $L$. In particular, $\phi_\rho^{-1}(\mathcal{L})$ is in $L$, so $\frac{1}{2} = \lambda\left(\phi_\rho^{-1}(\mathcal{L})\right) = \frac{1}{2}\lambda\left(\phi_\rho^{-1}(\mathcal{L})\right) = \frac{1}{4}$, a contradiction! ∎

We now return to the study of weak extensions and strong extensions. Our goal is to show $^s\text{ext}(F)$ is measurable for all measurable $F \subseteq [0,1]$. Then item 1 of Lemma 7 reduces it to showing $^w\text{ext}(F)$ is measurable for all Borel $F \subseteq [0,1]$.

**Lemma 17:** $^w\text{ext}(F)$ is measurable for all Borel $F \subseteq [0,1]$. Furthermore, $\lambda_\Lambda(^w\text{ext}(F)) = \lambda(F)$.

**Proof:** Let $L$ be the collection of all Borel $F \subseteq [0,1]$ for which $^w\text{ext}(F)$ is measurable and $\lambda_\Lambda(^w\text{ext}(F)) = \lambda(F)$. It is easy to see that $L$ is a $\lambda$-system. Let $P = \{(a,b): 0 \leq a < b \leq 1\} \cup \{\{0\}, \{1\}, \emptyset\} \subseteq \mathcal{P}([0,1])$. We see that $P$ generates the Borel $\sigma$-algebra and is a $\pi$-system. Hence, by the $\pi$-$\lambda$ theorem, it suffices to show $P \subseteq L$.

$\emptyset \in L$ is obvious. $\{0\} \in L$ follows from Theorem 15 that $\Omega_{\mathcal{D}_0} = \Omega_{\mathcal{R}_0}$ is closed with no interior. Similarly, $\{1\} \in L$. Finally, let $0 \leq a < b \leq 1$. Then, as in the proof of Theorem 16,

$$\mathfrak{b}^{-1}((a,b)) = \left(\bigcup_{t \in (a,b)} \Omega_{\mathcal{D}_t}\right) \cup \Omega_{\mathcal{R}_a} \cup \Omega_{\mathcal{L}_b}$$

But then $^w\text{ext}((a,b)) = \mathfrak{b}^{-1}((a,b)) \setminus (\Omega_{\mathcal{R}_a} \cup \Omega_{\mathcal{L}_b})$. Again, Theorem 15 implies that $\Omega_{\mathcal{R}_a}$ and $\Omega_{\mathcal{L}_b}$ are both closed with no interior, hence of zero measure. As $\mathfrak{b}^{-1}((a,b))$ is clopen with $\mu_\Lambda\left(\mathfrak{b}^{-1}((a,b))\right) = \mu((a,b))$, the result follows. ∎

**Remark:** We observe that the proof above actually implies $^w\text{ext}(F)$ is even Baire for Borel $F \subseteq [0,1]$. This follows from altering the definition of $L$ to require $^w\text{ext}(F)$ be Baire and observing that $\Omega_{\mathcal{R}_a}$ and $\Omega_{\mathcal{L}_b}$ are all $G_\delta$ sets. (For example, for $\Omega_{\mathcal{R}_a}$, $0 \leq a < 1$, the formula above for $\mathfrak{b}^{-1}((a,b))$ shows $\Omega_{\mathcal{R}_a} = \bigcap_{n=1}^\infty \mathfrak{b}^{-1}\left(\left(a, a + \frac{1}{n}\right) \cap [0,1]\right)$.) We also observe the above also implies $^w\text{ext}(F)$ is measurable whenever $F \subseteq [0,1]$ is measurable and $\lambda_\Lambda(^w\text{ext}(F)) = \lambda(F)$. This follows from observing that any measurable set differs from a Borel set by a null set and any null set is contained in a Borel null set. ∎

**Corollary 18:** $^s\text{ext}(F)$ is measurable for all measurable $F \subseteq [0, 1]$.

**Remark:** While $\lambda_\Lambda(^w\text{ext}(F)) = \lambda(F)$, the same is not true for strong extensions. We do have an inequality $\lambda_\Lambda(^s\text{ext}(F)) \geq \lambda(F)$. Indeed, fix a standard lifting map $\rho: L^\infty(\Omega, \mu) \to M^\infty(\Omega, \mu)$. Then we see that $\text{ints}(\mathfrak{b}^{-1}(F)) \supseteq \phi_\rho^{-1}(\mathfrak{b}^{-1}(F))$, so by Lemma 7, $^s\text{ext}(F) = {}^w\text{ext}\left(\text{ints}(\mathfrak{b}^{-1}(F))\right) \supseteq {}^w\text{ext}\left(\phi_\rho^{-1}(\mathfrak{b}^{-1}(F))\right)$. Hence, $\lambda_\Lambda(^s\text{ext}(F)) \geq {}^w\text{ext}\left(\phi_\rho^{-1}(\mathfrak{b}^{-1}(F))\right) = \lambda\left(\phi_\rho^{-1}(\mathfrak{b}^{-1}(F))\right) = \lambda(F)$. The inequality in the other direction, however, does not hold. Recall we have previously constructed $O \subseteq [0, 1]$ open dense with measure strictly smaller than 1. As we have seen, $K = [0, 1] \setminus O$ almost contains no open set, i.e., $^p\text{int}(K) = \emptyset$. But then Lemma 7 implies $^s\text{ext}(O) = \Lambda \setminus {}^w\text{ext}(^p\text{int}(K)) = \Lambda$. So while $\lambda(O) < 1$, $\lambda_\Lambda(^s\text{ext}(O)) = 1$.

This allows us to construct an example where condition A in Theorem 5 is satisfied but not condition B. Fix $O$ and $K$ as above. Then as $^s\text{ext}(O) = \Lambda$, $\mathfrak{b}^{-1}(O)$ intersects with $\Omega_{\mathcal{D}_t}$ for all $t$. Note that $\phi_\rho^{-1}(\mathfrak{b}^{-1}(K))$ has the same measure as $K$, which is positive. As such, there exists a non-measurable $A \subseteq \phi_\rho^{-1}(\mathfrak{b}^{-1}(K))$. Note that for all $t \in A \subseteq \phi_\rho^{-1}(\mathfrak{b}^{-1}(K))$, $\mathfrak{b}^{-1}(K)$ intersects with $\Omega_{\mathcal{D}_t}$. Hence, we may define a map $\phi: [0, 1] \to \Lambda$ by requiring $\phi(t) \in \mathfrak{b}^{-1}(K) \cap \Omega_{\mathcal{D}_t}$ for all $t \in A$ and $\phi(t) \in \mathfrak{b}^{-1}(O) \cap \Omega_{\mathcal{D}_t} = \Omega_{\mathcal{D}_t} \setminus \mathfrak{b}^{-1}(K)$ for all $t \notin A$. Then $\phi$ clearly satisfies condition A in Theorem 5, but as $\phi^{-1}(\mathfrak{b}^{-1}(K)) = A$, it does not satisfy condition B.

Moreover, we observe that, by altering the definition of $O$ somewhat we can have a decreasing sequence of open dense sets $O_n$ whose measures converge to 0. All of them satisfy $^s\text{ext}(O_n) = \Lambda$. Let $K = \bigcap_{n=1}^\infty \mathfrak{b}^{-1}(O_n)$. Then $K$ is compact and $G_\delta$. Furthermore, $K$ has measure zero, so it has no interior. By the same method as in the proof of Corollary 8, we see that $\text{ints}(K) = \bigcap_{n=1}^\infty \text{ints}(\mathfrak{b}^{-1}(O_n))$. But as $^s\text{ext}(O_n) = \Lambda$, we see that $\text{ints}(\mathfrak{b}^{-1}(O_n)) = [0, 1]$ for all $n$, so $\text{ints}(K) = [0, 1]$, i.e., $K$ intersects with $\Omega_{\mathcal{D}_t}$ for all $t$. Thus, we may construct a map $\phi: [0, 1] \to \Lambda$ by requiring $\phi(t) \in K \cap \Omega_{\mathcal{D}_t}$ for all $t$. Such a map satisfies condition A in Theorem 5, but as $\phi^{-1}(K) = [0, 1]$, it does not satisfy condition C. ∎

While condition A itself does not imply condition C in Theorem 5, it is unclear whether condition A plus B implies condition C. In fact, we have the following,

**Theorem 19:** Let $\Omega$ be a compact Hausdorff space s.t. $C(\Omega)$ is norm-separable, $\mu$ be a faithful Radon probability measure on $\Omega$. Then the following are equivalent,

1. Conditions A and B in Theorem 5 together imply condition C, as applied to $(\Omega, \mu)$;
2. Any *-homomorphism $\beta: L^\infty(\Omega, \mu) \to L^\infty(\Omega, \mu)$ with $\beta|_{C(\Omega)} = Id_{C(\Omega)}$ must be the identity map on the entire $L^\infty(\Omega, \mu)$.

**Proof:** (1 ⇒ 2) Fix a standard lifting map $\rho: L^\infty(\Omega, \mu) \to M^\infty(\Omega, \mu)$. Then $\rho \circ \beta: L^\infty(\Omega, \mu) \to M^\infty(\Omega, \mu)$ is a *-homomorphism that preserves $C(\Omega)$ pointwise. By the proof of 1 ⇒ 2 in Theorem 5, we see that $\phi_{\rho \circ \beta}$ satisfies conditions A and B. Then condition C must be satisfied as well, whence $\rho \circ \beta$ is a lifting map, which implies $\beta$ is the identity map.

(2 ⇒ 1) Suppose for some $\rho: L^\infty(\Omega, \mu) \to l^\infty(\Omega)$, $\phi_\rho$ satisfies conditions A and B. Let $\pi: M^\infty(\Omega, \mu) \to L^\infty(\Omega, \mu)$ be the canonical quotient map. As we have seen in the proof of 2 ⇒ 1 in Theorem 5, $\phi_\rho$ satisfying conditions A and B implies $\rho$ has range in $M^\infty(\Omega, \mu)$ and $\pi \circ \rho|_{C(\Omega)} = Id_{C(\Omega)}$. Hence $\pi \circ \rho = Id_{L^\infty(\Omega, \mu)}$, i.e., $\rho$ is a lifting map and thus $\phi_\rho$ satisfies condition C. ∎

As we have seen in the remark following Corollary 18, one of the main obstructions to proving the equivalent statements in the theorem above hold is the existence of open dense sets in $[0, 1]$ with arbitrarily small measures. As such, it is likely that any *-homomorphism $\beta\colon L^\infty(\Omega, \mu) \to L^\infty(\Omega, \mu)$ that serves as a counterexample to statement 2 in the theorem above (if one exists) will need to be constructed by requiring it to not act as the identity on the indicator functions of some open dense sets.

# Section IV: Tensor Products and Their Representations on the Spectrum

We shall revert to assuming $(\Omega, \Sigma, \mu)$ is a semifinite, decomposable measure space. We now consider tensor products with $L^\infty(\Omega, \mu)$. Let $M$ be a von Neumann algebra, $a \in L^\infty(\Omega, \mu) \overline{\otimes} M$, then,

**Proposition 1:** The function $\Lambda \ni \omega \mapsto \omega \otimes Id(a) \in M$ is uniformly bounded and continuous when $M$ is equipped with the weak* topology.

**Proof:** Uniform boundedness is obvious. To prove continuity, it suffices to show $\psi(\omega \otimes Id(a))$ is a continuous function in $\omega$ for all $\psi \in M_*$. We have $\psi(\omega \otimes Id(a)) = \omega(Id \otimes \psi(a))$. $Id \otimes \psi(a)$, being an element of $L^\infty(\Omega, \mu)$, is represented as a continuous function on $\Lambda$, whence the result follows. ∎

**Proposition 2:** $\overline{\text{span}\{\omega \otimes Id(a) | \omega \in \Lambda\}}^{weak^*} = E_M^{sp}(a)$.

**Proof:** We first show that $\omega \otimes Id(a) \in E_M^{sp}(a)$ for all $\omega \in \Lambda$. Assume otherwise. By Hahn-Banach, there exists $\psi \in M_*$ s.t. $\psi(\varphi \otimes Id(a)) = 0$ for all $\varphi \in L^\infty(\Omega, \mu)_*$ and $\psi(\omega \otimes Id(a)) = 1$. However, we have $\varphi(Id \otimes \psi(a)) = \varphi \otimes \psi(a) = \psi(\varphi \otimes Id(a)) = 0$ for all $\varphi \in L^\infty(\Omega, \mu)_*$, so $Id \otimes \psi(a) = 0$. Hence, $\psi(\omega \otimes Id(a)) = \omega(Id \otimes \psi(a)) = 0$, a contradiction!

Conversely, to show $E_M^{sp}(a) \subseteq \overline{\text{span}\{\omega \otimes Id(a) | \omega \in \Lambda\}}^{weak^*}$, it suffices to show $\varphi \otimes Id(a) \in \overline{\text{span}\{\omega \otimes Id(a) | \omega \in \Lambda\}}^{weak^*}$ for all $\varphi \in L^\infty(\Omega, \mu)_*$. Since normal states span $L^\infty(\Omega, \mu)_*$, it suffices to consider the case where $\varphi$ is a normal state. Now, assume $\varphi \otimes Id(a) \notin \overline{\text{span}\{\omega \otimes Id(a) | \omega \in \Lambda\}}^{weak^*}$. Again by Hahn-Banach, there exists $\psi \in M_*$ s.t. $\psi(\omega \otimes Id(a)) = 0$ for all $\omega \in \Lambda$ and $\psi(\varphi \otimes Id(a)) = 1$. We have $0 = \psi(\omega \otimes Id(a)) = \omega(\psi \otimes Id(a))$ and $1 = \psi(\varphi \otimes Id(a)) = \varphi(\psi \otimes Id(a))$. Recall that the $\sigma(L^\infty(\Omega, \mu)^*, L^\infty(\Omega, \mu))$-closed convex hall of $\Lambda$ in $L^\infty(\Omega, \mu)^*$ is the collection of all states (see [Mur90, Theorem 5.1.6(2)] and [Mur90, Corollary 5.1.10]). Since $\psi \otimes Id(a)$ is a fixed element of $L^\infty(\Omega, \mu)$ and $\omega(\psi \otimes Id(a)) = 0$ for all $\omega \in \Lambda$, we must have $\varphi(\psi \otimes Id(a)) = 0$, a contradiction! ∎

**Definition:** For any topological space $X$, let $C_w(X; M)$ be the space of all uniformly bounded, continuous functions from $X$ to $M$ where $M$ is equipped with the weak* topology. It is a Banach space under the uniform norm. Note that it is a *-closed subspace of $l^\infty(X) \overline{\otimes} M$ (where in $l^\infty(X)$, $X$ is understood to be a discrete set) containing the identity, so it has a natural operator system structure. For $\varphi \in M_*$ and $f \in C_w(X; M)$, we let $Id \otimes \varphi(f) = \varphi \circ f$. Note that since $f$ is a uniformly bounded continuous function into $M$ equipped with the weak* topology, $Id \otimes \varphi(f) \in C_b(X)$.

**Lemma 3:** Let $M$ acts on $H$, then elements of $C_w(X; M)$ are in one-to-one correspondence with sesquilinear forms $\eta: H \times H \to C_b(X)$ satisfying:

1. $\eta$ is bounded, i.e., there exists $C \geq 0$ s.t. $\|\eta(h, k)\| \leq C\|h\|\|k\|$;
2. $\eta$ commutes with $M'$, i.e., $\eta(h, xk) = \eta(x^*h, k)$ for all $x \in M'$.

The correspondence associates $a \in C_w(X; M)$ to the sesquilinear form $\eta$ defined by $\eta(h, k) = Id \otimes |k\rangle\langle h|(a)$.

**Proof:** We first show that, given $a \in C_w(X; M)$, $\eta(h, k) = Id \otimes |k\rangle\langle h|(a)$ indeed satisfies the given conditions. The first condition follows from $a$ being uniformly bounded. For the second condition, let $\omega \in X$, then by definition,

$$\eta(h, xk)(\omega) = |xk\rangle\langle h|(a(\omega)) = \langle h, a(\omega)xk\rangle = \langle h, xa(\omega)k\rangle = \langle x^*h, a(\omega)k\rangle = \eta(x^*h, k)(\omega)$$

For all $x \in M'$.

Now, given a bounded sesquilinear form $\eta: H \times H \to C_b(X)$ that commutes with $M'$, we shall construct an element of $C_w(X; M)$. Let the norm of the sesquilinear form be $C$. For any $\omega \in X$, $\eta_\omega(h, k) = \eta(h, k)(\omega)$ gives a sesquilinear form $H \times H \to \mathbb{C}$ with norm bounded by $C$. Thus, it gives a bounded operator on $H$. Since $\eta$ commutes with $M'$, $\eta_\omega$ does so as well, so the corresponding operator lies in $M$. We thus have a uniformly bounded function $f$ from $X$ to $M$. Since $|k\rangle\langle h| \circ f = \eta(h, k) \in C_b(X)$ by construction, we see that $f$ is WOT-continuous. On bounded sets, WOT and the weak* topology coincide, so $f \in C_w(X; M)$. It is not hard to see that the two maps defined above are inverses to each other, whence the result follows. ∎

We now look at $C_w(\Lambda; M)$ where $\Lambda$ is the spectrum of $L^\infty(\Omega, \mu)$:

**Proposition 4:** Let $\iota: L^\infty(\Omega, \mu) \overline{\otimes} M \to C_w(\Lambda; M)$ be defined by,

$$\iota(a)(\omega) = \omega \otimes Id(a)$$

Then $\iota$ is an injective ucp map.

**Proof:** Proposition 1 shows the map indeed has range within $C_w(\Lambda; M)$. Proposition 2 implies injectivity. Going through the definitions, it is not hard to see the map is unital. To show positivity, one may adapt the first paragraph of the proof of Lemma II.10 to see that $Id \otimes \psi(a) \geq 0$ whenever $\psi \in (M_*)_+$ and $a \geq 0$. But then $\psi(\iota(a)(\omega)) = \psi(\omega \otimes Id(a)) = \omega(Id \otimes \psi(a)) \geq 0$. Since $\psi \in (M_*)_+$ is arbitrary, this gives the result.

Now, to show complete positivity, one observes that $\mathbb{M}_n \otimes (L^\infty(\Omega, \mu) \overline{\otimes} M) \cong L^\infty(\Omega, \mu) \overline{\otimes} (\mathbb{M}_n \otimes M)$ and $\mathbb{M}_n \otimes C_w(\Lambda; M) \cong C_w(\Lambda; \mathbb{M}_n \otimes M)$ naturally. So the positivity of $Id_{\mathbb{M}_n} \otimes \iota: \mathbb{M}_n \otimes (L^\infty(\Omega, \mu) \overline{\otimes} M) \to \mathbb{M}_n \otimes C_w(\Lambda; M)$ follows from the positivity of $\iota: L^\infty(\Omega, \mu) \overline{\otimes} (\mathbb{M}_n \otimes M) \to C_w(\Lambda; \mathbb{M}_n \otimes M)$. (The precise technical details proving these two maps are indeed the same are left for readers to verify.) ∎

We now focus on the case where $L^\infty(\Omega, \mu)$ is separable. Then as $L^\infty(\Omega, \mu)$ is countably decomposable, it admits a normal faithful tracial state, which induces a probability measure on $\Omega$ equivalent to $\mu$. Hence, we may assume $\mu$ is a probability measure. (These assumptions shall be in effect for the rest of this section.) It therefore induces a normal Radon probability measure on $\Lambda$, which we shall also denote by $\mu$. Since $\mu$ is faithful, we see that its support is $\Lambda$.

**Lemma 5:** Given $f \in C_w(\Lambda; M)$, its range lies within a countably generated von Neumann subalgebra of $M$.

**Proof:** We prove this by adapting the proof of Proposition 2 and Corollary I.4.

We need to first define $\varphi \otimes Id(f) \in M = (M_*)^*$ for $\varphi \in L^\infty(\Omega, \mu)_*$. This is given by $(\varphi \otimes Id(f))(\psi) = \varphi(\psi \circ f)$, where we note $L^\infty(\Omega, \mu) = C(\Lambda)$. An alternative way of understanding this would be to consider $f$ as a bounded sesquilinear form $\eta: H \times H \to C(\Lambda)$ commuting with $M'$, in which case $\varphi \otimes Id(f)$ is simply given by $\varphi \circ f$ as a sesquilinear form.

We now claim $\{f(\omega): \omega \in \Lambda\} \subseteq \overline{\{\varphi \otimes Id(f) | \varphi \in L^\infty(\Omega, \mu)_*\}}^{weak^*}$. Indeed, assume to the contrary. Then by Hahn-Banach there exists $\omega \in \Lambda$ and $\psi \in M_*$ s.t. $\psi(\varphi \otimes Id(f)) = 0$ for all $\varphi \in L^\infty(\Omega, \mu)_*$ and

$\psi(f(\omega)) = 1$. But $0 = \psi(\varphi \otimes Id(f)) = \varphi(\psi \circ f)$. Since this holds for all $\varphi \in L^\infty(\Omega, \mu)_*$, $\psi \circ f = 0$, contradicting $\psi(f(\omega)) = 1$.

Observing that $\|\varphi \otimes Id(f)\| \leq \|\varphi\| \|f\|$ and following the proof of Corollary I.4 yields the result. ∎

**Corollary 6:** Elements of $C_w(\Lambda; M)$, regarded as functions from $(\Lambda, \mu)$ to $M$, are strongly measurable.

**Proof:** Let $f \in C_w(\Lambda; M)$. $\Lambda \ni \omega \mapsto \langle h, f(\omega)k \rangle = (Id \otimes |k\rangle\langle h|(f))(\omega)$ is continuous by Lemma 3. Thus, it suffices to show that, for all $h \in H$, $\{f(\omega)h : \omega \in \Lambda\}$ lies within a separable subspace of $H$. By Lemma 5, all $f(\omega)$ lie within a countably generated von Neumann subalgebra $M_0 \subseteq M$. In particular, $M_0$ admits a countable subset $S$ which is SOT-dense. Hence, $\{f(\omega)h : \omega \in \Lambda\} \subseteq M_0 h \subseteq \overline{Sh}$, which is separable. ∎

By the proof of (i) $\Rightarrow$ (iii) in [Tak79, Theorem III.1.18], we see that the canonical map $L^\infty(\Omega, \mu) = C(\Lambda) \to L^\infty(\Lambda, \mu)$ is an isomorphism. Therefore, $L^\infty(\Lambda, \mu) \overline{\otimes} M \cong L^\infty(\Omega, \mu) \overline{\otimes} M$. Now, by the previous corollary and the remark preceding Lemma I.5, each element of $C_w(\Lambda; M)$ yields an element of $L^\infty(\Lambda, \mu) \overline{\otimes} M \cong L^\infty(\Omega, \mu) \overline{\otimes} M$ via pointwise multiplication on $L^2(\Lambda, \mu; H)$. Let $J: C_w(\Lambda; M) \to L^\infty(\Omega, \mu) \overline{\otimes} M$ denote this map. Representing $L^\infty(\Omega, \mu) \overline{\otimes} M$ on $L^2(\Lambda, \mu; H)$, it is not hard to see that the map is norm-decreasing, unital, and positive. Complete positivity can be proved following the same method as in the second paragraph of the proof of Proposition 4.

**Theorem 7:** The map $\iota$ defined in Proposition 4 and the map $J$ just defined are inverses to each other. Hence, $L^\infty(\Omega, \mu) \overline{\otimes} M$ and $C_w(\Lambda; M)$ are isomorphic as operator systems, with $\iota$ and $J$ both isomorphisms.

**Proof:** Since $\iota: L^\infty(\Omega, \mu) \overline{\otimes} M \to C_w(\Lambda; M)$ is injective, it suffices to show $\iota \circ J = Id$. Let $L^\infty(\Omega, \mu) \overline{\otimes} M$ acts on $L^2(\Lambda, \mu; H)$. Fix $f \in C_w(\Lambda; M)$, $\omega_0 \in \Lambda$. We need to show $\omega_0 \otimes Id(J(f)) = f(\omega_0)$.

Since $M$ acts on $H$, it suffices to show $\langle h, \omega_0 \otimes Id(J(f))k \rangle = \langle h, f(\omega_0)k \rangle$ for all $h, k \in H$. Observe that $\langle h, \omega_0 \otimes Id(J(f))k \rangle = |k\rangle\langle h| \big(\omega_0 \otimes Id(J(f))\big) = \omega_0 \big(Id \otimes |k\rangle\langle h|(J(f))\big)$. $f'_{h,k} = Id \otimes |k\rangle\langle h|(J(f))$ is an element of $L^\infty(\Omega, \mu) = C(\Lambda)$. Thus, for any $\varepsilon > 0$, there is a clopen neighborhood $F_1$ of $\omega_0$ s.t. $|\omega(f'_{h,k}) - \omega_0(f'_{h,k})| < \varepsilon$ for all $\omega \in F_1$. Moreover, $f_{h,k}(\omega) = \langle h, f(\omega)k \rangle$ is also a continuous function as $f \in C_w(\Lambda; M)$, so there is a clopen neighborhood $F_2$ of $\omega_0$ s.t. $|f_{h,k}(\omega) - f_{h,k}(\omega_0)| < \varepsilon$ for all $\omega \in F_2$. Let $F = F_1 \cap F_2$.

Since $\mu$ has full support, $\mu(F) > 0$. Thus, $\frac{1_F}{\mu(F)} \in L^1(\Lambda, \mu)$. It is easy to see it induces a normal state on $L^\infty(\Lambda, \mu)$, which we shall denote by $\chi_F$. We have,

$$\left|\chi_F(f'_{h,k}) - \omega_0(f'_{h,k})\right| = \left|\left(\frac{1}{\mu(F)} \int_F f'_{h,k}\, d\mu\right) - \omega_0(f'_{h,k})\right| \leq \varepsilon$$

Similarly, $\left|\chi_F(f_{h,k}) - \omega_0(f_{h,k})\right| \leq \varepsilon$. Now,

$\chi_F(f'_{h,k}) = \chi_F \otimes |k\rangle\langle h|(J(f))$

$= \frac{1}{\mu(F)} |1_F\rangle\langle 1| \otimes |k\rangle\langle h|(J(f))$

$$= \frac{1}{\mu(F)} \langle 1 \otimes h, J(f)(1_F \otimes k) \rangle$$

As a measurable function, $(J(f)(1_F \otimes k))(\omega) = 1_F(\omega) f(\omega) k$, so,

$$\chi_F(f'_{h,k}) = \frac{1}{\mu(F)} \int_F \langle h, f(\omega)k \rangle \, d\mu = \frac{1}{\mu(F)} \int_F f_{h,k} \, d\mu = \chi_F(f_{h,k})$$

Thus, $|\omega_0(f'_{h,k}) - \omega_0(f_{h,k})| \leq 2\varepsilon$. As $\varepsilon > 0$ is arbitrary, $\omega_0(f'_{h,k}) = \omega_0(f_{h,k})$, i.e., $\langle h, \omega_0 \otimes Id(J(f))k \rangle = \langle h, f(\omega_0)k \rangle$. ∎

**Remark:** For an explicit example of the spectrum representation of an element of $L^\infty(\Omega, \mu) \overline{\otimes} M$, see the example following Theorem III.11. ∎

This isomorphism means that there is a multiplication on $C_w(\Lambda; M)$ consistent with its operator system structure, given by $f * g = \iota(J(f)J(g))$. In the general case, the best we can say about this multiplication is follows:

**Proposition 8:** Let $f, g \in C_w(\Lambda; M)$, $\varphi \in M_*$. Then there exists an open dense subset $O \subseteq \Lambda$ s.t. $\varphi((f * g)(\omega)) = \varphi(f(\omega)g(\omega))$ for all $\omega \in O$.

**Proof:** We represent elements of $C_w(\Lambda; M)$ as operators on $L^2(\Lambda, \mu; H)$. Then $f * g$, being the composition of the operators $f$ and $g$, is given by pointwise multiplication by $f(\omega)g(\omega)$.

We first consider the case where $\varphi = |k\rangle\langle h|$ for some $h, k \in H$. We first note that, since $f^* \in C_w(\Lambda; M)$, $\langle h, f(\omega)g(\omega)k \rangle = \langle f(\omega)^*h, g(\omega)k \rangle$ is measurable. Hence $\eta(\omega) = \langle h, (f(\omega)g(\omega) - (f * g)(\omega))k \rangle = \langle h, f(\omega)g(\omega)k \rangle - \langle h, (f * g)(\omega)k \rangle$ is a function in $L^\infty(\Lambda, \mu)$. Let $d \in L^1(\Lambda, \mu) = L^\infty(\Lambda, \mu)_*$. Let the polar decomposition of $d$ be $d = u|d|$, so,

$$d(\eta) = \int_\Lambda d(\omega)\eta(\omega) \, d\mu = \int_\Lambda \langle |d|^{1/2}(\omega)h, (f(\omega)g(\omega) - (f * g)(\omega))u(\omega)|d|^{1/2}(\omega)k \rangle \, d\mu$$

Observe that, by definition of the inner product on $L^2(\Lambda, \mu; H)$, this is the inner product of $|d|^{1/2} \otimes h$ with a vector that, in the form of a measurable function, is $(f(\omega)g(\omega) - (f * g)(\omega))u(\omega)|d|^{1/2}(\omega)k = f(\omega)g(\omega)u(\omega)|d|^{1/2}(\omega)k - (f * g)(\omega)u(\omega)|d|^{1/2}(\omega)k$. The vector $(f * g)(\omega)u(\omega)|d|^{1/2}(\omega)k$ is the result of the operator $f * g$ acting on the vector $u|d|^{1/2} \otimes k$. As we have observed, $f * g$, as an operator, is pointwise multiplication by $f(\omega)g(\omega)$. Thus, as a vector, $(f * g)(\omega)u(\omega)|d|^{1/2}(\omega)k = f(\omega)g(\omega)u(\omega)|d|^{1/2}(\omega)k$ (i.e., the two expressions coincide $\mu$-a.e.), so $(f(\omega)g(\omega) - (f * g)(\omega))u(\omega)|d|^{1/2}(\omega)k$, as a vector, is 0. Therefore, $d(\eta) = 0$. As $d \in L^1(\Lambda, \mu) = L^\infty(\Lambda, \mu)_*$ is arbitrary, $\eta = 0$ in $L^\infty(\Lambda, \mu)$, i.e., $\eta(\omega) = 0$ $\mu$-a.e., i.e., $\varphi((f * g)(\omega)) = \varphi(f(\omega)g(\omega))$ $\mu$-a.e. when $\varphi = |k\rangle\langle h|$. Since the linear span of $\varphi$ of such form is dense in $M_*$, we see that $\varphi((f * g)(\omega)) = \varphi(f(\omega)g(\omega))$ $\mu$-a.e. for all $\varphi \in M_*$. Finally, by [Tak79, Proposition III.1.15], the set on which $\varphi((f * g)(\omega)) \neq \varphi(f(\omega)g(\omega))$, being $\mu$-null, must be nowhere dense, so there is an open dense set on which $\varphi((f * g)(\omega)) = \varphi(f(\omega)g(\omega))$. ∎

**Remark:** One conclusion from the above proposition is that there is only one way to "alter" $f(\omega)g(\omega)$ to obtain $(f * g)(\omega)$. More precisely, for any fixed $\varphi \in M_*$, there is only one continuous function (namely

$\varphi((f * g)(\omega)))$ which coincides with $\varphi(f(\omega)g(\omega))$ on an open dense set. Indeed, let $h \in C(\Lambda)$ be such a function, then by [Tak79, Proposition III.1.15], it coincides with $\varphi(f(\omega)g(\omega))$ $\mu$-a.e. The same is true for $\varphi((f * g)(\omega))$, which means $\varphi((f * g)(\omega)) = h(\omega)$ $\mu$-a.e. Applying [Tak79, Proposition III.1.15] again, we see that $\varphi((f * g)(\omega)) = h(\omega)$ on an open dense set. Since both functions are continuous, they must be the same. One particular implication from this observation is that, if $f(\omega)g(\omega)$ is already weak* continuous, then $(f * g)(\omega) = f(\omega)g(\omega)$ for all $\omega \in \Lambda$. ∎

We may improve on the above proposition significantly if $M$ is assumed to be countably decomposable:

**Theorem 9:** Assume $M$ is countably decomposable. Let $f, g \in C_w(\Lambda; M)$, then there exists an open dense subset $O \subseteq \Lambda$ s.t. $(f * g)(\omega) = f(\omega)g(\omega)$ for all $\omega \in O$.

**Proof:** By Lemma 5, there is a countably generated subalgebra $M_0 \subseteq M$ s.t. $f(\omega), g(\omega), (f * g)(\omega) \in M_0$ for all $\omega \in \Lambda$. By Theorem I.8, $M_0$ is separable. We may regard $f, g$, and $f * g$ as elements of $C_w(\Lambda; M_0)$. Pick a countable dense subset $S \subseteq (M_0)_*$. Then applying Proposition 8 and again using the fact that in $\Lambda$, a set is $\mu$-null iff it is nowhere dense, we see that there exists an open dense subset $O \subseteq \Lambda$ s.t. $\varphi((f * g)(\omega)) = \varphi(f(\omega)g(\omega))$ for all $\varphi \in S$ and $\omega \in O$. Since $S \subseteq (M_0)_*$ is dense, the result follows. ∎

**Remark:** For an example where Theorem 9 fails without the assumption of countable decomposability, see the example preceding Theorem V.19. Even assuming countable decomposability, Theorem 9 is the best that can be done in most of the interesting cases. Indeed, consider the following example: Let $\Omega$ be the discrete set $\mathbb{N}$ and $\mu$ be the counting measure. Then $\Lambda$ is the Stone-Čech compactification of $\mathbb{N}$ and contains $\mathbb{N}$ as an open dense subset. In that case any $a \in l^\infty(\mathbb{N}) \overline{\otimes} M$ can be regarded as a function from $\mathbb{N}$ to $M$. And it is easy to verify that, when restricted to $\mathbb{N}$, $a$ as a function in $C_w(\Lambda; M)$ and $a$ as a function from $\mathbb{N}$ to $M$ coincide. Now, consider any separable, diffuse, and tracial von Neumann algebra $M$. By [AP18, Exercises 3.2 & 3.3], there exists a unitary $u \in M$ s.t. $u^n \to 0$ weak* as $n \to \infty$. Let $a \in l^\infty(\mathbb{N}) \overline{\otimes} M$ be defined, as a function from $\mathbb{N}$ to $M$, by $a(n) = u^n$. As a function in $C_w(\Lambda; M)$, it is not hard to see $a$ takes value 0 on all free ultrafilters. (See the remark following Corollary III.6.) The same holds for $a^*$. But $aa^* = 1$, so $(a * a^*)(\omega) = a(\omega)a^*(\omega)$ only when $\omega \in \mathbb{N}$. ∎

The situation is also easier to deal with if the elements under consideration lie within the min tensor product $L^\infty(\Omega, \mu) \otimes_{min} M \subseteq L^\infty(\Omega, \mu) \overline{\otimes} M$. Since $L^\infty(\Omega, \mu) = C(\Lambda)$, by [Mur90, Theorem 6.4.17], the min tensor product $L^\infty(\Omega, \mu) \otimes_{min} M = C(\Lambda) \otimes_{min} M$ is canonically isomorphic to $C(\Lambda; M)$, which consists of norm-continuous functions from $\Lambda$ to $M$, by sending $f \otimes m$ to the function $\Lambda \ni \omega \mapsto f(\omega)m$, for $f \in C(\Lambda), m \in M$. Our isomorphism $\iota: L^\infty(\Omega, \mu) \overline{\otimes} M \to C_w(\Lambda; M)$, when restricted to the min tensor product, does the exact same thing. Indeed, in light of [Mur90, Theorem 6.4.17], it is sufficient to verify that $\omega \otimes Id(f \otimes m) = f(\omega)m$, which follows easily from the definition. Thus, we have,

**Proposition 10:** $\iota: L^\infty(\Omega, \mu) \overline{\otimes} M \to C_w(\Lambda; M)$, when restricted to $L^\infty(\Omega, \mu) \otimes_{min} M$, maps onto $C(\Lambda; M) \subseteq C_w(\Lambda; M)$.

This also means that the multiplication on $C_w(\Lambda; M)$, when restricted to $C(\Lambda; M)$, actually does work pointwise. In fact, more can be said:

**Theorem 11:** Let $f, g \in C_w(\Lambda; M)$, if at least one of them belongs to $C(\Lambda; M)$, then $(f * g)(\omega) = f(\omega)g(\omega)$ for all $\omega \in \Lambda$.

**Proof:** We shall prove the case where $f \in C(\Lambda; M)$. The remainder can be proved similarly. By the remark after Proposition 8, it suffices to show $f(\omega)g(\omega)$ is weak* continuous. We first prove this when $f = h \otimes m$ for $h \in C(\Lambda), m \in M$. Let $\varphi \in M_*$, then $\varphi(f(\omega)g(\omega)) = h(\omega)\varphi(mg(\omega))$. Since multiplication by a fixed element is weak* continuous, this is a continuous function. Taking linear combinations and norm-limits yields the result. ∎

Recall that for a Radon probability measure, all bounded measurable functions are pointwise a.e. limits of sequences of continuous functions (see, for example, the proof of [Fol99, Theorem 7.10]). If $M$ is countably decomposable, then the same can be said for functions in $C_w(\Lambda; M)$:

**Theorem 12:** Assume $M$ is countably decomposable. Let $f \in C_w(\Lambda; M)$, then there exists a sequence $(f_n)_{n=1}^\infty \subseteq C(\Lambda; M)$ and an open dense subset $O \subseteq \Lambda$ s.t. $f_n(\omega) \to f(\omega)$ weak* for all $\omega \in O$. The sequence may be chosen so that $\|f_n\| \leq \|f\|$ for all $n$ and $f_n \to f$ weak*.

**Proof:** We may assume $\|f\| = 1$. By Lemma 5 and Theorem I.8, $f$ has its range within a separable subalgebra of $M$. Thus, we may assume $M$ is separable. As such, $L^\infty(\Omega, \mu) \overline{\otimes} M$ is separable, so its unit ball, under the weak* topology, is second-countable. Then, by Kaplansky density theorem, we may choose a sequence $(g_i)_{i=1}^\infty$ in the unit ball of $L^\infty(\Omega, \mu) \otimes_{min} M$ s.t. $g_i \to f$ weak*. Fix a countable dense subset $\{\varphi_j\}_{j=1}^\infty \subseteq M_*$. Then, for each fixed $\varphi_j$, $Id \otimes \varphi_j(g_i) \to Id \otimes \varphi_j(f)$ weak*. In particular, as functions in $C(\Lambda)$, $\int h(\omega) \left(Id \otimes \varphi_j(g_i)\right)(\omega) d\mu \to \int h(\omega) \left(Id \otimes \varphi_j(f)\right)(\omega) d\mu$ for each $h \in L^\infty(\Lambda, \mu) \subseteq L^1(\Lambda, \mu)$, i.e., $Id \otimes \varphi_j(g_i) \to Id \otimes \varphi_j(f)$ Banach-weakly as $L^1$-functions.

We construct $(f_{mn})_{m,n=1}^\infty$ inductively as follows. We begin with $m = 1$. By Hahn-Banach, the closure of the convex hall $co\{Id \otimes \varphi_1(g_i): i \geq N\} = Id \otimes \varphi_1(co\{g_i: i \geq N\})$ under $L^1$-norm contains $Id \otimes \varphi_1(f)$ for all $N$. Let $(g_{1N})_{N=1}^\infty$ be defined so that $g_{1N} \in co\{g_i: i \geq N\}$ and $\|Id \otimes \varphi_1(g_{1N}) - Id \otimes \varphi_1(f)\|_{L^1(\Lambda,\mu)} < \frac{1}{N}$. Then $Id \otimes \varphi_1(g_{1N}) \to Id \otimes \varphi_1(f)$ in $L^1$-norm, so a subsequence must converge a.e. This subsequence comes from a subsequence $(g_{1N_n})_{n=1}^\infty$ of $(g_{1N})_{N=1}^\infty$, and we let $f_{1n} = g_{1N_n}$. Observe that we have $f_{1n} \in co\{g_i: i \geq n\}$, so $f_{1n} \to f$ weak*.

Suppose we have constructed all $f_{mn}$ for $m \leq M$ s.t. $Id \otimes \varphi_m(f_{mn}) \to Id \otimes \varphi_m(f)$ a.e. as $n \to \infty$ for all $m \leq M$ and $f_{mN} \in co\{f_{(m-1)n}: n \geq N\}$ for all $1 < m \leq M$ and $N \in \mathbb{N}$. We now construct $f_{(M+1)n}$ as follows: By the convex hall condition, we have $f_{Mn} \to f$ weak* as $n \to \infty$. Similar to the arguments above, we see that the $L^1$-closure of the convex hall $co\{Id \otimes \varphi_{M+1}(f_{Mn}): n \geq N\} = Id \otimes \varphi_{M+1}(co\{f_{Mn}: n \geq N\})$ contains $Id \otimes \varphi_{M+1}(f)$ for all $N$. Hence, we may choose $g_{(M+1)N} \in co\{f_{Mn}: n \geq N\}$ and $\|Id \otimes \varphi_{M+1}(g_{(M+1)N}) - Id \otimes \varphi_{M+1}(f)\|_{L^1(\Lambda,\mu)} < \frac{1}{N}$. Then $Id \otimes \varphi_{M+1}(g_{(M+1)N}) \to Id \otimes \varphi_{M+1}(f)$ in $L^1$-norm, so a subsequence must converge a.e. This subsequence comes from a subsequence $(g_{(M+1)N_n})_{n=1}^\infty$ of $(g_{(M+1)N})_{N=1}^\infty$, and we let $f_{(M+1)n} = g_{(M+1)N_n}$. It is easy to see that this sequence satisfies the conditions required for induction.

Let $f_n = f_{nn}$. We claim this sequence satisfies the conditions of the theorem. Indeed, by definition $f_n \in co\{f_{mi}: i \geq n\}$ for all $m \leq n$. Since $Id \otimes \varphi_m(f_{mi}) \to Id \otimes \varphi_m(f)$ a.e. as $i \to \infty$, we have $Id \otimes \varphi_m(f_n) \to Id \otimes \varphi_m(f)$ a.e. as $n \to \infty$. This holds for all $m$. Since there are only countably many $m$, there is a shared set of measure 1 on which $Id \otimes \varphi_m(f_n) \to Id \otimes \varphi_m(f)$ holds pointwise for all $m$. By [Tak79, Proposition III.1.15], we may take this set of measure 1 as an open dense subset $O \subseteq \Lambda$. Since $f_{1n} \in co\{g_i: i \geq n\}$, we have $f_n \in co\{g_i: i \geq n\}$. Because $g_i \to f$ weak*, $f_n \to f$ weak* as well.

Moreover, as $g_i \in (L^\infty(\Omega, \mu) \otimes_{min} M)_1$ for all $i$, $f_n \in (L^\infty(\Omega, \mu) \otimes_{min} M)_1 = (C(\Lambda; M))_1$ for all $n$, so by the density of $\{\varphi_j\}_{j=1}^\infty$ in $M_*$, we have $Id \otimes \varphi(f_n) \to Id \otimes \varphi(f)$ pointwise on $O$ for all $\varphi \in M_*$. Hence, $\varphi(\omega \otimes Id(f_n)) \to \varphi(\omega \otimes Id(f))$ for all $\varphi \in M_*$ and $\omega \in O$, so $f_n(\omega) \to f(\omega)$ weak* for all $\omega \in O$. ∎

**Remark:** We note that the method used in proving this theorem actually showed the following: Assume $M$ is countably decomposable. Suppose there is a uniformly bounded sequence $(g_i)_{i=1}^\infty \subseteq C_w(\Lambda; M)$ converging to $f \in C_w(\Lambda; M)$ weak*, then we may choose a sequence $(f_n)_{n=1}^\infty \subseteq C_w(\Lambda; M)$ and an open dense subset $O \subseteq \Lambda$ s.t. $f_n \in \text{co}\{g_i : i \geq n\}$ for all $n$, $f_n \to f$ weak*, and $f_n(\omega) \to f(\omega)$ weak* for all $\omega \in O$. (We may still choose a separable subalgebra of $M$ needed for choosing $\{\varphi_j\}_{j=1}^\infty$ as we can use the algebra generated by $E_M^{sp}(g_i)$ and $E_M^{sp}(f)$, which is countably generated and therefore separable.) A particular consequence of this is that, suppose we already know $f$ is contained in a separable subalgebra $M_0 \subseteq L^\infty(\Omega, \mu) \overline{\otimes} M$, and that $A \subseteq M_0$ is a weak*-dense *-subalgebra, then there exists a sequence $(f_n)_{n=1}^\infty \subseteq A$ and an open dense subset $O \subseteq \Lambda$ s.t. $f_n(\omega) \to f(\omega)$ weak* for all $\omega \in O$. Furthermore, $\|f_n\| \leq \|f\|$ for all $n$ and $f_n \to f$ weak*. ∎

**Further Remark:** We note that certain parts of the theory in this section can actually be developed without assuming $L^\infty(\Omega, \mu)$ is separable. We would first need to extend it to the case where $L^\infty(\Omega, \mu)$ is countably decomposable, where we can still assume $\mu$ is a probability measure. In this case it is currently unknown whether Lemma 5 and consequently Corollary 6 still hold. As such, it would be necessary to define $C_w(\Lambda; M)$ so that for any $f \in C_w(\Lambda; M)$, $k \in H$, $f(\omega)k$ belongs to a separable subspace of $H$. Then we need to prove the range of $\iota: L^\infty(\Omega, \mu) \overline{\otimes} M \to C_w(\Lambda; M)$ actually lies within $C_w(\Lambda; M)$. To do this, we use methods similar to what we did in Proposition 8. We represent $a \in L^\infty(\Omega, \mu) \overline{\otimes} M \cong L^\infty(\Lambda, \mu) \overline{\otimes} M$ as an operator on $L^2(\Lambda, \mu; H)$. Since $\mu$ is a probability measure, $1 \otimes k \in L^2(\Lambda, \mu; H)$. So, $a(1 \otimes k) \in L^2(\Lambda, \mu; H)$ and can be represented as a measurable function $(\Lambda, \mu) \to H$. Then, for any fixed $h \in H$, $\langle h, (a(1 \otimes k))(\omega)\rangle$ is a measurable function in $\omega \in \Lambda$ and $\langle h, \iota(a)(\omega)k\rangle$ is a continuous function. Take any measurable $F \subseteq \Lambda$, then,

$$\int_F \langle h, \iota(a)(\omega)k\rangle \, d\mu = \int_F (Id \otimes |k\rangle\langle h|(a))(\omega) \, d\mu$$

$$= \int_F Id \otimes |k\rangle\langle h|(a) \, d\mu$$

$$= |1\rangle\langle 1_F| \otimes |k\rangle\langle h|(a)$$

$$= \langle 1_F \otimes h, a(1 \otimes k)\rangle$$

$$= \int_F \langle h, (a(1 \otimes k))(\omega)\rangle \, d\mu$$

Since $F \subseteq \Lambda$ is arbitrary, $\langle h, \iota(a)(\omega)k\rangle = \langle h, (a(1 \otimes k))(\omega)\rangle$ $\mu$-a.e. By definition, $(a(1 \otimes k))(\omega)$ lies a.e. in a separable subspace $H_0 \subseteq H$, so if $h \in H_0^\perp$, $\langle h, \iota(a)(\omega)k\rangle = \langle h, (a(1 \otimes k))(\omega)\rangle = 0$ $\mu$-a.e. So $\langle h, \iota(a)(\omega)k\rangle = 0$ on an open dense set. Since $\langle h, \iota(a)(\omega)k\rangle$ is continuous, $\langle h, \iota(a)(\omega)k\rangle = 0$ for all $\omega \in \Lambda$. Since this holds for all $h \in H_0^\perp$, $\iota(a)(\omega)k \in H_0$ for all $\omega \in \Lambda$.

After resolving this technical point, Theorem 7 can be proved just as before. The same goes for Proposition 8, Proposition 10, and Theorem 11.

We can then extend to arbitrary $L^\infty(\Omega,\mu)$ by noting it can be written as $L^\infty(\Omega,\mu) = \prod_{i\in I} L^\infty(\Omega_i,\nu_i)$ where each $L^\infty(\Omega_i,\nu_i)$ is countably decomposable. In fact, $\Omega_i$ can be chosen as pairwise disjoint clopen subsets of $\Lambda$ s.t. $\Lambda' = \bigcup_{i\in I}\Omega_i$ is dense in $\Lambda$ and each $\nu_i$ can be chosen as a Radon probability measure with full support on $\Omega_i$ (see the proof of [Tak79, Theorem III.1.18]). We then define $\nu$, a Radon measure on $\Lambda' = \bigcup_{i\in I}\Omega_i$, by $\nu(B) = \sum_i \nu_i(B\cap \Omega_i)$ for all Borel sets $B$. In this case we need to define strong measurability so that it only requires strong measurability on each $(\Omega_i,\nu_i)$. Then bounded strongly measurable functions still represent operators on $L^2(\Lambda',\nu;H)$, since each vector in $L^2(\Lambda',\nu;H)$ can only be supported on countably many $\Omega_i$. We similarly change the definition of $C_w(\Lambda;M)$. Here, we need to use the fact that any scalar-valued bounded continuous function defined on $\Lambda'$ can be extended to a continuous function on $\Lambda$ [Tak79, Theorem III.1.8], so using Lemma 3 we see that $C_w(\Lambda;M) = C_w(\Lambda';M)$ canonically. The proofs of Theorem 7, Proposition 8, Proposition 10, and Theorem 11 then reduce to proving the results on each $\Omega_i$, which follows from our discussions above.

In fact, if each $L^\infty(\Omega_i,\nu_i)$ is separable (for example, when the abelian von Neumann algebra we started with is $l^\infty(\kappa)\,\overline{\otimes}\,L^\infty([0,1],\lambda)$ for some possibly uncountable cardinal $\kappa$), then proving the results on each $\Omega_i$ actually shows that Theorem 9 also holds in this context. A modified version of Theorem 12 also holds, where instead of asking for $(f_n)_{n=1}^\infty \subseteq C(\Lambda;M)$, we require $(f_n)_{n=1}^\infty \subseteq C_b(\Lambda';M)$. It is easy to see that elements of $C_b(\Lambda';M)$ correspond to elements of $\prod_{i\in I} L^\infty(\Omega_i,\nu_i)\otimes_{min} M$. ∎

As an application of the theory in this section as well as the further remark above, we shall now finish what we started in Section II and prove that the ultrapower defined there, $(M,\tau)^{\mathcal{U}}$, is indeed a tracial von Neumann algebra and $\tau_{\mathcal{U}}$ is a normal faithful tracial state on it.

**Theorem 13:** Let $(\Omega,\mu)$ be a semifinite, decomposable measure space, $\mathcal{U}$ be an ultrafilter on $\mathfrak{B}\big((\Omega,\mu)\big)$, $M$ be a tracial von Neumann algebra with $\tau$ a normal faithful tracial state on $M$. Then $(M,\tau)^{\mathcal{U}}$ is a tracial von Neumann algebra and $\tau_{\mathcal{U}}$ is a normal faithful tracial state on it. Furthermore, if $M$ is a factor, then so is $M^{\mathcal{U}}$.

**Proof:** For the first part of the theorem, we have already seen that $(M,\tau)^{\mathcal{U}}$ is a well-defined $C^*$-algebra and $\tau_{\mathcal{U}}$ is a faithful tracial state. By [AP18, Proposition 2.6.4], it suffices to show $\big((M,\tau)^{\mathcal{U}}\big)_1$ is complete under the 2-norm defined by $\tau_{\mathcal{U}}$.

We first need to represent elements of $L^\infty(\Omega,\mu)\,\overline{\otimes}\,M$ as functions. Following the further remark above, we write $L^\infty(\Omega,\mu) = \prod_{i\in I} L^\infty(\Omega_i,\nu_i)$ where $\Omega_i$ are pairwise disjoint clopen subsets of $\Lambda$, the spectrum of $L^\infty(\Omega,\mu)$, s.t. $\Lambda' = \bigcup_{i\in I}\Omega_i$ is dense in $\Lambda$ and each $\nu_i$ is a Radon probability measure with full support on $\Omega_i$. We then define $\nu$, a Radon measure on $\Lambda' = \bigcup_{i\in I}\Omega_i$, by $\nu(B) = \sum_i \nu_i(B\cap\Omega_i)$ for all Borel sets $B$. We define $C_w(\Lambda;M)$ so that for each $f \in C_w(\Lambda;M)$, $k \in H$, $i \in I$, there exists a separable subspace $H_0 \subseteq H$ s.t. $f(\omega)k \in H_0$ for all $\omega \in \Omega_i$. Under these definitions, $L^\infty(\Omega,\mu)\,\overline{\otimes}\,M \cong C_w(\Lambda;M)$ as in Theorem 7.

We now follow similar lines of proof as [AP18, Proposition 5.4.1]. Let $(x_n)_{n=1}^\infty \subseteq (M,\tau)^{\mathcal{U}}$ be a sequence such that for every $n$, $\|x_n\|_\infty < 1$ and $\|x_{n+1} - x_n\|_{2,\tau_{\mathcal{U}}} < 2^{-(n+1)}$. We begin by lifting $x_1$ to an element $\mathbb{x}_1 \in C_w(\Lambda;M)$ with $\|\mathbb{x}_1\|_\infty \leq 1$. We then begin constructing lifts of $x_n$ inductively. Suppose we already have lifts $\mathbb{x}_n$ of $x_n$ for all $n < N$ s.t. $\|\mathbb{x}_n\|_\infty \leq 1$. Let $\mathbb{y}_N$ be an arbitrary lift of $x_N$ with $\|\mathbb{y}_N\|_\infty \leq 1$. We see that,

$$\|x_N - x_{N-1}\|_{2,\tau_{\mathcal{U}}}^2 = \omega_{\mathcal{U}}\left(Id \otimes \tau((\mathbb{y}_N^* - \mathbb{x}_{N-1}^*) * (\mathbb{y}_N - \mathbb{x}_{N-1}))\right)$$

$$= \tau((\mathbb{y}_N^* - \mathbb{x}_{N-1}^*) * (\mathbb{y}_N - \mathbb{x}_{N-1})(\omega_{\mathcal{U}}))$$

$$< 2^{-2N}$$

By continuity, there exists a clopen neighborhood $F$ of $\omega_{\mathcal{U}}$ s.t. $\tau((\mathbb{y}_N^* - \mathbb{x}_{N-1}^*) * (\mathbb{y}_N - \mathbb{x}_{N-1})(\omega)) < 2^{-2N}$ for all $\omega \in F$. Let $\mathbb{x}_N = 1_F \mathbb{y}_N + 1_{\Lambda \setminus F} \mathbb{x}_{N-1}$. (Recall that by Theorem 11, multiplication with norm-continuous functions such as $1_F$ and $1_{\Lambda \setminus F}$ is done pointwise.) Since $\|\mathbb{x}_{N-1}\|_\infty \leq 1$, $\|\mathbb{y}_N\|_\infty \leq 1$, we see that $\|\mathbb{x}_N\|_\infty \leq 1$. We note that,

$$(\mathbb{y}_N^* - \mathbb{x}_N^*) * (\mathbb{y}_N - \mathbb{x}_N) = \left(1_{\Lambda \setminus F} \mathbb{y}_N^* - 1_{\Lambda \setminus F} \mathbb{x}_{N-1}^*\right) * \left(1_{\Lambda \setminus F} \mathbb{y}_N - 1_{\Lambda \setminus F} \mathbb{x}_{N-1}\right)$$

$$= 1_{\Lambda \setminus F}(\mathbb{y}_N^* - \mathbb{x}_{N-1}^*) * (\mathbb{y}_N - \mathbb{x}_{N-1})$$

In particular, $\omega_{\mathcal{U}}\left(Id \otimes \tau((\mathbb{y}_N^* - \mathbb{x}_N^*) * (\mathbb{y}_N - \mathbb{x}_N))\right) = 0$, so $\mathbb{x}_N$ is a lift of $x_N$. We also have, for $\omega \in \Lambda \setminus F$, $\mathbb{x}_N(\omega) = \mathbb{x}_{N-1}(\omega)$, so $\|\mathbb{x}_N(\omega) - \mathbb{x}_{N-1}(\omega)\|_{2,\tau} < 2^{-N}$. For $\omega \in F$, we observe that,

$$\|\mathbb{x}_N(\omega) - \mathbb{x}_{N-1}(\omega)\|_{2,\tau}^2 = \tau((\mathbb{x}_N^* - \mathbb{x}_{N-1}^*)(\omega)(\mathbb{x}_N - \mathbb{x}_{N-1})(\omega))$$

By Proposition 8, there is an open dense subset $O' \subseteq F$ s.t.,

$$\tau((\mathbb{x}_N^* - \mathbb{x}_{N-1}^*)(\omega)(\mathbb{x}_N - \mathbb{x}_{N-1})(\omega)) = \tau((\mathbb{x}_N^* - \mathbb{x}_{N-1}^*) * (\mathbb{x}_N - \mathbb{x}_{N-1})(\omega))$$

$$= \tau(1_F(\omega)(\mathbb{x}_N^* - \mathbb{x}_{N-1}^*) * (\mathbb{x}_N - \mathbb{x}_{N-1})(\omega))$$

$$= \tau(1_F(\omega)(\mathbb{y}_N^* - \mathbb{x}_{N-1}^*) * (\mathbb{y}_N - \mathbb{x}_{N-1})(\omega))$$

$$= \tau((\mathbb{y}_N^* - \mathbb{x}_{N-1}^*) * (\mathbb{y}_N - \mathbb{x}_{N-1})(\omega))$$

$$< 2^{-2N}$$

For all $\omega \in O'$. Hence, there is an open dense subset $O_N \subseteq \Lambda$, namely $O_N = O' \cup (\Lambda \setminus F)$, s.t. $\|\mathbb{x}_N(\omega) - \mathbb{x}_{N-1}(\omega)\|_{2,\tau} < 2^{-N}$ for all $\omega \in O_N$.

Now, having constructed all the lifts $\mathbb{x}_n$, we see when $\omega$ is in the intersection of all $O_n$, $(\mathbb{x}_n(\omega))_{n=1}^\infty \subseteq (M)_1$ is a Cauchy sequence in $\|\cdot\|_{2,\tau}$, so it converges in 2-norm (and therefore SOT and weak*) to some $\mathbb{x}(\omega) \in (M)_1$. Recall that open dense subsets of $\Lambda$ (more precisely, their intersections with $\Lambda'$), in particular $O_n$, are co-null, so $\mathbb{x}$ (or rather, its restriction to $\Lambda'$) is defined $\nu$-a.e. We claim that $\mathbb{x} \colon (\Lambda', \nu) \to M$ is strongly measurable (in the sense as defined in the further remark above). For any $\varphi \in M_*$, since $\mathbb{x}_n(\omega) \to \mathbb{x}(\omega)$ a.e. weak*, $\varphi(\mathbb{x}_n(\omega)) \to \varphi(\mathbb{x}(\omega))$ a.e. Since $\mathbb{x}_n \in C_w(\Lambda; M)$, $\varphi(\mathbb{x}(\omega))$ is thus an a.e. limit of a sequence of continuous functions and hence measurable. On the other hand, for each $\Omega_i$, $k \in H$, and $n \in \mathbb{N}$, by definition of $C_w(\Lambda; M)$ there is a separable subspace $H_n \subseteq H$ s.t. $\mathbb{x}_n(\omega)k \in H_n$ for all $\omega \in \Omega_i$. Let $H_0$ be the subspace of $H$ generated by all $H_n$, which is separable. Then as $\mathbb{x}_n(\omega) \to \mathbb{x}(\omega)$ a.e. SOT, we see that $\mathbb{x}(\omega)k \in H_0$ a.e. This proves the claim.

As $\mathbb{x}(\omega) \in (M)_1$ a.e., it represents an element of $\left(L^\infty(\Lambda', \nu) \overline{\otimes} M\right)_1 = \left(L^\infty(\Omega, \mu) \overline{\otimes} M\right)_1$. Let $x = (\mathbb{x}) \in \left((M, \tau)^{\mathcal{U}}\right)_1$. We claim that $x_n \to x$ under $\|\cdot\|_{2,\tau_{\mathcal{U}}}$. To do so, we calculate the norm of $Id \otimes \tau((\mathbb{x}_n - \mathbb{x})^*(\mathbb{x}_n - \mathbb{x}))$. Fix $h \in L^1(\Lambda', \nu) = L^\infty(\Lambda', \nu)_*$. Let the polar decomposition of $h$ be $h = u|h|$, so,

$$h\Big(Id \otimes \tau\big((\mathbb{x}_n - \mathbb{x})^*(\mathbb{x}_n - \mathbb{x})\big)\Big) = |u|h|^{1/2}\rangle\langle|h|^{1/2}| \otimes |\hat{1}\rangle\langle\hat{1}|\big((\mathbb{x}_n - \mathbb{x})^*(\mathbb{x}_n - \mathbb{x})\big)$$

$$= \langle(\mathbb{x}_n - \mathbb{x})(|h|^{1/2} \otimes \hat{1}), (\mathbb{x}_n - \mathbb{x})(u|h|^{1/2} \otimes \hat{1})\rangle$$

$$= \int h(\omega)\langle(\mathbb{x}_n(\omega) - \mathbb{x}(\omega))\hat{1}, (\mathbb{x}_n(\omega) - \mathbb{x}(\omega))\hat{1}\rangle \, dv$$

$$= \int h(\omega)\|\mathbb{x}_n(\omega) - \mathbb{x}(\omega)\|_{2,\tau}^2 \, dv$$

Recall that, on the intersection of all $O_n$, $\|\mathbb{x}_{n+1}(\omega) - \mathbb{x}_n(\omega)\|_{2,\tau} < 2^{-(n+1)}$ for all $n$ and $\mathbb{x}_n(\omega) \to \mathbb{x}(\omega)$ in 2-norm, so $\|\mathbb{x}_n(\omega) - \mathbb{x}(\omega)\|_{2,\tau} \leq 2^{-n}$. Thus, $\|\mathbb{x}_n(\omega) - \mathbb{x}(\omega)\|_{2,\tau}^2 \leq 2^{-2n}$ a.e., so $\left|h\Big(Id \otimes \tau\big((\mathbb{x}_n - \mathbb{x})^*(\mathbb{x}_n - \mathbb{x})\big)\Big)\right| \leq 2^{-2n}\|h\|_1$, i.e., $\left\|Id \otimes \tau\big((\mathbb{x}_n - \mathbb{x})^*(\mathbb{x}_n - \mathbb{x})\big)\right\|_\infty \leq 2^{-2n}$. Therefore, $\|x_n - x\|_{2,\tau_u}^2 = \lim_{\mathcal{U}}\Big(Id \otimes \tau\big((\mathbb{x}_n - \mathbb{x})^*(\mathbb{x}_n - \mathbb{x})\big)\Big) \leq 2^{-2n}$, which concludes the proof.

We now prove the second half of the theorem. Assume $M$ is a factor and let $x \in \mathcal{Z}(M^\mathcal{U})$. Let $x = (\mathbb{x})$ for some $\mathbb{x} \in L^\infty(\Omega, \mu) \, \overline{\otimes} \, M$. By the general Dixmier averaging theorem [Dix81, Part III, Chapter 5, Corollary of Theorem 1], we may choose a sequence of convex combinations of $u\mathbb{x}u^*$ which converges to $E_{\mathcal{Z}(L^\infty(\Omega,\mu)\overline{\otimes}M)}(\mathbb{x}) = Id \otimes \tau(\mathbb{x})$ in norm, where $u \in U(L^\infty(\Omega,\mu) \, \overline{\otimes} \, M)$. Since $(u) \in U(M^\mathcal{U})$ and $(\mathbb{x}) \in \mathcal{Z}(M^\mathcal{U})$, we see that $(u\mathbb{x}u^*) = (\mathbb{x})$. Hence, $(Id \otimes \tau(\mathbb{x})) = (\mathbb{x})$. Now, for $f \in L^\infty(\Omega, \mu) \otimes 1$, regarding it as a continuous function on the spectrum, it is easy to see that $(f) = \lim_{\mathcal{U}} f$. Hence, $(\mathbb{x}) = \big(Id \otimes \tau(\mathbb{x})\big) \in \mathbb{C}$, i.e., $M^\mathcal{U}$ is a factor. ∎

The following corollary can be obtained immediately:

**Corollary 14:** $M \ni m \mapsto (1 \otimes m) \in M^\mathcal{U}$ is trace-preserving normal embedding. We shall call it the *diagonal embedding of $M$ into $M^\mathcal{U}$*.

In most cases, more can be said:

**Corollary 15:** Suppose $(\Omega, \mu)$ admits a decomposition into countably infinite many pairwise disjoint measurable sets, each with positive measure. Then for each free ultrafilter $\omega$ on $\mathbb{N}$, there exists an ultrafilter $\mathcal{U}$ on $\mathfrak{B}\big((\Omega, \mu)\big)$ s.t. there exists a trace-preserving normal embedding $M^\omega \hookrightarrow M^\mathcal{U}$.

**Proof:** This can be proved following the same general idea as the second half of the proof of Theorem III.14. The details are left to the readers. ∎

# Section V: Isomorphisms between Tensor Products with Abelian von Neumann Algebras

We now turn to the issue of isomorphisms between $L^\infty(\Omega,\mu) \overline{\otimes} M$ and $L^\infty(\Omega,\mu) \overline{\otimes} N$ which acts as the identity on $L^\infty(\Omega,\mu) \otimes 1$. We shall restrict our attention to the case where $L^\infty(\Omega,\mu)$ is separable. We may then assume $\mu$ is a probability measure. The induced normal Radon probability measure on $\Lambda$ shall also denoted by $\mu$, which has its support the entirety of $\Lambda$. We shall further assume $M$ is countably decomposable. These assumptions shall be in effect until indicated otherwise in this section.

Now, fix an isomorphism $\alpha: L^\infty(\Omega,\mu) \overline{\otimes} M \to L^\infty(\Omega,\mu) \overline{\otimes} N$ s.t. $\alpha|_{L^\infty(\Omega,\mu) \otimes 1} = Id$. Then it induces an isomorphism $C_w(\Lambda; M) \to C_w(\Lambda; N)$, which we shall also denote by $\alpha$. $\alpha|_{L^\infty(\Omega,\mu) \otimes 1} = Id$ means that $\alpha$ restricts to the identity on $C(\Lambda)$. We first show the following lemma:

**Lemma 1:** $\alpha$, when restricted to $C(\Lambda; M)$, acts pointwise, i.e., for all $f \in C(\Lambda; M)$ and $\omega \in \Lambda$, $\alpha(f)(\omega)$ only depends on $f(\omega)$.

**Proof:** It suffices to show $\alpha(f)(\omega) = 0$ when $f(\omega) = 0$. By continuity, for any $\varepsilon > 0$ there is a clopen neighborhood $F$ of $\omega$ s.t. $\|1_F f\| < \varepsilon$. Hence $\|1_F \alpha(f)\| = \|\alpha(1_F f)\| < \varepsilon$. In particular, $\|\alpha(f)(\omega)\| = \|(1_F \alpha(f))(\omega)\| < \varepsilon$. Since this holds for all $\varepsilon > 0$, $\alpha(f)(\omega) = 0$. ∎

**Definition:** For each $\omega \in \Lambda$, let $\alpha_\omega: M \to N$ be defined by $\alpha_\omega(m) = \alpha(1 \otimes m)(\omega)$. Since the operator system structure on $C_w(\Lambda; N)$ is defined pointwise, $\alpha_\omega$ is a ucp map.

By Lemma 1, when $f \in C(\Lambda; M)$, $\alpha(f)(\omega) = \alpha_\omega(f(\omega))$. The situation with general $f \in C_w(\Lambda; M)$ is more complicated, as the following example shows:

**Example:** Let $M = L(S_\infty(\mathbb{Z}))$. Then for each $n \in \mathbb{Z}$, sending any permutation $\sigma \in S_\infty(\mathbb{Z})$ to $\sigma_n$ defined by $\sigma_n(z) = \sigma(z-n) + n$ is an automorphism of $S_\infty(\mathbb{Z})$, which induces an automorphism of $M$, denoted by $\alpha_n$. We may then define an automorphism $\alpha: l^\infty(\mathbb{Z}) \overline{\otimes} M \to l^\infty(\mathbb{Z}) \overline{\otimes} M$ by,

$$\alpha((a_n)) = (\alpha_n(a_n))$$

Suppose $a \in l^\infty(\mathbb{Z}) \overline{\otimes} M$ is defined by $(a_n) = (\sigma_{-n})$ for some fixed nontrivial $\sigma \in S_\infty(\mathbb{Z})$. Then it is easy to see that $\lim_{n \to \pm\infty} a_n = 0$ weak* in $M$, so $a(\omega) = 0$ for all free ultrafilter $\omega$ (see the remark following Corollary III.6). But $\alpha(a) = 1 \otimes \sigma$, so $\alpha(a)(\omega) = \sigma$, i.e., $\alpha$ no longer acts pointwise in this case.

However, $\alpha$ does act pointwise on an open dense set. Indeed, when $\omega$ is a principal ultrafilter and corresponds to $n \in \mathbb{Z}$, then it is easy to see that $\alpha_\omega = \alpha_n$, and we indeed have $\alpha(a)(\omega) = \alpha_\omega(a(\omega))$ in this case. (Working a little bit harder, one can also see that $\alpha_\omega$ is the trace on $M$ when $\omega$ is a free ultrafilter.) We shall see later that under suitable conditions, similar phenomena can happen.

While $\alpha$ is an isomorphism, one cannot expect so much from $\alpha_\omega$. Being ucp maps, they are linear and preserve adjoints, but they are only "nearly isometric", "nearly surjective", and "nearly multiplicative", in the following sense:

**Proposition 2:** The collection of $\alpha_\omega$'s enjoys the following properties:

1. They are "nearly isometric", i.e., for each $m \in M$, there is an open dense subset $O \subseteq \Lambda$ s.t. $\|\alpha_\omega(m)\| = \|m\|$ for all $\omega \in O$;

2. They are "nearly surjective", i.e., for each $m \in N$, there is an open dense subset $O \subseteq \Lambda$ s.t. $m$ is contained in the weak*-closure of the range of $\alpha_\omega$ for all $\omega \in O$;
3. They are "nearly multiplicative", i.e., for each $m_1, m_2 \in M$, there is an open dense subset $O \subseteq \Lambda$ s.t. $\alpha_\omega(m_1 m_2) = \alpha_\omega(m_1)\alpha_\omega(m_2)$ for all $\omega \in O$.

**Proof:**

1. We observe that,
$$\|\alpha_\omega(m)\| = \sup_{\varphi \in (N_*)_1} |\varphi(\alpha_\omega(m))| = \sup_{\varphi \in (N_*)_1} |\varphi(\alpha(1 \otimes m)(\omega))|$$

Since $\alpha(1 \otimes m) \in C_w(\Lambda; N)$, $|\varphi(\alpha(1 \otimes m)(\omega))|$ is continuous, so $\|\alpha_\omega(m)\|$, being the supremum of continuous functions, is lower semicontinuous. In particular, it is measurable. Now, assume the proposition is false. Recalling that a subset of $\Lambda$ is null iff it is nowhere dense, we see that $\{\omega \in \Lambda: \|\alpha_\omega(m)\| \neq \|m\|\}$ has positive measure. As $\|\alpha_\omega(m)\| \leq \|m\|$ for all $\omega$, we see that there exists $\varepsilon > 0$ s.t. $F = \{\omega \in \Lambda: \|\alpha_\omega(m)\| \leq \|m\| - \varepsilon\}$ has positive measure. By inner regularity, there exists compact $K \subseteq F$ with positive measure. By [Tak79, Corollary III.1.13], $K°$, the interior of $K$, has the same positive measure as $K$. Since $K°$ is clopen, we have $\|1_{K°}\alpha(1 \otimes m)\| \leq \|m\| - \varepsilon$. But $\|1_{K°}\alpha(1 \otimes m)\| = \|\alpha(1_{K°} \otimes m)\| = \|1_{K°} \otimes m\| = \|m\|$, a contradiction!

2. By Theorem IV.12, we may choose a sequence $(g_n)_{n=1}^\infty \subseteq C(\Lambda; M)$ s.t. $\|g_n\| \leq \|m\|$ and $g_n \to \alpha^{-1}(1 \otimes m)$ weak*. Thus, $\alpha(g_n) \to 1 \otimes m$ weak*. Then by the remark following Theorem IV.12, we may choose a sequence $(f_n)_{n=1}^\infty \subseteq C(\Lambda; M)$ and an open dense subset $O \subseteq \Lambda$ s.t. $\alpha_\omega(f_n(\omega)) = \alpha(f_n)(\omega) \to m$ weak* for all $\omega \in O$. In particular, $m$ is contained in the weak*-closure of the range of $\alpha_\omega$.
3. This follows directly from Theorem IV.9. ∎

Items 1 and 3 of the above proposition lead directly to the following corollary:

**Corollary 3:** Let $A \subseteq M$ be a norm-separable $C^*$-subalgebra. Then there exists open dense subset $O \subseteq \Lambda$ s.t. $\alpha_\omega$ is an injective *-homomorphism when restricted to $A$ for all $\omega \in O$. In particular, a norm-separable $C^*$-algebra embeds into $M$ iff it embeds into $N$.

We shall now further restrict our attention to the case where $M$ and $N$ are $II_1$ factors. We shall use $\tau_M$ and $\tau_N$ to denote the traces on $M$ and $N$, respectively. These assumptions shall be in effect until indicated otherwise in this section.

**Lemma 4:** All $\alpha_\omega$ are trace-preserving.

**Proof:** Fix any $m \in M$. By the Dixmier averaging theorem, we may choose countably many unitaries $\{u_n\}_{n=1}^\infty$ s.t. a sequence of convex combinations of $u_n m u_n^*$ converges to $\tau_M(m)$ in norm. Let $A$ be the $C^*$-subalgebra of $M$ generated by $m$ and $\{u_n\}_{n=1}^\infty$. Then $A$ is norm-separable and any tracial state on $A$ must evaluate to $\tau_M(m)$ at $m$. By Corollary 3, there is an open dense subset $O \subseteq \Lambda$ s.t. $\alpha_\omega|_A$ is an injective *-homomorphism for all $\omega \in O$. Then $\tau_N \circ \alpha_\omega|_A$ is a tracial state on $A$, so $\tau_N(\alpha_\omega(m)) = \tau_M(m)$ for all $\omega \in O$. Note that $\tau_N(\alpha_\omega(m)) = \tau_N(\alpha(1 \otimes m)(\omega))$ is continuous in $\omega$ and $O \subseteq \Lambda$ is dense, so $\tau_N(\alpha_\omega(m)) = \tau_M(m)$ for all $\omega \in \Lambda$. ∎

**Remark:** Let $\tau UCP(M, N)$ consists of all trace-preserving ucp maps from $M$ to $N$, equipped with the topology of pointwise weak* convergence. Then we see $\Lambda \ni \omega \mapsto \alpha_\omega \in \tau UCP(M, N)$ is continuous. ∎

**Corollary 5:** For each $m \in M$, there exists open dense $O \subseteq \Lambda$ s.t. $\|\alpha_\omega(m)\|_2 = \|m\|_2$ for all $\omega \in O$. Furthermore, $\|\alpha_\omega(m)\|_2 \leq \|m\|_2$ holds for all $\omega \in \Lambda$.

**Proof:** The first half of the corollary follows directly from item 3 of Proposition 2 and Lemma 4. For the second half, as $\alpha_\omega$ is an ucp map, Stinespring dilation implies that $\alpha_\omega(m)^*\alpha_\omega(m) \leq \alpha_\omega(m^*m)$, whence the result follows from Lemma 4. ∎

**Corollary 6:** All $\alpha_\omega$ are normal.

**Proof:** Recall that for a positive map such as $\alpha_\omega$, to prove it is normal it suffices to show $\alpha_\omega(A_i) \to \alpha_\omega(A)$ SOT for a uniformly bounded increasing net $A_i$ converging strongly to $A$ [Con99, Definition 46.1 & Corollary 46.5]. Recall that for a tracial algebra, SOT on bounded sets coincide with the topology generated by the 2-norm (see the proof of [AP18, Proposition 2.6.4]). The result then follows from Corollary 5. ∎

**Remark:** Applying Corollary 6, we see that methods similar to those used to prove item 2 of Proposition 2 yield a stronger result: for each $m \in N$, there is an open dense subset $O \subseteq \Lambda$ s.t. $m$ is contained in the range of $\alpha_\omega$ for all $\omega \in O$. ∎

The main theorem of this discussion now easily follows:

**Theorem 7:** Let $M_0 \subseteq M$ be a separable subalgebra. Then there exists open dense subset $O \subseteq \Lambda$ s.t. $\alpha_\omega$ is a normal trace-preserving embedding when restricted to $M_0$ for all $\omega \in O$. In particular, a separable von Neumann algebra embeds into $M$ iff it embeds into $N$. Furthermore, if $M$ is separable to begin with, then we may choose $O$ s.t. $\alpha_\omega$ is an isomorphism between $M$ and $N$ for all $\omega \in O$.

**Proof:** For the first part of the theorem, combine Theorem I.8, Corollary 3, Lemma 4, and Corollary 6. For the "furthermore" part, observe that $N$ must be separable and therefore countably generated. Then apply the remark after Corollary 6. ∎

We shall now deliver upon the promised improvement to Lemma 1:

**Theorem 8:** $\alpha$ acts "nearly pointwise", i.e., for each $f \in C_w(\Lambda; M)$, there is an open dense subset $O \subseteq \Lambda$ s.t. $\alpha(f)(\omega) = \alpha_\omega(f(\omega))$.

**Proof:** The proof follows similar lines as the proof of item 2 of Proposition 2. By Theorem IV.12, we may choose a sequence $(g_n)_{n=1}^\infty \subseteq C(\Lambda; M)$ s.t. $g_n(\omega) \to f(\omega)$ weak* a.e., $\|g_n\| \leq \|f\|$, and $g_n \to f$ weak*. Thus, $\alpha(g_n) \to \alpha(f)$ weak*. Then by the remark following Theorem IV.12, we may choose a sequence $(f_n)_{n=1}^\infty \subseteq C(\Lambda; M)$ s.t. $f_n \in \mathrm{co}\{g_i : i \geq n\}$ and $\alpha_\omega(f_n(\omega)) = \alpha(f_n)(\omega) \to \alpha(f)(\omega)$ weak* a.e. Note that since $f_n \in \mathrm{co}\{g_i : i \geq n\}$, we have $f_n(\omega) \to f(\omega)$ weak* a.e., so by normality of $\alpha_\omega$, $\alpha_\omega(f_n(\omega)) \to \alpha_\omega(f(\omega))$ a.e. Therefore, $\alpha(f)(\omega) = \alpha_\omega(f(\omega))$ a.e. ∎

If $M$ is separable, as we have seen in Theorem 7, $\alpha_\omega$ are isomorphisms between $M$ and $N$ on an open dense subset $O \subseteq \Lambda$. If we let $\mathrm{Iso}(M, N)$ consists of isomorphisms from $M$ to $N$, equipped with the topology of pointwise weak* convergence, then in light of the remark after Lemma 4 we have a continuous map $O \ni \omega \mapsto \alpha_\omega \in \mathrm{Iso}(M, N)$. In light of Theorem 8, this map can be considered as providing all the structural information of $\alpha$. Furthermore, we have,

**Proposition 9:** Assuming $M$ is separable, then there is an open dense subset $O \subseteq \Lambda$ s.t. $\alpha_\omega$ are invertible and $\alpha_\omega^{-1} = (\alpha^{-1})_\omega$.

**Proof:** We have already seen that $\alpha_\omega$ are isomorphisms a.e. Now, as $M$ is separable, we may pick a countable weak*-dense subset $S \subseteq M$. For each $m \in S$, then, by Theorem 8, a.e.,

$$\alpha_\omega^{-1}(\alpha_\omega(m)) = m$$

$$= \alpha^{-1}(\alpha(1 \otimes m))(\omega)$$

$$= (\alpha^{-1})_\omega(\alpha(1 \otimes m)(\omega))$$

$$= (\alpha^{-1})_\omega(\alpha_\omega(m))$$

Since $S$ is countable, there is a shared co-null set on which $\alpha_\omega^{-1}(\alpha_\omega(m)) = (\alpha^{-1})_\omega(\alpha_\omega(m))$ for all $m \in S$. By weak*-density, $\alpha_\omega^{-1}(\alpha_\omega(m)) = (\alpha^{-1})_\omega(\alpha_\omega(m))$ holds for all $m \in M$ on this co-null set, whence $\alpha_\omega^{-1} = (\alpha^{-1})_\omega$ a.e. ∎

Hence, by restricting $O$ if necessary, $O \ni \omega \mapsto \alpha_\omega^{-1} \in \mathrm{Iso}(N, M)$ is also continuous. These results can be summarized as follows:

**Theorem 10 (Structure Theorem for Isomorphisms $L^\infty(\Omega, \mu) \overline{\otimes} M \to L^\infty(\Omega, \mu) \overline{\otimes} N$ when $M$ and $N$ are Separable II$_1$ Factors, Part 1):** Let $L^\infty(\Omega, \mu)$ be a separable abelian von Neumann algebra, $\Lambda$ be its spectrum, $M$ and $N$ be separable II$_1$ factors, $\beta: L^\infty(\Omega, \mu) \overline{\otimes} M \to L^\infty(\Omega, \mu) \overline{\otimes} N$ be an isomorphism. Then there exists a unique automorphism $\sigma$ of $L^\infty(\Omega, \mu)$, an open dense subset $O \subseteq \Lambda$, a continuous map $O \ni \omega \mapsto \alpha_\omega \in \mathrm{Iso}(M, N)$ whose pointwise inverse $O \ni \omega \mapsto \alpha_\omega^{-1} \in \mathrm{Iso}(N, M)$ is also continuous (the map is unique up to potentially changing $O$ to another open dense subset but maintaining the map on the intersection of domains), s.t. $(\sigma^{-1} \otimes \mathrm{Id}) \circ \beta$, as a map from $C_w(\Lambda; M)$ to $C_w(\Lambda; N)$, is given by $((\sigma^{-1} \otimes \mathrm{Id}) \circ \beta(f))(\omega) = \alpha_\omega(f(\omega))$ a.e.

**Proof:** Since the center of $L^\infty(\Omega, \mu) \overline{\otimes} M$ (or $L^\infty(\Omega, \mu) \overline{\otimes} N$) is $L^\infty(\Omega, \mu) \otimes 1$, both $\beta$ and $\beta^{-1}$ preserve $L^\infty(\Omega, \mu) \otimes 1$. As such, it is easy to see that $\sigma$ is simply the restriction of $\beta$ to $L^\infty(\Omega, \mu) \otimes 1$. The remainder of the theorem is mostly a restatement of Theorem 7, Theorem 8, and Proposition 9, with a bit more work using density and continuity to establish uniqueness of $\alpha_\omega$. ∎

**Remark:** We may also put the first part of Theorem 7, combined with Theorem 8, in a form similar to Theorem 10 as follows: Let $L^\infty(\Omega, \mu)$ be a separable abelian von Neumann algebra, $\Lambda$ be its spectrum, $M$ and $N$ be II$_1$ factors, $M_0 \subseteq M$ be a separable subalgebra, $\beta: L^\infty(\Omega, \mu) \overline{\otimes} M \to L^\infty(\Omega, \mu) \overline{\otimes} N$ be an isomorphism. Then there exists a unique automorphism $\sigma$ of $L^\infty(\Omega, \mu)$, an open dense subset $O \subseteq \Lambda$, a continuous map $O \ni \omega \mapsto \alpha_\omega \in \mathrm{Hom}(M_0, N)$ (the map is unique up to potentially changing $O$ to another open dense subset but maintaining the map on the intersection of domains), where $\mathrm{Hom}(M_0, N)$ is the space of normal trace-preserving embeddings of $M_0$ into $N$ equipped with the topology of pointwise weak* convergence, s.t. $(\sigma^{-1} \otimes \mathrm{Id}) \circ \beta$, when restricted to $L^\infty(\Omega, \mu) \overline{\otimes} M_0$ and regarded as a map from $C_w(\Lambda; M_0)$ to $C_w(\Lambda; N)$, is given by $((\sigma^{-1} \otimes \mathrm{Id}) \circ \beta(f))(\omega) = \alpha_\omega(f(\omega))$ a.e. ∎

The converse to Theorem 10 is also true, i.e., given a continuous map $O \ni \omega \mapsto \alpha_\omega \in \mathrm{Iso}(M, N)$ with continuous pointwise inverse, we can define an isomorphism $\alpha: L^\infty(\Omega, \mu) \overline{\otimes} M \to L^\infty(\Omega, \mu) \overline{\otimes} N$:

**Theorem 11 (Structure Theorem for Isomorphisms $L^\infty(\Omega, \mu) \overline{\otimes} M \to L^\infty(\Omega, \mu) \overline{\otimes} N$ when $M$ and $N$ are Separable II$_1$ Factors, Part 2):** Let $L^\infty(\Omega, \mu)$ be a separable abelian von Neumann algebra, $\Lambda$ be its spectrum, $M$ and $N$ be separable II$_1$ factors, $O \subseteq \Lambda$ be an open dense subset, $O \ni \omega \mapsto \alpha_\omega \in \mathrm{Iso}(M, N)$ be a continuous map s.t. its pointwise inverse $O \ni \omega \mapsto \alpha_\omega^{-1} \in \mathrm{Iso}(N, M)$ is also continuous. Then there

exists an isomorphism $\alpha: L^\infty(\Omega, \mu) \overline{\otimes} M \to L^\infty(\Omega, \mu) \overline{\otimes} N$, which acts as the identity on $L^\infty(\Omega, \mu) \otimes 1$ and which is defined by, as a map from $C_w(\Lambda; M)$ to $C_w(\Lambda; N)$, $\alpha(f)(\omega) = \alpha_\omega(f(\omega))$ a.e. Furthermore, we have $\alpha_\omega(m) = \alpha(1 \otimes m)(\omega)$ for all $m \in M$ and $\omega \in O$.

**Proof:** The furthermore part follows from continuity. So it suffices to show that $\alpha(f)(\omega) = \alpha_\omega(f(\omega))$ a.e. indeed defines an isomorphism $\alpha: L^\infty(\Omega, \mu) \overline{\otimes} M \to L^\infty(\Omega, \mu) \overline{\otimes} N$. Unambiguity of the definition follows again from continuity. Now, choose a countable dense subset $\{\varphi_i\}_{i=1}^\infty \subseteq N_*$. We first observe that if $f = g \otimes m$ where $g \in L^\infty(\Omega, \mu)$ and $m \in M$, then $\varphi_i\left(\alpha_\omega(f(\omega))\right) = g(\omega)\varphi_i(\alpha_\omega(m))$ is continuous. Recall that any bounded continuous function defined on an open dense subset of $\Lambda$ can be extended to a continuous function on $\Lambda$ [Tak79, Theorem III.1.8], so $\alpha_\omega(f(\omega))$ indeed defines a function in $C_w(\Lambda; N)$. By taking linear combinations and norm-limits, this holds for all $f \in C(\Lambda; M)$. Now, for $f \in C_w(\Lambda; M)$, by Theorem IV.12, there exists a sequence $(f_n)_{n=1}^\infty \subseteq C(\Lambda; M)$ s.t. $f_n(\omega) \to f(\omega)$ weak* a.e. Since $\alpha_\omega$ are isomorphisms and as such normal, $\alpha_\omega(f_n(\omega)) \to \alpha_\omega(f(\omega))$ a.e. Hence, for each $\varphi_i$, $\varphi_i\left(\alpha_\omega(f(\omega))\right)$ is an a.e. limit of continuous functions and therefore measurable. By [Tak79, Proposition III.1.12], $\varphi_i\left(\alpha_\omega(f(\omega))\right)$ coincides with a unique continuous functions a.e. By taking norm limits of $\varphi_i$, we see that these continuous functions indeed define an element $\alpha(f)$ of $C_w(\Lambda; N)$. Since there are only countably many $\varphi_i$, there is a shared co-null set on which $\varphi_i\left(\alpha_\omega(f(\omega))\right) = \varphi_i(\alpha(f)(\omega))$ for all $i$, whence $\alpha(f)(\omega) = \alpha_\omega(f(\omega))$ a.e.

Therefore, we see that $\alpha(f)(\omega) = \alpha_\omega(f(\omega))$ a.e. indeed defines a linear map from $C_w(\Lambda; M)$ to $C_w(\Lambda; N)$. It is easily seen to preserve the adjoint. Multiplicativity follows from the following observation: By repeated uses of the definition of $\alpha$ and Theorem IV.9, we have, a.e.,

$$\alpha(f * g)(\omega) = \alpha_\omega(f * g(\omega))$$
$$= \alpha_\omega(f(\omega)g(\omega))$$
$$= \alpha_\omega(f(\omega))\alpha_\omega(g(\omega))$$
$$= (\alpha(f)(\omega))(\alpha(g)(\omega))$$
$$= (\alpha(f) * \alpha(g))(\omega)$$

Hence, $\alpha$ defines a *-homomorphism. Finally, we see that by assumption $O \ni \omega \mapsto \alpha_\omega^{-1} \in \mathrm{Iso}(N, M)$ is also continuous, so it defines a *-homomorphism $\beta: L^\infty(\Omega, \mu) \overline{\otimes} N \to L^\infty(\Omega, \mu) \overline{\otimes} M$. It is easy to see that $\beta$ is the inverse of $\alpha$, so $\alpha$ is an isomorphism. ∎

**Remark:** The preceding theorem actually holds without assuming $M$ and $N$ are separable. The proof proceeds in similar ways up to showing $\alpha_\omega(f(\omega))$ indeed defines an element of $C_w(\Lambda; N)$ when $f \in C(\Lambda; M)$. Then, for general $f \in C_w(\Lambda; M)$, we still have, by Theorem IV.12, that there exists a sequence $(f_n)_{n=1}^\infty \subseteq C(\Lambda; M)$ s.t. $f_n(\omega) \to f(\omega)$ weak* a.e. and then $\alpha_\omega(f_n(\omega)) \to \alpha_\omega(f(\omega))$ a.e. Now, by Lemma IV.5, as $\alpha_\omega(f_n(\omega))$ defines an element of $C_w(\Lambda; N)$ for each $n$ and there are only countably many $f_n$, we see that $\alpha_\omega(f_n(\omega))$ all lie within a separable subalgebra $N_0 \subseteq N$. Hence, $\alpha_\omega(f(\omega)) \in N_0$ a.e. Thus, regarding $\alpha_\omega(f(\omega))$ as a function from $O$ to $N_0$, the same proof as before shows that $\alpha_\omega(f(\omega))$ coincides a.e. with a function in $C_w(\Lambda; N_0) \subseteq C_w(\Lambda; N)$. The rest of proof can then proceed as before.

However, we have chosen to state this theorem only for separable algebras because part 1 of the structure theorem, Theorem 10, does not hold without separability assumption, as we will now see. ∎

We shall now show that part 1 of the structure theorem does not hold without the assumption that $M$ and $N$ are separable. To do so, we note that while representing elements of tensor products as functions on the spectrum is useful for theoretical purposes, it is hard to utilize for the purpose of constructing counterexamples, as it is impossible to define any point on the spectrum without the axiom of choice. As such, we shall introduce a "measurable" version of pointwise definitions of maps $L^\infty(\Omega,\mu) \overline{\otimes} M \to L^\infty(\Omega,\mu) \overline{\otimes} N$, similar to how we can use strongly measurable functions to define elements of $L^\infty(\Omega,\mu) \overline{\otimes} M$. And then we will use Theorem III.11 to obtain its representation on the spectrum. To develop this theory, we shall revert back to the assumptions that $(\Omega,\mu)$ is an arbitrary semifinite, decomposable measure space without assuming $L^\infty(\Omega,\mu)$ is separable or $\mu$ is a probability measure, and that $M$ is an arbitrary von Neumann algebra without assuming it is a $II_1$ factor or countably decomposable. We shall gradually add back those assumptions as we progress.

We start by showing each element of $L^\infty(\Omega,\mu) \overline{\otimes} M$ can indeed by represented by a bounded strongly measurable function $(\Omega,\mu) \to M$:

**Lemma 12:** Given $a \in L^\infty(\Omega,\mu) \overline{\otimes} M$, then there exists a bounded strongly measurable function $f: (\Omega,\mu) \to M$ s.t. $a$ is the operator defined by $f$.

**Proof:** This is a modified version of [Tak79, Theorem IV.7.17]. We shall reproduce a somewhat more streamlined version of the proof here. Following the notation in Section III, let $\rho: L^\infty(\Omega,\mu) \to M^\infty(\Omega,\mu)$ be a lifting map. Now, let $a \in L^\infty(\Omega,\mu) \overline{\otimes} M$ and let $M$ be represented on a Hilbert space $H$. Similar to Lemma IV.3, we see that $\eta: H \times H \to L^\infty(\Omega,\mu)$ defined by $\eta(h,k) = Id \otimes |k\rangle\langle h|(a)$ is a bounded sesquilinear form. Then $\rho \circ \eta: H \times H \to M^\infty(\Omega,\mu)$ is a sesquilinear form. Since $\eta$ commutes with $M'$, $\rho \circ \eta$ does so as well. Evaluation at each point then shows this corresponds to a bounded function $f: (\Omega,\mu) \to M$ s.t. for all $h,k \in H$, $\Omega \ni \omega \mapsto \langle h, f(\omega)k\rangle$ is measurable.

Now, to show that $f$ is strongly measurable, we first observe that, based on the further remark after Theorem IV.12, it suffices to consider the case where $\mu$ is a probability measure. Fix $k \in H$. Then $a(1 \otimes k) \in L^2(\Omega,\mu) \otimes H$. Hence, there exists a separable subspace $H_0 \subseteq H$ s.t. $a(1 \otimes k) \in L^2(\Omega,\mu) \otimes H_0$. We may then write $a(1 \otimes k) = \sum_{i=1}^\infty g_i \otimes k_i$, where $g_i \in L^2(\Omega,\mu)$, $k_i \in H_0$, and the convergence is in $L^2$. Now, fix $h \in H_0^\perp$. We claim that $\eta(h,k) = 0$. Indeed, for any measurable $F \subseteq \Omega$, we have,

$$1_F(\eta(h,k)) = |1\rangle\langle 1_F| \otimes |k\rangle\langle h|(a)$$

$$= \langle 1_F \otimes h, a(1 \otimes k)\rangle$$

$$= \sum_{i=1}^\infty \langle 1_F, g_i\rangle\langle h, k_i\rangle$$

$$= 0$$

As $h \in H_0^\perp$, $k_i \in H_0$. As $F \subseteq \Omega$ is arbitrary, we get $\eta(h,k) = 0$, so $\rho \circ \eta(h,k) = 0$. Because this holds for all $h \in H_0^\perp$, we have $f(\omega)k \in H_0$ for all $\omega \in \Omega$. This proves strong measurability. That $f$ indeed defines the operator $a$ can be proved similarly. We shall leave the precise details to readers. ∎

Despite the result above, in general it may not be a good idea to represent elements of the tensor products in this way and work with them without additional assumptions. Indeed, while one might expect, for example, that two bounded strongly measurable functions would define the same operator iff they coincide a.e., that is actually not true, as the following example demonstrates:

**Example:** Consider $f: ([0, 1], \lambda) \to l^\infty([0, 1])$ defined by $f(t) = E_{\{t\}}$, where $E_X$ is the projection in $l^\infty([0, 1])$ onto $X \subseteq [0, 1]$, and where $l^\infty([0, 1])$ is represented canonically on $H = l^2([0, 1])$. Readers will easily verify that $f$ is indeed strongly measurable. While $f$ is nowhere zero, the operator it defines actually is zero. Indeed, for any $h \in L^2(\Omega, \mu; H)$, we may take its range to lie within a separable subspace of $l^2([0, 1])$, which will then be contained in $l^2(K)$ for a countable subset $K \subseteq [0, 1]$. Hence the support of $f(t)h(t)$ is contained in $K$, which is countable and therefore null, i.e., $f(t)h(t)$ is the zero vector in $L^2(\Omega, \mu; H)$. As this holds for all $h \in L^2(\Omega, \mu; H)$, we see that $f$ defines the zero operator.

In fact, similar to Proposition IV.8 and proved using essentially the same method, one can see that,

**Proposition 13:** Two bounded strongly measurable functions $f, g: (\Omega, \mu) \to M$ defines the same operator iff for all $\varphi \in M_*$, $\varphi \circ f = \varphi \circ g$ a.e.

We shall now reinstate the assumptions that $L^\infty(\Omega, \mu)$ is separable, $\mu$ is a probability measure, and that $M$ is countably decomposable. Furthermore, we shall assume the representation of $M$ is a GNS representation associated with a faithful normal state $\varphi_M$. (These assumptions shall be in effect for the rest of this section.) Then, similar to how Proposition IV.8 can be improved under these conditions, we will also have that the desired two bounded strongly measurable functions defining the same operator iff they coincide a.e. To prove this, we need the analogue of Lemma IV.5,

**Lemma 14:** Given a bounded strongly measurable function $f: (\Omega, \mu) \to M$, its range lies a.e. within a separable subalgebra of $M$.

**Proof:** The main observation needed is that, similar to the case of tracial algebras, on bounded sets SOT and the topology generated by the 2-norm associated with $\varphi_M$ coincide. This is because $\hat{1}$ is separable for $M$ as $\varphi_M$ is faithful, so $\hat{1}$ is cyclic for $M'$. Then the proof proceeds as in the proof of [AP18, Proposition 2.6.4]. Then as $f$ is strongly measurable, $f(\omega)\hat{1}$ lies a.e. within a separable subspace of $H$. As $H$ is a metric space, this implies the range of $f$, restricted to a co-null set, is itself separable under the 2-norm, whence under SOT. As $M$ is countably decomposable, the result follows from Theorem I.8. ∎

Similar to how Theorem IV.9 follows from Lemma IV.5 and Proposition IV.8, we then also have,

**Theorem 15:** Two bounded strongly measurable functions $f, g: (\Omega, \mu) \to M$ defines the same operator iff $f = g$ a.e.

We record some consequences of Lemma 14 and Theorem 15:

**Proposition 16:** Let $f: (\Omega, \mu) \to M$ be a bounded strongly measurable function,

1. The pointwise adjoint $f^*$ of $f$ is also strongly measurable. Furthermore, the operators defined by $f$ and $f^*$ are adjoints of each other;
2. $f$ defines a positive operator iff $f(\omega) \geq 0$ a.e.;
3. $\Omega \ni \omega \mapsto \|f(\omega)\|$ is measurable;
4. The norm of the operator defined by $f$ equals the essential supremum norm of $f$.

**Proof:**

1. Since $f(\omega)$ lies within a separable subalgebra of $M$ a.e., $f^*(\omega)$ does so as well. Strong measurability of $f^*$ then easily follows as in the proof of Corollary IV.6. That $f$ and $f^*$ define adjoint operators is easy to verify.
2. The backwards direction is easy to verify. For the forward direction, let $f$ define $a \in \left(L^\infty(\Omega,\mu) \overline{\otimes} M\right)_+$. We note that the sesquilinear form $\eta\colon H \times H \to L^\infty(\Omega,\mu)$ defined in the proof of Lemma 12 is positive. Indeed, for any positive $\varphi \in L^\infty(\Omega,\mu)_*$, $\varphi(\eta(h,h)) = \varphi \otimes |h\rangle\langle h|(a) \geq 0$. Hence, $\rho \circ \eta\colon H \times H \to M^\infty(\Omega,\mu)$ is also positive. The order structure on $M^\infty(\Omega,\mu)$ is defined pointwise, so the bounded strongly measurable function $f'$ associated to $a$ that was constructed in the proof of Lemma 12 takes positive values everywhere. Theorem 15 implies $f = f'$ a.e., so $f(\omega) \geq 0$ a.e.
3. By altering $f$ on a co-null set if necessary, we may assume $f(\omega)$ all lie within a separable subalgebra $M_0 \subseteq M$. Choose a countable dense subset $S$ of the unit ball of $(M_0)_*$. Then $\|f(\omega)\| = \sup_{\varphi \in S}|\varphi(f(\omega))|$ is the supremum of countably many measurable functions and therefore measurable.
4. Let $a$ be the operator defined by $f$. Then by item 1, $a^*a$ is defined by $f^*(\omega)f(\omega)$. Note that as $a^*a \geq 0$, $\|a\|^2 = \|a^*a\| = \inf\{r \in \mathbb{R}_+ : r - a^*a \geq 0\}$. By item 2, $r - a^*a \geq 0$ iff $r - f^*(\omega)f(\omega) \geq 0$ a.e. iff $r \geq \|f(\omega)\|^2$ a.e., whence the result follows. ∎

**Remark:** None of the items in this proposition holds without assuming $M$ is countably decomposable. For item 1, let $f\colon ([0,1],\lambda) \to \mathbb{B}(l^2([0,1]))$ be defined by $f(t) = E_{0t}$, where $E_{st}$ is the partial isometry with left support $E_{\{s\}}$ and right support $E_{\{t\}}$ for $s,t \in [0,1]$. Then $f$ is strongly measurable but $f^*(t) = E_{t0}$ is not. For item 2, let $f\colon ([0,1],\lambda) \to l^\infty([0,1])$ be defined by $f(t) = -E_{\{t\}}$. Then the operator defined by $f$ is zero, which is positive but $f(t) < 0$ for all $t$. This is also a counterexample for item 4 since the essential supremum norm of $f$ is 1 but the operator it defines is zero. Finally, for item 3, fix a non-Lebesgue-measurable set $A \subseteq [0,1]$ and let $f\colon ([0,1],\lambda) \to l^\infty([0,1])$ be defined by $f(t) = 1_{t \in A}E_{\{t\}}$. Then $f$ is strongly measurable but $\|f(t)\| = 1_{t \in A}$ is not. ∎

Armed with these results, it is now possible to prove an analogue of Theorem IV.12 and then an analogue of Theorem 11. We shall record these results here and left the details of their proofs to readers.

**Theorem 17:** Let $f\colon (\Omega,\mu) \to M$ be a bounded strongly measurable function, then there exists a sequence $(f_n)_{n=1}^\infty$ of bounded strongly measurable functions, representing elements of $L^\infty(\Omega,\mu) \otimes_{min} M$, s.t. $f_n(\omega) \to f(\omega)$ weak* a.e. The sequence may be chosen so that $\|f_n\| \leq \|f\|$ for all $n$ and $f_n \to f$ weak* as elements of $L^\infty(\Omega,\mu) \overline{\otimes} M$.

**Theorem 18:** Assume $M$ and $N$ are $II_1$ factors. Let $\Omega \ni \omega \mapsto \alpha_\omega \in \tau\text{UCP}(M,N)$ be a *strongly measurable map*, in the sense that $\Omega \ni \omega \mapsto \alpha_\omega(m)$ is strongly measurable for all $m \in M$. Then there exists a normal trace-preserving ucp map $\alpha\colon L^\infty(\Omega,\mu) \overline{\otimes} M \to L^\infty(\Omega,\mu) \overline{\otimes} N$, which acts as the identity on $L^\infty(\Omega,\mu) \otimes 1$ and which is defined by, mapping from bounded strongly measurable functions $(\Omega,\mu) \to M$ to bounded strongly measurable functions $(\Omega,\mu) \to N$, $\alpha(f)(\omega) = \alpha_\omega(f(\omega))$. Furthermore, $\alpha$ is a *-homomorphism iff, for any fixed $m \in M$, $m$ is in the multiplicative domain of $\alpha_\omega$ for a.e.-$\omega$.

**Remark:** The converse of Theorem 18 is also true. Namely, given a normal trace-preserving ucp map $\alpha\colon L^\infty(\Omega,\mu) \overline{\otimes} M \to L^\infty(\Omega,\mu) \overline{\otimes} N$ which acts as the identity on $L^\infty(\Omega,\mu) \otimes 1$, one can define a strongly measurable $\Omega \ni \omega \mapsto \alpha_\omega \in \tau\text{UCP}(M,N)$ so that $\alpha(f)(\omega) = \alpha_\omega(f(\omega))$. This is because, after we fix a lifting map $\rho\colon L^\infty(\Omega,\mu) \to M^\infty(\Omega,\mu)$, the construction in Lemma 12 gives a canonical bounded strongly measurable function $\tilde{a}\colon (\Omega, \mu) \to N$ for each $a \in L^\infty(\Omega, \mu) \overline{\otimes} N$. We can then just define $\alpha_\omega(m) = \widetilde{\alpha(1 \otimes m)}(\omega)$. Linearity of $\alpha_\omega$ follows from linearity of $\rho$. Since $\rho$, being a *-homomorphism, acts as the identity on $\mathbb{C}$, we see that $\alpha_\omega$ is unital. Complete positivity of $\alpha_\omega$ follows from similar arguments as in the proof of item 2 of Proposition 16. To show they preserve the trace, we observe that, for any measurable $F \subseteq \Omega$,

$$1_F\left(Id \otimes \tau_N(\alpha(1 \otimes m))\right) = |1_F\rangle\langle 1| \otimes |\hat{1}\rangle\langle \hat{1}|(\alpha(1 \otimes m))$$

$$= \langle 1 \otimes \hat{1}, \alpha(1 \otimes m) 1_F \otimes \hat{1}\rangle$$

$$= \langle 1 \otimes \hat{1}, \alpha(1_F \otimes 1)\alpha(1 \otimes m) 1 \otimes \hat{1}\rangle$$

$$= \langle 1 \otimes \hat{1}, \alpha(1_F \otimes m) 1 \otimes \hat{1}\rangle$$

$$= \mu \otimes \tau_N(\alpha(1_F \otimes m))$$

$$= \mu \otimes \tau_M(1_F \otimes m)$$

$$= 1_F(\tau_M(m)1)$$

Where the equality on the fourth line uses the fact that, as $\alpha$ acts as the identity on $L^\infty(\Omega, \mu) \otimes 1$, $L^\infty(\Omega, \mu) \otimes 1$ is in the multiplicative domain of $\alpha$ so $\alpha(1_F \otimes 1)\alpha(1 \otimes m) = \alpha(1_F \otimes m)$. Hence, $Id \otimes \tau_N(\alpha(1 \otimes m))$ is the constant function $\tau_M(m)1$. Again, as $\rho$ acts as the identity on $\mathbb{C}$, we see that $\tau_N \circ \widetilde{\alpha(1 \otimes m)}$ is the constant function $\tau_M(m)1$, i.e., all $\alpha_\omega$ preserve the trace.

Another way to see $\alpha_\omega$'s enjoy such properties is to observe that they are actually a subset of the corresponding maps on the spectrum. This follows from the following observation: Adopting the notation in Section III, $\rho$ followed by pointwise evaluation at $\omega \in \Omega$ is a character on $L^\infty(\Omega, \mu)$, denoted by $\phi_\rho(\omega)$. Then it is easy to see that $\alpha_\omega = \alpha_{\phi_\rho(\omega)}$, where the latter map is as defined after Lemma 1. This also means that the method used here cannot produce any map between $M$ and $N$ not obtained by working on the spectrum.

However, while the above implies a strong connection to the corresponding maps on the spectrum, as opposed to the spectrum case where we can define $\alpha_\omega$ in a unique and unambiguous way, the same is not true here. In fact, even if $\alpha_\omega \neq \beta_\omega$ everywhere, the induced maps $\alpha$ and $\beta$ could still be the same. We can even obtain a quite "wild" counterexample assuming the continuum hypothesis, with $M$ being non-separable, and a strongly measurable map $([0, 1], \lambda) \ni t \mapsto \alpha_t \in \tau\mathrm{UCP}(M, M)$ s.t. the multiplicative domain of $\alpha_t$ is always a separable subalgebra of $M$, but the induced map $\alpha\colon L^\infty([0, 1], \lambda) \overline{\otimes} M \to L^\infty([0, 1], \lambda) \overline{\otimes} M$ is the identity map. Indeed, since we are assuming the continuum hypothesis, let $\omega_1$ be the first uncountable ordinal, and we may fix a bijection $k\colon [0, 1] \to \omega_1$. Let $P$ be any separable $\mathrm{II}_1$ factor, $M = P^{\overline{\otimes}\omega_1}$, i.e., the tracial tensor product of continuum many copies of $P$, indexed by all countable (including finite) ordinals. We then define, for any $t \in [0, 1]$, $\alpha_t$ to be the conditional expectation onto $P^{\overline{\otimes}k(t)}$. Then as $k(t)$ is countable, the multiplicative domain of $\alpha_t$ is indeed always separable. Now, for any fixed $x \in P^{\overline{\otimes}\omega_1}$, it necessarily belongs to $P^{\overline{\otimes}K}$ for some countable subset $K \subseteq \omega_1$. (This can be seen by representing $x$ in terms of an $L^2$ basis consisting of simple tensors.) As $K$ is a countable set of countable ordinals, it has a countable ordinal $\omega$ as its supremum, namely $\omega = \bigcup K$. Hence for all $t \in [0, 1]$, as long as $k(t) > \omega$, $x$ is in the multiplicative domain of $\alpha_t$ and in fact $\alpha_t(x) = x$. Since there are only countably many $t$ with $k(t) \leq \omega$ and countable sets are null in $([0, 1], \lambda)$, it is now easy to deduce that the map $([0,1], \lambda) \ni t \mapsto \alpha_t \in \tau\mathrm{UCP}(M, M)$ is indeed strongly measurable and that it defines the identity map on $L^\infty([0,1], \lambda) \overline{\otimes} M$.

Now, one might argue that this example is quite "unnatural" and relies on assuming the continuum hypothesis. We shall later describe another example that does not do so and is somewhat more natural. ∎

We shall now construct the promised counterexample to Theorem 10 without assuming the separability of the $\mathrm{II}_1$ factors involved. First, fix a sequence of automorphisms of the hyperfinite $\mathrm{II}_1$ factor $R$, $(\alpha_n)_{n=1}^\infty$, s.t. $\alpha_n \to \tau_R$ and $\alpha_n^{-1} \to \tau_R$ pointwise weak*. (This is possible. For example, $R \cong L(S_\infty(\mathbb{Z}))$, so the automorphisms in the example preceding Proposition 2 satisfy the requirements.) Write $\alpha_\infty = \tau_R$.

Now, for each $r \in [0,1]$ and $n \in \mathbb{N}$, let,

$$L_n^r = \left[r - \frac{r}{2^{n-1}}, r - \frac{r}{2^n}\right)$$

$$R_n^r = \left(r + \frac{1-r}{2^n}, r + \frac{1-r}{2^{n-1}}\right]$$

Note that these sets form a partition of $[0,1] \setminus \{r\}$. Now, for each $t \in [0,1]$, let,

$$\alpha_t^r = \begin{cases} \alpha_n, & \text{if } t \in L_n^r \text{ or } t \in R_n^r \\ \alpha_\infty, & \text{if } t = r \end{cases}$$

$$\beta_t^r = \begin{cases} \alpha_n^{-1}, & \text{if } t \in L_n^r \text{ or } t \in R_n^r \\ \alpha_\infty, & \text{if } t = r \end{cases}$$

Now, let $M = R^{\overline{\otimes}[0,1]}$, i.e., the tracial tensor product of continuum many copies of $R$, indexed by elements of $[0,1]$. Define, for each $t \in [0,1]$,

$$\alpha_t : M \to M, \alpha_t = \bigotimes_{r \in [0,1]} \alpha_t^r$$

$$\beta_t : M \to M, \beta_t = \bigotimes_{r \in [0,1]} \beta_t^r$$

Readers can verify that these are indeed well-defined trace-preserving ucp maps. Note that, for any fixed $m \in M$, there is a countable set $K \subseteq [0,1]$ s.t. $m \in R^{\overline{\otimes}K} \subseteq R^{\overline{\otimes}[0,1]}$. Furthermore, $\alpha_t(R^{\overline{\otimes}K}) \subseteq R^{\overline{\otimes}K}$ and $\beta_t(R^{\overline{\otimes}K}) \subseteq R^{\overline{\otimes}K}$ and for fixed $r \in [0,1]$ both $\alpha_t^r$ and $\beta_t^r$ only take countably many values as $t$ ranges over $[0,1]$. Using these facts, readers can verify that $[0,1] \ni t \mapsto \alpha_t \in \tau\mathrm{UCP}(M, M)$ is a strongly measurable map by first checking this on simple tensors and then approximating weak* an arbitrary element of $M$ by a sequence of linear combinations of simple tensors. We also observe that, if $m \in R^{\overline{\otimes}K}$, then $m$ is in the multiplicative domain of $\alpha_t$ for all $t \in [0,1] \setminus K$, so for any fixed $m \in M$, $m$ is in the multiplicative domain of $\alpha_t$ for a.e.-$t$. It is clear that these results also hold for $\beta_t$. Hence, by Theorem 18, they define normal *-homomorphisms $\alpha, \beta : L^\infty([0,1], \lambda) \overline{\otimes} M \to L^\infty([0,1], \lambda) \overline{\otimes} M$. They are also inverses to each other. Indeed, for any $a \in L^\infty([0,1], \lambda) \overline{\otimes} M$, it is contained in $L^\infty([0,1], \lambda) \overline{\otimes} M_0$ for some separable subalgebra $M_0 \subseteq M$. But as $M_0$ is countably generated, it is therefore contained in $R^{\overline{\otimes}K}$ for some countable $K \subseteq [0,1]$. Representing $a$ as a bounded strongly measurable function $f : ([0,1], \lambda) \to R^{\overline{\otimes}K}$, we see that whenever $t \in [0,1] \setminus K$, $\beta_t(\alpha_t(f(t))) = f(t)$, so $\beta_t(\alpha_t(f(t))) = f(t)$ a.e. and $\beta(\alpha(a)) = a$. Similarly, $\alpha(\beta(a)) = a$. Hence $\alpha$ is an automorphism of $L^\infty([0,1], \lambda) \overline{\otimes} M$ acting as the identity on $L^\infty([0,1], \lambda) \otimes 1$.

**Remark:** Continuing the discussion started in the remark following Theorem 18, we see that the above example also provides a relatively more "natural" (and independent of the continuum hypothesis) case of two strongly measurable maps $\Omega \to \tau\mathrm{UCP}(M, N)$ differing everywhere yet nevertheless inducing the same map $L^\infty(\Omega, \mu) \overline{\otimes} M \to L^\infty(\Omega, \mu) \overline{\otimes} N$. Indeed, if we change the definition of $\alpha_\infty$ in the above example to any other normal trace-preserving ucp maps from $R$ to $R$, the resulting $\alpha_t$'s would differ from the ones in the above example everywhere. But for any $m \in R^{\overline{\otimes}K}$ where $K \subseteq [0, 1]$ is countable, $\alpha_t(m)$ is unchanged so long as $t \in [0, 1] \setminus K$, so $\alpha_t(m)$ is unchanged a.e., which means the induced maps $L^\infty([0, 1], \lambda) \overline{\otimes} M \to L^\infty([0, 1], \lambda) \overline{\otimes} N$ are exactly the same. ∎

Let $\Lambda$ be the spectrum of $L^\infty([0, 1], \lambda)$. At this point it is already relatively straightforward using Theorem III.11 to show $\alpha_\omega$ is never an isomorphism for all $\omega \in \Lambda$. In fact, they are not even *-homomorphisms. For any $x \in R, r \in [0, 1]$, we shall write $x^{(r)}$ for the simple tensor in $R^{\overline{\otimes}[0,1]}$ that takes value $x$ at $r \in [0, 1]$ and 1 at all other points. Then $\lim_{t \to r} \alpha_t(x^{(r)}) = \lim_{t \to r} \alpha_t^r(x)^{(r)} = \lim_{n \to \infty} \alpha_n(x)^{(r)} = \tau_R(x)$. Hence, for any $\omega \in \Omega_{\mathcal{D}_r}$, $\alpha_\omega(x^{(r)}) = \omega \otimes \mathrm{Id}\left(\alpha(1 \otimes x^{(r)})\right) = \lim_{t \to r} \alpha_t(x^{(r)}) = \tau_R(x)$, whence $\alpha_\omega$ is not a *-homomorphism. Since $\Omega_{\mathcal{D}_r}, r \in [0, 1]$ form a partition of $\Lambda$ by Proposition III.3, we see that $\alpha_\omega$ is never a *-homomorphism for all $\omega \in \Lambda$.

Of course, using Theorem III.11 we can go further and compute $\alpha_\omega$. Because $\alpha_t^r$ has countably many points at which it jumps, we will need to distinguish between left limits and right limits. To clarify notations, we define, for each $n \in \mathbb{N}, r \in [0, 1]$,

$$L_n^{r,L} = \left(r - \frac{r}{2^{n-1}}, r - \frac{r}{2^n}\right]$$

$$R_n^{r,L} = \left(r + \frac{1-r}{2^n}, r + \frac{1-r}{2^{n-1}}\right]$$

$$L_n^{r,R} = \left[r - \frac{r}{2^{n-1}}, r - \frac{r}{2^n}\right)$$

$$R_n^{r,R} = \left[r + \frac{1-r}{2^n}, r + \frac{1-r}{2^{n-1}}\right)$$

$$\alpha_t^{r,L} = \begin{cases} \alpha_n, & \text{if } t \in L_n^{r,L} \text{ or } t \in R_n^{r,L} \\ \alpha_\infty, & \text{if } t = r \end{cases}, 0 < t \leq 1$$

$$\alpha_t^{r,R} = \begin{cases} \alpha_n, & \text{if } t \in L_n^{r,R} \text{ or } t \in R_n^{r,R} \\ \alpha_\infty, & \text{if } t = r \end{cases}, 0 \leq t < 1$$

$$\alpha_t^L : M \to M, \alpha_t^L = \bigotimes_{r \in [0,1]} \alpha_t^{r,L}, 0 < t \leq 1$$

$$\alpha_t^R : M \to M, \alpha_t^R = \bigotimes_{r \in [0,1]} \alpha_t^{r,R}, 0 \leq t < 1$$

We can then calculate that, for any $\omega \in \Omega_{\mathcal{L}_t}, 0 < t \leq 1$, we have $\alpha_\omega = \alpha_t^L$; and for any $\omega \in \Omega_{\mathcal{R}_t}, 0 \leq t < 1$, we have $\alpha_\omega = \alpha_t^R$. This can be done by checking it on simple tensors, the precise details of which are left to readers. We then see that $\alpha_\omega$ is never injective, multiplicative, or surjective. In fact it is always the composition of the conditional expectation onto a proper subalgebra (namely $R^{\overline{\otimes}[0,1]\setminus\{t\}}$ when $\omega \in \Omega_{\mathcal{D}_t}$) followed by an automorphism of the said subalgebra (namely $\bigotimes_{r\in[0,1]\setminus\{t\}} \alpha_t^{r,L}$ when $\omega \in \Omega_{\mathcal{L}_t}$ and $\bigotimes_{r\in[0,1]\setminus\{t\}} \alpha_t^{r,R}$ when $\omega \in \Omega_{\mathcal{R}_t}$).

**Example:** Using these notations, we can also construct a counterexample to Theorem IV.9 without the assumption of countable decomposability. Indeed, let $M = l^\infty([0,1]) \overline{\otimes} R$. Fix a traceless unitary $u \in R$. As a strongly measurable function $([0,1], \lambda) \to M$ and regarding elements of $M$ as functions from $[0,1]$ to $R$, we define $a \in L^\infty([0,1], \lambda) \overline{\otimes} M$ by $(a(t))(r) = \alpha_t^r(u)$. Then $a^*$ is given by $(a^*(t))(r) = \alpha_t^r(u^*)$. We left it to readers to verify these are indeed strongly measurable functions. Readers may then verify, applying Proposition 13, that $aa^* = 1$. But using similar methods as in the example above, we can prove that, regarding $a$ and $a^*$ as continuous functions on the spectrum, $a(\omega)a^*(\omega) \neq 1$ for all $\omega \in \Lambda$.

This counterexample to Theorem 10 in the non-separable case raises the question of what exactly can be known when $M$ and $N$ are $II_1$ factors and $L^\infty(\Omega, \mu) \overline{\otimes} M \cong L^\infty(\Omega, \mu) \overline{\otimes} N$, without assuming $M$ and $N$ are separable. Now, of course if $(\Omega, \mu)$ has atoms, we easily see that $M \cong N$. As any atomless separable abelian von Neumann algebra is isomorphic to $L^\infty([0,1], \lambda)$ [AP18, Theorem 3.2.4], the only interesting case of this question is when $L^\infty([0,1], \lambda) \overline{\otimes} M \cong L^\infty([0,1], \lambda) \overline{\otimes} N$. While in the example above, we still have $M \cong N$, this can not be simply obtained by looking at $\alpha_\omega$'s. This presents an obstacle to answering the question whether $L^\infty([0,1], \lambda) \overline{\otimes} M \cong L^\infty([0,1], \lambda) \overline{\otimes} N$ implies $M \cong N$ for $II_1$ factors $M$ and $N$, which is a special case of a question raised by Ozawa [Oza17].

We shall prove the following theorem, which will answer another special case of Ozawa's question as well as provide an extra condition on $M$ and $N$ when they are $II_1$ factors and $L^\infty(\Omega, \mu) \overline{\otimes} M \cong L^\infty(\Omega, \mu) \overline{\otimes} N$.

**Theorem 19:** Let $P, Q$ be von Neumann algebras. Suppose $L^\infty(\Omega, \mu) \overline{\otimes} P \cong L^\infty(\Omega, \mu) \overline{\otimes} Q$ and that $P$ can be generated by $\kappa$ many elements where $\kappa$ is an infinite cardinal. Then $Q$ can be generated by $\kappa$ many elements as well.

**Proof:** As $L^\infty(\Omega, \mu)$ is separable and therefore countably generated and $P$ can be generated by $\kappa$ many elements, we see that $L^\infty(\Omega, \mu) \overline{\otimes} P$ and therefore $L^\infty(\Omega, \mu) \overline{\otimes} Q$ can be generated by $\kappa$ many elements. Hence, there exists a generating set $S \subseteq L^\infty(\Omega, \mu) \overline{\otimes} Q$ of cardinality $\kappa$. For each $s \in S$, by Corollary I.4, $E_Q^{alg}(s)$ is countably generated, so we may choose a countable generating set $K_s \subseteq E_Q^{alg}(s)$. We claim that $\bigcup_{s \in S} K_s$ generates $Q$, whence the result follows.

To prove the claim, let $K$ be the algebra generated by $\bigcup_{s \in S} K_s$, which then contains $E_Q^{alg}(s)$ for all $s \in S$. Thus, $S \subseteq L^\infty(\Omega, \mu) \overline{\otimes} K$. But $S$ generates $L^\infty(\Omega, \mu) \overline{\otimes} Q$, so $K = Q$. ∎

**Corollary 20:** Suppose $L^\infty(\Omega, \mu) \overline{\otimes} \mathbb{B}(l^2(I)) \cong L^\infty(\Omega, \mu) \overline{\otimes} \mathbb{B}(l^2(J))$. Then $|I| = |J|$.

**Proof:** When either $I$ or $J$ is finite or countably infinite, the result is well-known, see [Con99, Proposition 51.1]. Now, assume both are uncountable. Then $\mathbb{B}(l^2(I))$ is generated by $|I|$ many elements, as there are $|I \times I| = |I|$ many standard matrix elements and they generate $\mathbb{B}(l^2(I))$. Hence, by Theorem 19, $\mathbb{B}(l^2(J))$ can also be generated by $|I|$ many elements. Then there is an SOT-dense subset $S \subseteq \mathbb{B}(l^2(J))$ of cardinality at most $|I|$. Hence $l^2(J) = \mathbb{B}(l^2(J))h = \overline{Sh}$ for any nonzero $h \in l^2(J)$. But $S$ is of cardinality at most $|I|$, so $l^2(J) = \overline{Sh}$ is of dimension at most $|I|$. Hence, $|J| \leq |I|$. By symmetry $|I| \leq |J|$ as well, so $|I| = |J|$. ∎

**Remark:** A different proof of a more general version of this result can be found in [Tak79, Theorem V.1.27]. ∎

As such, when $M$ and $N$ are $II_1$ factors and $L^\infty([0,1], \lambda) \overline{\otimes} M \cong L^\infty([0,1], \lambda) \overline{\otimes} N$, we have two known restrictions on $M$ and $N$: By Theorem 7, they must share the same collection of separable subalgebras. And by Theorem 19, they must be generated by the same number of elements. Another restriction can be given using the concept of ultrapowers defined in Section II:

**Theorem 21:** Let $M$ and $N$ be $II_1$ factors. If $L^\infty(\Omega, \mu) \overline{\otimes} M \cong L^\infty(\Omega, \mu) \overline{\otimes} N$, then $M^\mathcal{U} \cong N^\mathcal{U}$ for all ultrafilters $\mathcal{U}$ on $\mathfrak{B}((\Omega, \mu))$.

**Proof:** Since the center of $L^\infty(\Omega, \mu) \overline{\otimes} M$ (or $L^\infty(\Omega, \mu) \overline{\otimes} N$) is $L^\infty(\Omega, \mu) \otimes 1$, any isomorphism between $L^\infty(\Omega, \mu) \overline{\otimes} M$ and $L^\infty(\Omega, \mu) \overline{\otimes} N$ preserves $L^\infty(\Omega, \mu) \otimes 1$. So by composing with an automorphism of $L^\infty(\Omega, \mu)$, we may choose an isomorphism $\alpha: L^\infty(\Omega, \mu) \overline{\otimes} M \to L^\infty(\Omega, \mu) \overline{\otimes} N$ that acts as the identity on $L^\infty(\Omega, \mu) \otimes 1$. We claim then that $\alpha(I_{\mathcal{U}, \tau_M}) = I_{\mathcal{U}, \tau_N}$, which implies the result. By symmetry, it suffices to show $\alpha(I_{\mathcal{U}, \tau_M}) \subseteq I_{\mathcal{U}, \tau_N}$. Hence, it suffices to show $Id \otimes \tau_N(\alpha(x)) = Id \otimes \tau_M(x)$ for all $x \in L^\infty(\Omega, \mu) \overline{\otimes} M$. Representing elements of the tensor products as weak* continuous functions on the spectrum, we see that,

$$Id \otimes \tau_N(\alpha(x))(\omega) = \tau_N(\alpha(x)(\omega))$$

$$= \tau_N\left(\alpha_\omega(x(\omega))\right)$$

$$= \tau_M(x(\omega))$$

$$= Id \otimes \tau_M(x)(\omega)$$

For a.e. $\omega$. By continuity, $Id \otimes \tau_N(\alpha(x)) = Id \otimes \tau_M(x)$. ∎

**Remark:** When $(\Omega, \mu)$ has atoms and $\mathcal{U}$ is the principal ultrafilter at an atom, it is easy to see that $M^\mathcal{U} \cong M$, so the above result just means $M \cong N$. When $(\Omega, \mu)$ has no atom, it is unclear what the above result implies. In particular, while formally the result may seem connected to elementary equivalence, it is unclear whether that is actually the case, as elementary equivalence needs an isomorphism between ultrapowers with respect to free ultrafilters on discrete sets, not ultrafilters on $\mathfrak{B}((\Omega, \mu))$ for atomless $(\Omega, \mu)$. ∎

**Closing Remark:** Theorems 10 and 11 can be seen as a more concrete (and representation-independent) version of the uniqueness of direct integral disintegration in the setting of tensor products of separable $II_1$ factors with separable abelian von Neumann algebras (see [Tak79, Theorem IV.8.23]). In this context, one can regard the theory developed in Section IV as a generalization of the theory of direct integral, as applied to tensor products, to non-separable von Neumann algebras. It might be interesting to ask whether this theory can be extended beyond tensor products to provide a more general theory of direct integrals of a "measurable field" (or possibly "continuous field") of non-separable von Neumann algebras. ∎

**Acknowledgement:** The author wishes to thank Adrian Ioana and Gregory Patchell for their help during the writing of this note.